\documentclass[a4paper,10pt]{article}
\usepackage{a4, amsxtra,amsmath,amsfonts,amscd, amssymb, mathrsfs}
\usepackage[all]{xy}
\usepackage{graphicx}
\usepackage{ifpdf}
\newtheorem{thm}{Theorem}[section]
\newtheorem{prop}[thm]{Proposition}
\newtheorem{cor}[thm]{Corollary}
\newtheorem{lem}[thm]{Lemma}
\newtheorem{ex}[thm]{Example}
\newtheorem{Def}[thm]{Definition}
\newtheorem{rem}[thm]{Remark}
\newtheorem{art}[thm]{}

\newcommand{\Hom}{{\rm Hom}}
\newcommand{\codim}{{\rm codim}}

\newcommand{\cyc}{{\rm cyc}}
\newcommand{\Pic}{{\rm Pic}}

\newcommand{\Spec}{{\rm Spec}}
\newcommand{\Spf}{{\rm Spf}}

\newcommand{\ve}{{\varepsilon}}
\newcommand{\id}{{\rm id}}

\newcommand{\Acal}{{\mathscr A}}

\newcommand{\Ccal}{{\mathscr C}}
\newcommand{\Dcal}{{\mathscr D}}
\newcommand{\Ecal}{{\mathscr E}}
\newcommand{\Fcal}{{\mathscr F}}
\newcommand{\Gcal}{{\mathscr G}}

\newcommand{\Lcal}{{\mathscr L}}
\newcommand{\Mcal}{{\mathscr M}}

\newcommand{\Ocal}{{\mathscr O}}

\newcommand{\Ucal}{{\mathscr U}}

\newcommand{\Vcal}{{\mathscr V}}
\newcommand{\Xcal}{{\mathscr X}}

\newcommand{\Ycal}{{\mathscr Y}}
\newcommand{\Ymcal}{{\mathcal Y}}

\newcommand{\qdop}{{\mathbb Q}}
\newcommand{\ndop}{{\mathbb N}}
\newcommand{\rdop}{{\mathbb R}}
\newcommand{\pdop}{{\mathbb P}}
\newcommand{\kdop}{{\mathbb K}}
\newcommand{\ldop}{{\mathbb L}}
\newcommand{\bdop}{{\mathbb B}}
\newcommand{\adop}{{\mathbb A}}
\newcommand{\zdop}{{\mathbb Z}}

\newcommand{\tdop}{{\mathbb T}}

\newcommand{\proof}{\noindent {\bf Proof: \/}}
\newcommand{\qed}{{ \hfill $\square$}}

\newcommand{\Trop}{{\rm Trop}}
\newcommand{\trop}{{\rm trop}}

\newcommand{\ub}{{\mathbf u}}

\newcommand{\mb}{{\mathbf m}}

\newcommand{\rb}{{\mathbf r}}
\newcommand{\Sb}{{\mathbf s}}
\newcommand{\xb}{{\mathbf x}}
\newcommand{\yb}{{\mathbf y}}
\newcommand{\zb}{{\mathbf z}}
\newcommand{\Tor}{{\mathbb G}_m^n}

\newcommand{\Xan}{{X^{\rm an}}}
\newcommand{\Yan}{{Y^{\rm an}}}
\newcommand{\Tan}{{T^{\rm an}}}
\newcommand{\Uan}{{U^{\rm an}}}

\newcommand{\relint}{{\rm relint}}

\newcommand{\kcirc}{{K^\circ}}
\newcommand{\ktilde}{{ \tilde{K}}}

\newcommand{\In}{{\rm in}}
\newcommand{\rec}{{\rm rec}}
\newcommand{\Stab}{{\rm Stab}}
\newcommand{\Wt}{{\rm Wt}}
\newcommand{\Univ}{{\rm Univ}}
\newcommand{\Hilb}{{\rm Hilb}}

\title{A guide to tropicalizations}
\author{Walter Gubler}
\date{\today}
\begin{document}

\maketitle

\begin{abstract}
Tropicalizations form a bridge between algebraic and convex geometry. We generalize basic results from tropical geometry which are well-known for special ground fields to arbitrary non-archimedean valued fields. To achieve this, we develop a theory of toric schemes over valuation rings of rank $1$. As a basic tool, we  use techniques from non-archimedean analysis.

\vspace{2mm}

\end{abstract}

\section{Introduction}

Let us consider a field $K$ endowed with a non-archimedean absolute value $|\phantom{a}|$. We fix coordinates $x_1, \dots, x_n$ on the split multiplicative torus $\Tor$ over $K$. Using logarithmic coordinates $-\log|x_1|, \dots , -\log|x_n|$, any closed subscheme $X$ of $\Tor$ transforms into a finite union $\Trop(X)$ of polyhedra in $\rdop^n$. This process is called {\it tropicalization} and it can be used to transform a problem from algebraic geometry into a corresponding problem in convex geometry which is usually easier. If the toric coordinates are well suited to the problem, it is sometimes possible to use a solution of the convex problem to solve the original algebraic problem. Another strategy is to vary the ambient torus to compensate the loss of information due to the tropicalization process.

Tropicalization originates from a paper of Bergman \cite{Berg} on logarithmic limit sets. The convex structure of the tropical variety $\Trop(X)$  was worked out by Bieri--Groves \cite{BG} with applications to geometric group theory in mind. Sturmfels \cite{Stu} pointed out that $\Trop(X)$ is a subcomplex of the Gr\"obner complex. In fact, the polyhedral complex $\Trop(X)$ has some natural weights satisfying a balancing condition which appears first in Speyer's thesis \cite{Spe}. This relies on the description of the Chow cohomology of a toric variety given by Fulton and Sturmfels \cite{FS}. An intrinsic approach to tropical geometry was proposed by Mikhalkin. The idea is to develop tropical geometry as some sort of algebraic geometry based on the min-plus algebra where every $\Trop(X)$ occurs as a natural object. This approach was used by Mikhalkin to prove celebrated results in enumerative geometry (see \cite{Mik}). These results popularized tropical geometry generating a huge amount of interesting results and applications. For the relation of tropicalizations to idempotent mathematics and Maslov dequantization, we refer to \cite{Lit}.

In tropical literature, one  usually considers tropicalizations under severe restrictions for the ground  field suited to the situation at hand. These restrictions can be subdivided into four groups. First, many papers are written in case of the trivial valuation on $K$. Then $\Trop(X)$ is a finite union of rational cones. This is a setting which occurs very often in algebraic geometry. A second group of people is working under the assumption that $K$ is the field of Puiseux series or a related field. This is often used by people interested in the combinatorial structure and effective computation of $\Trop(X)$. In this case, $K$ has a natural grading,  contains the residue field and the valuation has a canonical section which  makes many arguments easier. A third group is assuming that the valuation is discrete. Most of the valuations occurring in applications to certain fields such as number theory are discrete. This makes it possible to use noetherian models over the valuation ring. Finally, a fourth group of people is working with  algebraically closed ground fields endowed with a non-trivial complete absolute value. This is suitable for using arguments from the theory of rigid analytic spaces. As an excellent source for this case, we refer the reader to the recent paper of Baker--Payne--Rabinoff \cite{BPR}.  

The goal of this paper is to survey basic results about tropicalizations and to generalize them to arbitrary non-archimedean absolute values on $K$. This will make these results accessible to all kind of applications. To handle the difficulties mentioned above, we will use methods from the theory of Berkovich analytic spaces which are very well-suited for this general framework. This is not surprising as even in the original Bieri--Groves paper, the analytification $\Xan$ of $X$ with its Berkovich analytic topology was implicitly used  before Berkovich introduced his new concept in rigid analytic geometry. Most parts of the paper can be read having just a topological understanding of $\Xan$ which is rather elementary. The paper is not meant as an introduction  to the subject of tropical geometry. For this purpose, we refer the reader to the forthcoming book of Maclagan--Sturmfels \cite{MS}. 

The paper is organized as follows: In Section 2, we introduce the analytification of an algebraic scheme $X$ over $K$ and we sketch how this fits into the theory of Berkovich analytic spaces. In Section 3, we define the tropicalization map and the tropicalization of a closed subscheme $X$ of the split torus $\Tor$ over $K$. In Section 4, we study models over the valuation ring. For a potentially integral point of the generic fibre, we define its reduction to the special fibre. We compare this with the reduction map from the theory of strictly affinoid algebras. In Section 5, the initial degeneration of a closed subscheme of $\Tor$ at $\omega 
\in \rdop^n$ is studied leading to an alternative characterization of the tropicalization.

In Section 6, we investigate normal affine toric schemes over a valuation ring associated to polyhedra. In Section 7, we  globalize these results assigning a normal toric scheme to every admissible fan in $\rdop^n \times \rdop_+$. This is rather new and it generalizes the theory from \cite{KKMS} worked out in the special case of a discrete valuation. In Section 8, we introduce the tropical cone of $X$ as a subset of $\rdop^n \times \rdop_+$. This new notion can be seen as a degeneration of the tropical variety with respect to the given valuation to the tropical variety with respect to the trivial valuation. It is very convenient to work with the tropical cone in the framework of tropical schemes over a valuation ring.

In Section 9, we study  projectively embedded toric varieties over the valuation ring which are not necessarily normal. This generalizes work of Eric Katz. In Section 10, we show that the tropical variety is a subcomplex of the Gr\"obner complex. In Section 11, we study the closure of $X$ in a toric scheme over the valuation ring and in Section 12, we generalize Tevelev's tropical compactifications to our setting. We introduce tropical multiplicities in Section 13 leading to the Sturmfels--Tevelev multiplicity formula for tropical cycles. In Section 14, we characterize proper compactifications of $X$ in a toric scheme which intersect all orbits properly. In the appendix, we collect results from convex geometry which are needed in the paper.

\vspace{3mm}

\centerline{\it Terminology}

In $A \subset B$, $A$ may be  equal to $B$. The complement of $A$ in $B$ is denoted by $B \setminus A$ \label{setminus} as we reserve $-$ for algebraic purposes. The zero is included in $\ndop$ and in $\rdop_+$.

All occurring rings and algebras are commutative with $1$. If $A$ is such a ring, then the group of multiplicative units is denoted by $A^\times$. A variety over a field is a separated reduced scheme of finite type. We denote by $\overline{F}$ an algebraic closure of the field $F$.

For $\mb \in \zdop^n$, let $\xb^\mb:=x_1^{m_1} \cdots x_n^{m_n}$ and $|\mb|:=|m_1|+ \cdots + |m_n|$. The standard scalar product of $\ub,\ub' \in \rdop^n$ is denoted by  $\ub \cdot \ub':=u_1 u_1' + \dots + u_n u_n'$. The terminology from convex geometry is explained in the appendix.

In the whole paper, the base field is a {\it valued field} $(K,v)$ which means that the field $K$ is endowed with a non-archimedean absolute value $|\phantom{a}|$ which might be trivial. The corresponding valuation is $v:=-\log(|\cdot|)$ with value group $\Gamma:=v(K^*)$. We get a valuation ring $K^\circ:=\{\alpha \in K \mid v(\alpha)\geq 0\}$ with maximal ideal $K^{\circ \circ} :=\{\alpha \in K \mid v(\alpha) > 0\}$ and residue field $\tilde{K}:=K^\circ / K^{\circ \circ}$. Note that $K=K^\circ=\tilde{K}$ if the valuation is trivial. We call $v$ a {\it discrete valuation} if $\Gamma \cong \zdop$.

We fix   a free abelian group $M$ of rank $n$ and $N:= \Hom(M, \zdop)$ the dual abelian group. An element of $M$ is usually denoted by $u$ and an element of $N$ is usually denoted by $\omega$. We get the duality pairing $\langle u, \omega \rangle := \omega( u)$. 
We have the split torus  $\tdop = \Spec(K^\circ[M])$  over $K^\circ$ with generic fibre $T$. Then  $M$ might be seen as the character group of this torus and the character corresponding to $u \in M$ is denoted by $\chi^u$. If $G$ is an abelian group, then $N_G:=N \otimes_\zdop G$ denotes the base change of $N$ to $G$. Similarly, $\tdop_A$ denotes the base change of $\tdop$ to a $K^\circ$-algebra $A$.

\vspace{3mm}
\small
The author thanks Matt Baker, Jose Burgos, Dustin Cartwright, Antoine Chambert-Loir, Qing Liu,  Sam Payne, C\'edric P\'epin, Joe Rabinoff, Martin Sombra, Alejandro Soto and Bernd Sturmfels for helpful comments. Special thanks also for Alejandro Soto for producing the pictures of the paper. I am very grateful to the referee for his careful reading and his suggestions which improved the paper a lot.
\normalsize

\section{Analytification} \label{Section:Analytification}

In this section, we recall the construction of the Berkovich analytic space $\Xan$ associated to an algebraic scheme $X$ over the  field $K$ with  non-archimedean complete absolute value $|\phantom{a}|$ and corresponding valuation $v:=-\log(|\cdot|)$. Note that completeness is no restriction of generality as analytic constructions are always performed over complete fields. In general, this may be achieved by base change to the completion of $K$. The topological part \ref{affine analytification}--\ref{analytification of a morphism}  of the construction  is elementary and essential for the understanding of the whole paper. The remaining analytic part is more technical and may be skipped in a first reading. Details and proofs for this section may be found in \cite{Berk1} and \cite{Tem}.

\begin{art} \rm \label{affine analytification}
We start with the construction for an affine scheme $X=\Spec(A)$ of finite type over $K$.  Then the topological space underlying the {\it Berkovich analytic space associated to $X$} is
the set of multiplicative seminorms on $A$ extending the  given absolute value on $K$.  In other words, we consider maps $p: A \rightarrow \rdop_+$ characterized by the properties
\begin{itemize}
 \item[(a)] $p(fg)=p(f)p(g)$
 \item[(b)] $p(1)=1$
 \item[(c)] $p(f+g) \leq p(f)+p(g)$
 \item[(d)] $p(\alpha)=|\alpha|$
\end{itemize}
for all $f,g \in A$ and $\alpha \in K$. It is easy to see that the triangle inequality (c) is equivalent to the ultrametric triangle inequality 
$$p(f+g) \leq \max(p(f),p(g)).$$ 
We denote this analytification of $X$  by $\Xan$ and we endow it with the coarsest topology such that the maps $\Xan \rightarrow \rdop, \, p \mapsto p(f)$ are continuous for every $f \in A$. We embed the set of closed points of $X$  into $\Xan$ by mapping $P$ to the seminorm $p$ given by $p(f)=|f(P)|$.
\end{art}

\begin{rem} \rm \label{rationality of seminorms}
Let $X=\Spec(A)$ be an affine scheme of finite type over $K$. For $p \in \Xan$, the integral domain $A/\{a \in A \mid p(a)=0\}$ is endowed with a canonical multiplicative norm induced by $p$. We conclude that its quotient field $L$ is endowed with an absolute value $|\phantom{a}|_w$ extending $|\phantom{a}|$. The canonical homomorphism $A \rightarrow L$ gives an $L$-rational point $P$ of $\Spec(A)$ and we may retrieve the seminorm by $p(a)=|a(P)|_w$ for any $a \in A$.

Conversely, any valued field $(L,w)$ extending $(K,v)$ and any $L$-rational point $P$ of $X$ give rise to an element $p \in \Xan$ by $p(a):=|a(P)|_w$, Obviously, different $L$-valued points might induce the same seminorm on $A$.
\end{rem}

\begin{lem} \label{analytic base change}
Let $X=\Spec(A)$ be an affine scheme of finite type over $K$ and let $(F,u)$ be a complete valued field extending $(K,v)$. Then the restriction map of seminorms gives a continuous surjective map from $(X_F)^{\rm an}$ onto $\Xan$.
\end{lem}

\proof Continuity is obvious from the definitions. For $p \in \Xan$, there is a valued field $(L,w)$ and $P \in X(L)$ as in Remark \ref{rationality of seminorms}. By Lemma \ref{composite of valued fields} below, there is a valued field $(F',u')$ extending both $(L,w)$ and $(F,u)$. We conclude that $P$ is also an $F'$-rational point of $X_F$ and hence it gives rise to a seminorm $p' \in (X_F)^{\rm an}$. By construction, $p'$ extends $p$ proving surjectivity. \qed

\begin{art} \rm \label{analytification of an algebraic variety} 
For any scheme $X$ of finite type over $K$, we choose an open affine covering $\{U_i\}_{i \in I}$. Then we define the topological space underlying the {\it Berkovich analytic space $\Xan$ associated to $X$} by glueing the spaces $U_i^{\rm an}$ obtained in \ref{affine analytification}. We get a topological space which is locally compact. It is Hausdorff if and only if $X$ is separated over $K$. 
\end{art}

\begin{art} \rm \label{analytification of a morphism}
If $\varphi:X \rightarrow Y$ is a morphism of schemes of finite type over $K$, then we have a canonical map $\varphi^{\rm an}: \Xan \rightarrow \Yan$ between the associated Berkovich analytic spaces. It is easy to see that it is enough to define the map locally, i.e. we may assume that $X$ and $Y$ are affine. Then we set $\varphi^{\rm an}(p):= p \circ \varphi^\sharp$ for any multiplicative seminorm $p$ on $\Ocal(X)$. 
\end{art}

\begin{art} \rm \label{closed points}
We will use analytifications also in situations where the valuation $v$ is not complete. Then we define $\Xan$ over the completion $K_v$ of $K$. Let $\overline{K}$ be the algebraic closure of $K$ in an algebraic closure of $K_v$. The absolute value on $\overline{K}$ depends on the embedding of $K$ into the algebraic closure of $K_v$ if $K$ is not complete.  Then every $\overline{K}$-rational point of $X$  induces a point in $\Xan$ as in Remark \ref{rationality of seminorms}. 

If the valuation is non-trivial, then the image of $X(\overline{K})$ is dense in $\Xan$. To see this, we may assume that $X$ is an irreducible affine variety of dimension $d$. By Noetherian normalization, there is a finite map $\varphi$ from $X$ onto the affine space $\adop_K^d$. Let $U$ be a non-empty open subset of $\Xan$. By  Lemma 3.2.7 in \cite{Berk1}, the map $\varphi^{\rm an}$ is open. Since the claim is obvious for $\adop^d$, there is a $\overline{K}$-rational point in $\varphi^{\rm an}(U)$. We choose a preimage in $U$ which has to be $\overline{K}$-rational  by finiteness of $\varphi^{\rm an}$. 
\end{art}

\begin{art} \rm \label{Tate algebra}
Now $K$ is complete again. We will explain below how $\Xan$ is endowed with an analytic structure. Of course, the analytic structure will depend on the underlying scheme structure.
First we handle the case of the affine space $\adop^n:= \Spec(K[x_1, \dots, x_n])$. For $f(\xb)=\sum_\mb \alpha_\mb \xb^\mb \in K[x_1, \dots, x_n]$, we have the {\it Gauss norm}
\begin{equation} \label{Gauss norm}
|f(\xb)|=\max_\mb|\alpha_\mb|. 
\end{equation}
The {\it Tate algebra} is defined as
$$K \langle x_1, \dots, x_n \rangle:= \{\sum_\mb \alpha_\mb \xb^\mb \in K[[x_1, \dots, x_n]] \mid \lim_{|\mb| \to \infty} |a_\mb|= 0\}$$
and it is the completion of $K[x_1, \dots, x_n]$ with respect to the Gauss norm. The corresponding Banach norm $|\phantom{a}|$ on $K \langle x_1, \dots, x_n \rangle$ is also defined by \eqref{Gauss norm}. The closed ball $\bdop^n$ of radius $1$ in $\adop^n$ is defined as the set of multiplicative seminorms on $K \langle x_1, \dots, x_n \rangle$ which are bounded by the Gauss norm, i.e. we have again properties (a)--(d) from \ref{affine analytification} for all $a,b \in K \langle x_1, \dots, x_n \rangle$ and the additional property $p(f) \leq |f|$. Note that a closed point of $\adop^n$ is in $\bdop^n$ if and only if all its coordinates have absolute value at most $1$.  It is easy to see that the supremum norm on $\bdop^n$ is equal to the Gauss norm.

More generally, we may consider  $\rb=(r_1, \dots, r_n)$ for strictly positive real numbers $r_1, \dots, r_n$. Then the Banach algebra $K\langle r_1^{-1}x_1, \dots, r_n^{-1}x_n \rangle $ is given by completion of $K[x_1, \dots, x_n]$ with respect to the weighted Gauss norm 
$$|f(\xb)|_\rb := \max_\mb |\alpha_\mb| \rb^\mb.$$
If we repeat the above construction, we get the closed ball $\bdop_\rb^n$ of radius $\rb$ in $\adop^n$. It is easy to see that $(\adop^n)^{\rm an}$ may be covered by a union of such balls. They serve as compact charts for the analytic structure of $(\adop^n)^{\rm an}$.
\end{art}

\begin{art} \rm \label{affinoid algebras}
An {\it affinoid algebra} is a Banach algebra $(\Acal, \| \phantom{a} \|)$ which is  isomorphic to $K\langle r_1^{-1}x_1, \dots, r_n^{-1}x_n \rangle / I$ for an ideal $I$ and such that the norm $\| \phantom{a} \|$ is equivalent to the quotient norm 
$$\|f+I\|_{\rm quot}:=\inf\{\|g\| \mid g \in f+I\}$$
on $\Acal/I$. It is called a {\it strictly affinoid algebra} if we may choose $r_i=1$ for all $i=1, \dots,n$. 

The Banach norm does not matter if the valuation is non-trivial (see \cite{Tem}, Fact 3.1.15). In this case, all such Banach norms on $\Acal$ are equivalent and every homomorphism between affinoid algebras is bounded. In classical rigid geometry as in \cite{BGR}, one considers only strictly affinoid algebras and they are called affinoid algebras there. 
\end{art}

\begin{art} \rm \label{spectral radius}
For an affinoid algebra as above, the {\it spectral radius} is defined by
$\rho(a):= \inf\{\|a^n\|^{1/n} \mid n \geq 1 \}$
for $a \in \Acal$. We set
$$\Acal^\circ:= \{a \in \Acal \mid \rho(a) \leq 1 \}, 
\quad \Acal^{\circ \circ}:=\{a \in \Acal \mid \rho(a) < 1 \}.$$
and the {\it residue algebra} is defined by $\tilde{\Acal}:= \Acal^\circ / \Acal^{\circ \circ}$. 
\end{art}

\begin{art} \rm \label{Berkovich spectrum} 
The {\it Berkovich spectrum} $\Mcal(\Acal)$ of a $\kdop$-affinoid algebra $\Acal$ is defined as the set of multiplicative bounded seminorms $p$ on $\Acal$, i.e. 
for all $a,b \in \Acal$, we have
\begin{itemize}
 \item[(a)] $p(ab)=p(a)p(b)$
 \item[(b)] $p(1)=1$
 \item[(c)] $p(a+b) \leq p(a)+p(b)$
 \item[(d)] $p(a)\leq  \rho(a).$
\end{itemize}

It is endowed with the coarsest topology such that the maps $p \mapsto p(a)$ are continuous for all $a \in \Acal$. We get a  compact space. The spectral radius $\rho(a)$ turns out to be equal to the supremum seminorm $\sup\{p(a) \mid p \in \Mcal(\Acal)\}$ of $a \in \Acal$.
\end{art}

\begin{ex} \rm \label{rational domain}
Let $\Acal = K\langle r_1^{-1}x_1, \dots, r_n^{-1}x_n \rangle / I$. A {\it rational subdomain} of $X:=\Mcal(\Acal)$ is  defined by
$$X \left(\Sb^{-1}\frac{\mathbf f}{g}\right) := \{x \in X \mid |f_j(x)| \leq s_j|g(x)|, \, j=1,\dots,m\}$$
where $g,f_1, \dots, f_m$ generate the unit ideal in $\Acal$ and $s_1, \dots, s_m>0$. The corresponding affinoid algebra is
$$\Acal \left\langle \Sb^{-1}\frac{\mathbf f}{g}\right\rangle := K \langle \rb^{-1}\mathbf x, s_1^{-1}y_1, \dots, s_m^{-1}y_m \rangle / \langle I, g(\mathbf x) y_j - f_j(\mathbf x) \mid j=1, \dots ,m \rangle$$
(see \cite{Berk1}, Remarks 2.2.2). 
\end{ex}

\begin{art} \rm \label{Berkovich analytic spaces}
We will not give the precise definition of a {\it Berkovich analytic space} $X$ (see \cite{Berk2} for details). Roughly speaking it is a topological space endowed with an atlas such that each chart is homeomorphic to the Berkovich spectrum of an affinoid algebra and then there are some compatibility conditions. Analytic functions on such a chart $\Mcal(\Acal)$ are given by the elements of $\Acal$.

A morphism $\varphi: X_1 \rightarrow X_2$ between Berkovich spaces $X_1$ and $X_2$  is a continuous map such that for every chart $U_1$ of $X_1$ with $\varphi(U_1)$ contained in a chart $U_2$ of $X_2$ and every analytic function $f_2$ on $U_2$, the function $\varphi^\sharp(f_2):=f_2 \circ \varphi$ is an analytic function on $U_1$ and $\varphi$ is induced by $\varphi^\sharp$.
\end{art}

\begin{art} \rm \label{Algebraic example}
If $X$ is a scheme of finite type over $K$ as at the beginning, then $\Xan$ is a Berkovich analytic space. As charts, we may choose $U^{\rm an} \cap \bdop_\rb^n$, where $U$ is an affine open subset of $X$  realized as a closed subset of $\adop^n$ and $\bdop_\rb^n$ is a closed ball in $\adop^n$. 
Serre's GAGA principle holds also in the non-archimedean framework. For details, we refer to \cite{Berk1}, \S 3.4 and \S 3.5. \end{art}

\begin{rem} \rm \label{analytic domains}
In Example \ref{rational domain}, we have defined rational subdomains of the Berko\-vich spectrum $X:=\Mcal(\Acal)$. More generally, one can define {\it affinoid subdomains} of $X$ by a certain universal property. They are Berkovich spectra contained in $X$ which are used for localization arguments on Berkovich spaces. For details, we refer to \cite{Berk1}, Section 2.2. By the Gerritzen--Grauert theorem, every affinoid subdomain is a union of rational domains if the valuation $v$ is non-trivial.

Roughly speaking, an {\it analytic subdomain} of a Berkovich analytic space $X$ is a subset which behaves locally like an affinoid subdomain. For a precise definition and for properties, we refer to \cite{Berk1}, Section 3.1. In this paper, we need only analytic functions on affinoid subdomains of $\Xan$ where they are just elements of the corresponding affinoid algebra. However, it should be noted that analytic functions form a sheaf on open subsets giving $\Xan$ the structure of a locally ringed space (see \cite{Berk1}, or \cite{BPS}, \S 1.2, for a neat description).
\end{rem}

\section{Tropicalization}

In this section, we consider a closed subscheme $X$ of the multiplicative torus $T$ over the valued field $(K,v)$ and we will define the tropical variety $\Trop_v(X)$ associated to $X$.  The tropicalization process is a transfer from algebraic geometry  to convex geometry in $\rdop^n$. We will use the analytifications $\Xan$ and $\Tan$ from the previous section which are always performed over the completion of $K$. 

\begin{art} \rm \label{multiplicative torus}
Let $M$ be a free abelian group of rank $n$ and let $N = \Hom(M,\zdop)$ be the dual group. Then we consider the multiplicative torus $T:=\Spec(K[M])$ with character group $M$. We have the {\it tropicalization map} 
$${\trop_v}: \Tan \rightarrow N_\rdop, \quad p \mapsto {\trop_v}(p),$$
where ${\trop_v}(p)$ is the element of $N_\rdop=\Hom(M,\rdop)$ given by 
$$\langle u, {\trop_v}(p) \rangle:= - \log(p(\chi^u))$$
with $\chi^u$ the character of $T$ corresponding to $u\in M$. Choosing coordinates $x_1, \dots, x_n$ on $T=\Tor$, we may identify $M$ and $N$ with $\zdop^n$ and we get an explicit description 
$${\trop_v}: \Tan \rightarrow \rdop^n, \quad p \mapsto (-\log(p(x_1)), \dots, -\log(p(x_n))).$$
It is immediate from the definitions that the map ${\trop_v}$ is continuous. This is the big advantage of working with Berkovich analytic spaces in this framework as we may use their nice topological properties.
\end{art}

\begin{Def} \rm \label{tropical variety}
We define the {\it tropical variety associated to $X$} by $\Trop_v(X):={\trop_v}(\Xan)$. In Section \ref{tropical mult},   we will  complete the definition of a tropical variety by assigning  certain weights.
\end{Def}

In the following result, we refer the reader to the appendix for the terminology borrowed from convex geometry.

\begin{thm}[Bieri--Groves] \label{Bieri-Groves theorem}
$Trop_v(X)$ is a finite union of $\Gamma$-rational polyhedra in $N_\rdop$. If $X$ is of pure dimension $d$, then we may choose all the polyhedra $d$-dimen\-sional. 
\end{thm}

\proof The  proof is given in \cite{BG}, Theorem A. Note that even the definition of $\Xan$ occurs implicitly in this paper.  For a translation into tropical language, we refer to \cite{EKL}, Theorem 2.2.3. In Theorem \ref{support of Groebner}, we will give a proof of the first statement using the Gr\"obner fan. A proof for  dimensionality is given in \cite{Gu3}, Proposition 5.4, which generalizes to closed analytic subvarieties. \qed

\begin{rem} \rm \label{tropical variety for trivial absolute value}
If the absolute value on $K$ is trivial, then a $\Gamma$-rational polyhedron is just a rational polyhedral cone. In this case, we conclude that $\Trop_v(X)$ is a finite union of such cones.
\end{rem}

We illustrate the advantage of Berkovich spaces by giving the proof of the following well-known result (see \cite{EKL}, Theorem 2.2.7). 

\begin{prop} \rm \label{connected}
If $K$ is complete and $X$ is connected, then $\Trop_v(X)$ is connected.
\end{prop}

\proof If $X$ is connected, then $\Xan$ is connected (\cite{Berk1}, Theorem 3.4.8 and Theorem 3.5.3). This is a rather nontrivial fact. By continuity of the tropicalization map, we conclude that $\Trop_v(X)$ is connected. \qed

\begin{rem} \rm \label{connected in codimension 1}
This is wrong if $K$ is non-complete, but still true if $K$ is algebraically closed. For more details, we refer to the paper \cite{CP} of Cartwright and Payne, where they prove also that $\Trop_v(X)$ is connected in codimension $1$ if $X$ is irreducible over a complete or algebraically closed field $K$.
\end{rem}

\begin{prop} \label{tropical base change}
Let $(L,w)$ be a valued field extending $(K,v)$. Then we have $\Trop_w(X_L)=\Trop_v(X)$.
\end{prop}

\proof Let $\varphi: (X_L)^{\rm an} \rightarrow \Xan$ be the restriction map of seminorms. We have seen in Lemma \ref{analytic base change} that $\varphi$ is surjective. Using $\trop_w=\trop_v \circ \varphi$, we get the claim. \qed

\vspace{2mm}
The following result shows that our definition of a tropical variety agrees with the usual one.

\begin{prop} \label{classical tropical characterization}
Let $(L,w)$ be an algebraically closed valued field extending $(K,v)$ endowed with a non-trivial  absolute value $|\phantom{a}|_w$ and let $x_1, \dots, x_n$ be torus coordinates on $T$. Then $\Trop_v(X)$ is  equal to the closure of $$\{(-\log|x_1|_w, \dots , -\log|x_n|_w) \mid \xb \in X(L)\}$$ in $\rdop^n$.
\end{prop}

\proof By base change and Proposition \ref{tropical base change}, we may assume that $K$ is algebraically closed and that $v$ is non-trivial. We have seen in \ref{closed points}, that $X(K)$ is dense in $\Xan$ and hence continuity of the tropicalization map yields the claim. \qed

\section{Models over the valuation ring and reduction} \label{algebraic and formal models}

In this section, $(K,v)$ is  a valued field with valuation ring $K^\circ$ and residue field $\tilde{K}$. We will study models of a  scheme $X$ of finite type over $K$. The models are flat schemes over $K^\circ$ but not necessarily of finite type. We will obtain a model of a closed subscheme of $X$ by taking the closure. For integral points of a model, there is always a reduction modulo the maximal ideal $K^{\circ \circ}$ which is a point in the special fibre. We will compare it with the reduction from the theory of strictly affinoid algebras.

\begin{Def} \label{algebraic model} \rm 
A {\it  $K^\circ$-model}  of a scheme $X$ over $K$ is a flat scheme $\Xcal$ over $K^\circ$ with generic fibre $\Xcal_\eta:= \Xcal_K=X$. The special fibre $\Xcal_{\tilde{K}}$ of $\Xcal$ is denoted by $\Xcal_s$. The model $\Xcal$ is called {\it algebraic} if it is of finite type over $K^\circ$. 
\end{Def}

\begin{lem} \label{flat module}
A module $Q$ over $K^\circ$ is flat if and only if $Q$ is torsion-free.
\end{lem}

\proof Any flat module is obviously torsion-free. If the base is a valuation ring, then the converse holds. 
It is enough to show that the map $I \otimes_{K^\circ} Q \rightarrow Q$ is injective for every finitely generated ideal $I$ of $K^{\circ}$. As such an ideal is generated by a single element $\alpha$, injectivity follows immediately from  $Q$ torsion free.  \qed

\begin{art} \label{affine closure} \rm
Let $\Xcal =\Spec(A)$ be a flat scheme over $K^\circ$ with generic fibre $X=\Xcal_\eta$. Then we have $X=\Spec(A_K)$ for $A_K:=A \otimes_{K^\circ} K$. Note that flatness implies $A \subset A_K$. A closed subscheme $Y$ of $X$ is given by an ideal $I_Y$ in $A_K$. The {\it closure $\overline{Y}$ of $Y$ in $\Xcal$} is defined as the closed subscheme of $\Xcal$ given by the ideal $I_Y \cap A$.
\end{art}

\begin{prop} \label{characterization of affine closure}
The closure of $Y$ is the unique closed subscheme of $\Xcal$ with generic fibre $Y$ which is flat over $K^\circ$. 
\end{prop}

\proof It is clear that $A/(I_Y \cap A)$ is $K^\circ$-torsion free and hence flat over $K^\circ$ by Lemma \ref{flat module}. For every $f \in A_K$, there is a non-zero $\lambda \in K^\circ$ with $\lambda f \in A$. We conclude that $I_Y \cap A$ generates $I_Y$ as an ideal in $A_K$ and hence the generic fibre of $\overline{Y}$ is $Y$.

Let $\Ycal$ be any closed subscheme of $\Xcal$ with generic fibre $Y$ which is flat over $K^\circ$. Then $\Ycal$ is given by an ideal $J$ in $A$ such that $J$ generates $I_Y$ as an ideal in $A_K$. We conclude that $J \subset I_Y \cap A$. Hence we have a canonical homomorphism $h:A/J \rightarrow A/(I_Y \cap A)$. By flatness over $K^\circ$, we have $A/J \subset A_K/I_Y$ and $h$ factors  through this inclusion. We deduce that $h$ is one-to-one proving $J=I_Y \cap A$. \qed


\begin{cor} \label{affine closure and flat pull-back}
Let $\psi:\Xcal' \rightarrow \Xcal$ be a flat morphism of flat affine $K^{\circ}$-schemes with generic fibre $\psi_\eta:X' \rightarrow X$. Then we have 
$\overline{(\psi_\eta)^{-1}(Y)}=\psi^{-1}(\overline{Y})$ for a closed subscheme $Y$ of $X$, where 
 the closures are taken in $\Xcal'$ and $\Xcal$.
\end{cor}


\proof As $\psi^{-1}(\overline{Y})$ is a closed subscheme of $\Xcal'$ with generic fibre ${(\psi_\eta)^{-1}(Y)}$ which is flat over $K^\circ$, the claim follows from Proposition \ref{characterization of affine closure}. \qed

\begin{rem} \label{closure in general} \rm 
In particular, this shows that localization is compatible with taking the closure. Therefore  the closure may be defined in any flat scheme $\Xcal$ over $K^\circ$. Indeed, let $Y$ be a closed subscheme of $X:=\Xcal_\eta$. First, we define $\overline{Y}$ locally on affine charts as in \ref{affine closure} and then we glue the affine pieces to get a closed subscheme $\overline{Y}$ of $\Xcal$ by compatibility of the affine construction with localization. The closure is still characterized by Proposition \ref{characterization of affine closure}. Moreover,  Corollary \ref{affine closure and flat pull-back} immediately yields that the formation of the closure is compatible with flat pull-back. Note also that the underlying set of $\overline{Y}$ is the topological closure of $Y$ in $\Xcal$. 
\end{rem}

\begin{cor} \label{closure and base change}
Let $(L,w)$ be a valued field extension of $(K,v)$ and let $\Xcal$ be a flat scheme over $K^\circ$. For a closed subscheme $Y$ of $X=\Xcal_\eta$, we have $(\overline{Y})_{L^\circ}=\overline{Y_L}$ with closures taken in $\Xcal$ and $\Xcal_{L^\circ}$. 
\end{cor}

\proof Note that the base change morphism $\Xcal_{L^\circ} \rightarrow \Xcal$ is flat. Taking the closure depends only on the model and not on the base and hence compatibility with flat pull-back (Corollary \ref{affine closure and flat pull-back}) yields the claim. \qed

\begin{art} \rm \label{reduction of integral points}
For an $L^\circ$-integral point $P$ of $\Xcal$, the reduction $\pi(P) \in \Xcal_s$ is defined as the image of the closed point of $\Spec(L^\circ)$ with respect to the morphism $\Spec(L^\circ)\rightarrow \Xcal$ defining $P$. If $\Xcal =\Spec(A)$ is affine, then $\pi(P)$ is given by the prime ideal $\{a \in A \mid |a(P)|_w <1\}$ in $A$. 
\end{art}

\begin{art} \rm \label{reduction for Xan}
Let $X$ be a scheme of finite type over $K$ with $K^\circ$-model $\Xcal$. Our goal is to introduce a reduction map $\pi$ from $X$ to the special fibre of $\Xcal$. Such a map can be defined only at integral points and it turns out that it is better to work analytically.

We handle first the affine case, i.e. $\Xcal =\Spec(A)$ for an algebra $A$ of finite type over $K$. Then we define $X^\circ:=\{p \in \Xan \mid p(f) \leq 1 \; \forall f \in A\}$.  Note that $X^\circ$ is the set of points in $\Xan$ which are induced by an $L^\circ$-integral point of $X$ for some valued field extension $(L,w)$. Such points of $X$ are called {\it potentially integral}. If $p \in X^\circ$, then the {\it reduction} $\pi(p) \in \Xcal_s$ is given by the prime ideal $\{a \in A \mid p(a) <1\}$ in $A$.

In general, we define $X^\circ$ as the union of all $U^\circ:=\{p \in \Uan \mid p(f) \leq 1 \; \forall f \in A\}$, where $\Ucal=\Spec(A)$ ranges over the affine open subsets of $\Xcal$ and $U:=\Ucal_\eta$. It is clear that the notions coincide in the affine case. The points of $X^\circ$ are induced by the potentially integral points of $X$ as above. Proceeding locally, we get the {\it reduction map} $\pi:X^\circ \rightarrow \Xcal_s$.

Note that if $\Xcal$ is an algebraic $K^\circ$-model, then $X^\circ$ is a compact analytic subdomain of $\Xan$. Indeed, in the affine case we get an affinoid subdomain in $\Xan$ and in general, $X^\circ$ is a finite union of affinoids. If $\Xcal$ is a proper scheme over $K^\circ$, then rational and integral points are the same and hence $X^\circ = \Xan$.

If we assume that $K$ is endowed with a non-trivial complete valuation $v$ and if we assume that $\Xcal$ is an algebraic $K^\circ$-model, then we will see in \ref{comparision of reduction} that $X^{\circ}$ is the generic fibre of the completion of $\Xcal$ along the special fibre.



\end{art}

\begin{art} \rm \label{algebraic description of reduction map}
For a scheme $X$ of finite type over $K$ with algebraic $K^\circ$-model $\Xcal$, the reduction map $\pi$  can be described algebraically in the following way: We consider an $L^\circ$-integral point $P$ of $X$ for a valued field  $(L,w)$ extending $(K,v)$. Integrality means here that there is an affine chart $\Ucal$ of $\Xcal$ with affine coordinates $x_1, \dots, x_n$ such that $x_1(P), \dots, x_n(P) \in L^\circ$. Then $\pi(P)$ is the point of the special fibre $\Ucal_s$ given by using the coordinates modulo the maximal ideal $L^{\circ \circ}$ of $L^\circ$. Note that this point is not closed if the residue field $\tilde{L}$ is an infinite extension of $\tilde{K}$.
\end{art}

In the theory of strictly affinoid algebras introduced in \ref{affinoid algebras}, there is a similar concept of reduction which we study next. For this, we assume that the valuation $v$ on $K$ is non-trivial and complete. 

\begin{art} \rm \label{reduction for strictly affinoids}
Let $\Acal$ be a strictly affinoid $K$-algebra with Berkovich spectrum $Y=\Mcal(\Acal)$. We define the {\it reduction} of $Y$ by $\tilde{Y}:=\Spec(\tilde{\Acal})$ and the {\it special fibre} of $Y$ by $Y_s:=\Spec(\Acal^\circ/(K^{\circ \circ}\Acal^\circ))$. The reduction is an algebraic variety over the residue field $\tilde{K}$ (see \cite{BGR}, Corollary 6.4.3/1). Since the maximal ideal $K^{\circ \circ}$ of $K^\circ$ generates an ideal in $\Acal^\circ$ contained in $\Acal^{\circ \circ}$, we get a canonical surjective homomorphism $\Acal^\circ/(K^{\circ \circ}\Acal^\circ) \rightarrow \tilde{\Acal}$. This induces a canonical morphism $\tilde{Y} \rightarrow Y_s$. Since the spectral radius is power-multiplicative, it is clear that this morphism is a bijection.
 
We have a  map $Y \rightarrow \tilde{Y}$, given by mapping the seminorm $p$ to the prime ideal $\{a \in \Acal^\circ \mid p(a) < 1\}/\Acal^{\circ \circ}$ of $\tilde{\Acal}$. It induces a {\it reduction map} $\pi: Y \rightarrow Y_s$.

As in algebraic geometry,  a {\it Zariski closed} subset of $Y$ is the zero set of a subset of $\Acal$ and this leads to the {\it Zariski topology} on $Y$. Note that a Zariski open subset is dense in the Berkovich topology if and only if it is dense in the Zariski topology. 
\end{art}

\begin{lem} \label{lifting lemma} 
For a Zariski open  and dense subset $S$ of $Y$, we have $\pi(S)=Y_s$. 
\end{lem}

\proof We first note that the reduction map $\pi$ is surjective (see \cite{Berk1}, Proposition 2.4.4). If $z$ is a closed point of $Y_s$, then $\pi^{-1}(z)$ is an open non-empty subset of $Y$ (\cite{Berk1}, Lemma 2.4.1). By density, there is $y \in S$ with $\pi(y)=z$. In general, there is a complete valued field $(L,w)$ extending $(K,v)$ and an $L$-rational point $P$ of $Y$ such that $\pi(P)=z$ (see Remark \ref{rationality of seminorms}). Then the reduction $\pi_L(P)$ of $P$ in the special fibre of $Y_L=\Mcal(\Acal \widehat{\otimes}_K L)$ is $\tilde{L}$-rational.  The preimage $S_L$ of $S$ is Zariski open and dense in $Y_L$. Using the above, there is $y_L \in S_L$ with  $\pi_L(y_L)=\pi_L(P)$. Then we have $\pi(y)=z$ for the image $y$ of $y_L$ in $Y$. \qed

\begin{art} \rm \label{comparision of reduction}
We compare the two concepts for a reduction map in the following situation: Let $(K,v)$ be an arbitrary valued field and let $(L,w)$ be a complete valued field extending $(K,v)$ with $w$ non-trivial. We consider a flat affine scheme $\Xcal=\Spec(A)$ of finite type over $K^\circ$ with generic fibre $X=\Spec(A_K)$. For convenience, we choose coordinates $x_1, \dots ,x_n$ on $\Xcal$, i.e. $A=K^\circ[x_1, \dots, x_n]/I$ for an ideal $I$ in $K^\circ[x_1, \dots, x_n]$. Then we complete the base change $\Xcal_{L^\circ}$ along the special fibre (more precisely, we take the $\nu$-adic completion for some non-zero $\nu \in K^{\circ \circ}$) to get a flat formal scheme $\Ymcal=\Spf(L^\circ \langle x_1, \dots ,x_n \rangle/ \langle I \rangle)$ over $L^\circ$ (see \cite{Ul}). The generic fibre $Y$ of $\Ymcal$ is the Berkovich spectrum of the strictly affinoid algebra $\Acal$ defined by
$$\Acal := (L^\circ \langle x_1, \dots ,x_n \rangle/ \langle I \rangle) \otimes_{L^\circ} L =L \langle x_1, \dots, x_n \rangle / \langle I \rangle.$$
By construction, $Y$ is the affinoid subdomain $(X_L)^{\circ}=\{p \in (X_L)^{\rm an} \mid p(x_1) \leq 1, \dots, p(x_n) \leq 1 \}$ in $(X_L)^{\rm an}$. It is easy to see that we have a commutative diagram
\begin{equation} 
\begin{CD} \label{commutative square}
Y @>\pi>> Y_s\\
@VV{}V    @VV{}V\\
X^\circ @>{\pi }>>\Xcal_s
\end{CD}
\end{equation}
where the vertical maps are induced by base change and the horizontal maps are the reduction maps. Applying Theorem 6.3.4/2 of \cite{BGR} to the surjective homomorphism $L\langle x_1, \dots ,x_n \rangle \rightarrow \Acal$, it follows that the canonical morphism $\tilde{Y} \rightarrow \Ymcal_s$ is a finite map. Both spaces have dimension equal to $\dim(X)$ and an easy localization argument shows that this finite map is surjective. As the base change morphism $\Ymcal_s= (\Xcal_L)_s \rightarrow \Xcal_s$ is also surjective and since the canonical morphism $\tilde{Y} \rightarrow Y_s$ is a bijection, we deduce that the morphism $Y_s \rightarrow \Xcal_s$ is surjective. By Lemma \ref{lifting lemma}, we conclude that the reduction map $\pi:X^\circ \rightarrow \Xcal_s$ is surjective.
\end{art}

\begin{prop} \label{finite generated lifting lemma}
Let  $\Xcal$ be a flat scheme of finite type over $K^\circ$ with generic fibre $X$ and let $U$ be an open dense subset of $X$. Then we have $\pi(\Uan \cap X^\circ) = \Xcal_s$. If $K$ is algebraically closed and $v$ is non-trivial, then every $\tilde{K}$-rational point of $\Xcal_s$ is the reduction of a $K^\circ$-integral point contained in $U$.
\end{prop}

\proof We may assume that $\Xcal$ is affine and hence we are in the situation of  \ref{comparision of reduction}. We choose a valued field $(L,w)$ extending $(K,v)$ with $w$ complete and non-trivial. Since the base change $U_L$ is open and dense in $X_L$, we conclude that $S:= (U_L)^{\rm an} \cap Y$ is Zariski open and dense in the analytic space $Y$. Using surjectivity of the map $Y_s \rightarrow \Xcal_s$ and Lemma \ref{lifting lemma}, we deduce $\pi(\Uan \cap X^\circ) = \Xcal_s$ from the commutative diagram \eqref{commutative square}.

If $K$ is algebraically closed and $v$ is non-trivial, then we note first that $\tilde{K}$ is algebraically closed (see \cite{BGR}, 3.4.1). For any closed point $z$ of the special fibre, the above and anticontinuity of the reduction map show that $\pi^{-1}(z) \cap \Uan$ is a non-empty open subset of $X^{\circ}$. If we embed $U$ into affine space, then we see that $\pi^{-1}(z) \cap \Uan$ is the intersection of $\Uan$ with an open ball and hence $\pi^{-1}(z) \cap \Uan$ is even an open subset of $\Xan$. Density of the $K$-rational points yields the claim (see \ref{closed points}). \qed


\begin{ex} \rm \label{projective line and reduction}
We assume that the absolute value is trivial on $K$. Let $X=\pdop^1_K$ with projective coordinates $x_0,x_1$. For $i=0,1$, we consider the affine charts $U_i:=\{x \in X \mid x_i \neq 0\}$ isomorphic to $\adop^1_K$. For any $r >0 $, we get an element $p_r \in U_0^{\rm an}$ given as the seminorm $p_r(f):=\max_i|a_i| r^i$ for $f(y)=\sum_i a_i y^i \in K[y]$ with $y:=\frac{x_1}{x_0}$. If $r \leq 1$, then $p_r(f)= r^j$ for $j$  minimal with $a_j \neq 0$. Then we have $p_r \leq 1$ and hence the reduction of $p_r$ is defined by $\pi(p_r)=\{f \in K[y] \mid p_r(f)<1\} \in \Spec(K[y])\subset \pdop^1_K$. If $r<1$, then $\pi(p_r)=(1:0)$. If $r=1$, then $\pi(p_r)$  is the generic point of $\pdop^1_K$.  If $r >1$, then we use the other chart $U_1$ with affine coordinate $z:=\frac{x_0}{x_1}$. For $g(z)=\sum_i a_i z^i \in K[z]$, we have $p_r(g)=\max_i|a_i| r^{-i}$ and hence $\pi(p_r)=(0:1)$.
\end{ex}

\section{Initial degeneration}

In this section, we study the initial degeneration $\In_\omega(X)$ of a closed subscheme $X$  of the multiplicative torus $T=\Tor$ over the valued field $(K,v)$ at $\omega \in N_\rdop$. We follow here the original definition of the initial degeneration using a translation to the origin of the torus. Then $\In_\omega(X)$ is a closed subscheme of the torus $\tdop_{\tilde{K}}$ which is only well-defined up to translations. This approach fits very well to  Hilbert schemes as we will see in Section \ref{Section: The Groebner complex}. For an intrinsic approach, we refer to \cite{OP}.



\begin{Def} \rm \label{initial degeneration at t}
Let $(L,w)$ be a valued field extending $(K,v)$ and let $t \in T(L)$. Then the {\it initial degeneration of $X$ at $t$} is defined as the special fibre of the closure of $t^{-1}X_L$ in the split multiplicative torus $\tdop_{L^\circ}$ over the valuation ring $L^\circ$. It is a closed subscheme of the split torus $\tdop_{\tilde{L}}$ over the residue field $\tilde{L}$ which we denote by $\In_t(X)$.
\end{Def}

\begin{lem} \label{composite of valued fields}
Let $(L,w)$ and $(L',w')$ be valued fields extending $(K,v)$. Then there is a valued field $(L'',w'')$ extending $(L,w)$ and $(L',w')$. 
\end{lem}

\proof This is proved in \cite{Du}, \S 0.3.2, using Berkovich's theory. \qed

\begin{prop} \label{dependence on tropicalization}
Let $(L,w)$ and $(L',w')$ be valued fields extending $(K,v)$. Suppose that there is $\omega \in N_\rdop$ such that $\omega = \trop_w(t)=\trop_{w'}(t')$ for $t \in T(L)$ and $t' \in T(L')$. For any field $(L'',w'')$ as in Lemma \ref{composite of valued fields}, there is $g \in \tdop(\tilde{L}'')$ with 
\begin{equation} \label{equivalence relation}
\In_{t'}(X)_{\tilde{L}''}=g \cdot \In_t(X)_{\tilde{L}''}.
\end{equation}
\end{prop}

\proof Since  $t$, $t'$ have the same tropicalizations, the point $t/t' \in T(L'')$ is in fact an $(L'')^\circ$-integral point of $\tdop$ and hence it has a well-defined reduction $g \in \tdop(\tilde{L}'')$. The relation 
$$(t')^{-1}X_{L''}=(t/t') \cdot t^{-1}X_{L''}$$
and Corollary \ref{closure and base change} give immediately the claim. \qed

\begin{art} \rm \label{initial degeneration at omega}
The proposition shows that the initial degeneration depends essentially only on $\omega$. 
For any $\omega \in N_\rdop$, there is a valued field $(L,w)$ extending $(K,v)$ and $t \in T(L)$ with $\trop_w(t)=\omega$. We define the {\it initial degeneration $\In_\omega(X)$ of $X$ at $\omega$} as  $\In_t(X)$ which is well-defined as an equivalence class for the equivalence relation \eqref{equivalence relation}. We call the residue field $\tilde{L}$ or any extension of it a {\it  field of definition} for $\In_\omega(X)$.
\end{art}

\begin{prop} \label{initial degeneration and base change}
Let $(L,w)$ be a valued field extending $(K,v)$ and let $\omega \in N_\rdop$. Then we have $\In_\omega(X_L)=\In_\omega(X)$ up to the equivalence relation given by \eqref{equivalence relation}.
\end{prop}

\proof By Corollary \ref{closure and base change}, the formation of the closure is compatible with  base change and this yields easily the claim. \qed

\vspace{2mm}

The next result is called the {\it fundamental theorem} of tropical algebraic geometry. It is due to Kapranov in the hypersurface case (unpublished manuscript, later incorporated in \cite{EKL}) and to Speyer--Sturmfels \cite{SS}, Draisma \cite{Dr}, Payne \cite{Pay} in general.

\begin{thm} \label{fundamental theorem}
For a closed subscheme $X$ of $T$, the tropical variety $\Trop_v(X)$ may be characterized in the following two equivalent ways:
\begin{itemize}
\item[(a)] $\Trop_v(X)= {\trop_v}(\Xan)$
\item[(b)] The set $\{\omega \in N_\rdop \mid \In_\omega(X) \neq \emptyset\}$ in $N_\rdop$ is equal to $\Trop_v(X)$. 
\end{itemize}
\end{thm}

\proof We have to prove that $\omega \in N_\rdop$ is in $\trop_v(\Xan)$ if and only if $\In_\omega(X) \neq \emptyset$. By base change, we may assume that $(K,v)$ is a non-trivially valued complete algebraically closed field such that $\omega = \trop_v(t)$ for some $t \in T(K)$ (see Propositions \ref{tropical base change} and \ref{initial degeneration and base change}). Passing to $t^{-1}X$, we may assume that $t=e$ and $\omega =0$. Let $\Xcal$ be the closure of $X$ in $\tdop$. It is an algebraic $K^\circ$-model of $X$. We recall that the reduction map $\pi$ to the special fibre $\Xcal_s$ is defined on the affinoid subdomain $X^\circ$ of $\Xan$ from \ref{reduction for Xan}. Using that the regular function on $\Xcal$ are generated by the characters $\chi^u$ with $u$ ranging over a basis of $M$, we deduce easily that $X^\circ=\trop_v^{-1}(0)\cap \Xan$. We have seen in Proposition \ref{finite generated lifting lemma} that $\pi(X^\circ)=\Xcal_s$ and hence $\In_0(X)=\Xcal_s$ is empty if and only if $\trop_v^{-1}(0)\cap \Xan=\emptyset$. The latter is equivalent to $0 \not \in \trop_v(\Xan)$ proving the claim. \qed

\begin{rem} \label{initial forms and ideals} \rm
Initial degenerations may be studied using methods from the theory of Gr\"obner bases. Let $(L,w)$ be a valued field extending $(K,v)$ and let $t \in T(L)$ with $\omega = \trop_w(t)$. For a Laurent polynomial $f=\sum_{u \in M} \alpha_u \chi^u \in K[M]$, we define {\it the initial form} $\In_t(f) \in \tilde{L}[M]$ in the following way. If $f=0$, then  $\In_t(f)$ is the zero polynomial in $\tilde{L}[M]$. If $f \neq 0$, then we choose $\lambda \in L$ with $v(\lambda)= \min_{u \in M} v(\alpha_u) + \langle u, \omega \rangle$ and we set $\In_t(f):=\sum_{u \in M} \pi(\lambda^{-1}\alpha_u \chi^u(t)) \chi^u$, where $\pi: L^\circ \rightarrow \tilde{L}$ is the reduction map. Note that the initial form is only well-defined up to multiplication by $\tilde{L}^\times$. 

For a closed subscheme $X$ of $T$ given by the ideal $I_X$ in $K[M]$, we define the {\it initial ideal} of $X$ at $t$ as the ideal $\In_t(I_X)$ in $\tilde{L}[M]$ generated by $\{\In_t(f) \mid f \in I_X\}$. By construction, the initial degeneration $\In_t(X)$ from \ref{initial degeneration at t} is given by the initial ideal $\In_t(I_X)$ in $\tilde{L}[M]$. If there is a canonical homomorphism $\tau:v(L^\times) \rightarrow L^\times$ with $v(\tau(r))=r$ for all $r \in \rdop$, then the initial form $\In_\omega(f)$ of $f$ can be defined canonically at $\omega \in N_\rdop$ using $t:=\tau(\omega)$ and $\lambda := \tau(\min_{u \in M} v(\alpha_u) + \langle u, \omega \rangle)$. For more details, we refer to \cite{MS}, \S 2.4. For an intrinsic approach to initial forms $\In_\omega(f)$ and initial degenerations $\In_\omega(X)$  without using translations to the origin, we refer to \cite{OP}. 
\end{rem}

\begin{ex} \label{Hyperflaeche} \rm
Suppose that $X$ is a hypersurface in $T$. Then $I_X$ is generated by some $f\in K[M)\setminus \{0\}$. For $t \in T(K)$, the initial ideal $\In_t(X)$ is generated by the initial form $\In_t(f)$. This is clear from $\In_t(fg)=\In_t(f)\In_t(g)$ (up to multiplication by $K^\times$) for every $g \in K[M]$. 
\end{ex}

\section{Affine toric schemes over a valuation ring}

First, we recall some facts from the theory of normal toric varieties which will be  very important in the sequel. We refer to \cite{CLS}, \cite{Fu2}, \cite{KKMS} or \cite{Oda} for details. They are independent of any valuations on the field $K$. Then 
we assume that $K$ is endowed with a  non-archimedean absolute value $|\phantom{a}|$ with valuation $v:=-\log|\phantom{a}|$ and value group $\Gamma:=v(K^\times)$.  
We consider the split torus $\tdop = \Spec(K^\circ[M])$ over the valuation ring $K^\circ$ with generic fibre $T$. The main focus will be laid on the theory of affine  $\tdop$-toric schemes over $K^\circ$ associated to a pointed $\Gamma$-rational polyhedron. While the generic fibre of such a scheme is a $T$-toric variety over $K$, the geometry of the special fibre is more complicated and is closely related to the combinatorics of the   polyhedron. This section can be seen as a generalization of \S 4.3 in \cite{KKMS}, where the case of a discrete valuation is handled. Further references: \cite{Rab}, \cite{BPR}.


\begin{Def} \rm \label{toric variety} Let $K$ be a field and let $T$ be a split torus over $K$. 
A $T$-{\it  toric variety}  is a  variety $Y$ over $K$ 
containing $T$ as an open dense subset such that the translation action of $T$ on itself extends to an algebraic $T$-action on $Y$. 
\end{Def}

\begin{art} \rm \label{affine toric variety}
There are bijective correspondences between
\begin{itemize}
\item[(a)] rational polyhedral cones $\sigma$ in $N_\rdop$ which do not contain a line;
\item[(b)] finitely generated saturated semigroups $S$ in $M$ which generate $M$ as a group;
\item[(c)] affine normal $T$-toric varieties $Y$ over $K$ (up to equivariant isomorphisms restricting to the identity on $T$).
\end{itemize}
The correspondences are given by $S=\check{\sigma}\cap M$ and $Y=\Spec(K[S])$. We refer the reader to the appendix for the terminology from convex geometry.
\end{art}

\begin{art} \rm \label{toric varieties and fans}  
In general, there is a bijective correspondence between normal $T$-toric varieties $Y$ over $K$ (up to equivariant isomorphisms restricting to the identity on $T$) and pointed rational fans in $N_\rdop$. We denote the toric variety associated to the fan $\Sigma$ by $Y_\Sigma$. Every cone $\sigma$ of $\Sigma$ induces an open affine toric subset $U_\sigma$ of $Y_\Sigma$ by the affine case considered above and $Y_\Sigma$ is covered by such affine charts. 
\end{art}



We extend the above definition to the case of valuation rings:

\begin{Def} \rm \label{toric scheme}
A $\tdop$-{\it toric scheme} over  the valuation ring $K^\circ$ is an integral separated flat scheme $\Ycal$  over $K^\circ$ such that the generic fibre $\Ycal_\eta$ contains $T$ as an open  subset and such that the translation action of $T$ on $T$ extends to an algebraic action of $\tdop$ on $\Ycal$ over $K^\circ$. We call it a $\tdop$-{\it toric variety} if $\Ycal$ is of finite type over $K^\circ$.
\end{Def}

\begin{Def} \rm \label{polyhedral algebra}
For a $\Gamma$-rational polyhedron $\Delta$ in $N_\rdop$, we set  
$$K[M]^\Delta:= \{\sum_{u \in M} \alpha_u \chi^u \in K[M] \mid  v(\alpha_u) + \langle u, \omega \rangle \geq 0 \; \forall \omega \in \Delta\}.$$
\end{Def}

\begin{prop}  \label{properties of polyhedral algebras}
We get a $K^\circ$-subalgebra $K[M]^\Delta$ of $K[M]$ which is flat over $K^\circ$. Moreover, $K[M]^\Delta$ is an integral domain with $K[M]^\Delta \otimes_{K^\circ} K =K[\check{\sigma} \cap M]$ and quotient field $K(\check{\sigma} \cap M)=K(\rho^\bot \cap M)$ where $\rho = \sigma \cap (-\sigma)$ is the largest linear subspace contained in the recession cone $\sigma$ of $\Delta$.
\end{prop}

\proof It is easy to show that $K[M]^\Delta$ is a $K^\circ$-subalgebra of $K[M]$ and hence it is an integral domain. In particular, $K[M]^\Delta$ has no $K^\circ$-torsion and hence it  is flat over $K^\circ$ (see Lemma \ref{flat module}). For $\sum_{u \in M} \alpha_u \chi^u \in K[M]^\Delta$, it follows from the Minkowski-Weil theorem (see \ref{recession cone}) that $\alpha_u = 0$ for every $u \not
 \in M \cap \check{\sigma}$ and that $K[M]^\Delta \otimes_{K^\circ} K =K[\check{\sigma} \cap M]$. The last claim is now obvious. \qed 
 
\vspace{2mm}
The algebra $K[M]^\Delta$ was studied by \cite{KKMS} in case of a  discrete valuation and by \cite{BPR} in case of an algebraically closed ground field endowed with a non-trivial complete absolute value. We will see in the following that most of their results hold in our more general setting. If the valuation $v$ is trivial, then $\Delta$ is a rational cone $\sigma$ and the above shows that $K[M]^{\Delta}=K[M \cap \check{\sigma}]$ leading to the classical case in \ref{affine toric variety}.


\begin{prop}  \label{finite generation}
If the value group $\Gamma$ is either discrete or divisible in $\rdop$, then the algebra $K[M]^\Delta$ is of finite presentation  over $K^\circ$. 
\end{prop}

\proof It is enough to prove that $K[M]^\Delta$ is  a finitely generated $K^\circ$-algebra. This follows from the fact that every finitely generated flat algebra over an integral domain is of finite presentation (\cite{RG}, Corollaire 3.4.7). If $\Gamma$ is discrete, then either $v$ is a discrete valuation 
 or $v$ is trivial. The latter is covered by the divisible case. 

If $v$ is a discrete valuation, then we may assume $\Gamma= \zdop$. We consider the closure $\sigma$ of the cone  in $N_\rdop \times \rdop_+$ generated by $\Delta \times \{1\}$. It is a rational polyhedral cone (see \ref{affine toric variety and cone} and ref{admissible cones} for an argument). If $\pi$ is a uniformizing parameter for $K^\circ$, then $K[M]^\Delta$ is generated by $\pi^k \chi^u$ with $(u,k)$ ranging over the semigroup $S_\sigma:=\check{\sigma} \cap (M\times \zdop)$. This semigroup is finitely generated (see \ref{affine toric variety}) and hence we get the claim in the case of a discrete valuation. 

If the value group is divisible in $\rdop$, we argue as follows: We reduce to the case of a pointed $\Gamma$-rational polyhedron by the procedure described in \ref{polyhedral scheme}  below.  Then the same proof as for  Proposition 4.11 in \cite{BPR} works. Indeed, the crucial point in this proof is that the vertices $\omega_1, \dots ,\omega_r$ of $\Delta$ are in $N_\Gamma$ which is always the case for $\Gamma$ divisible in $\rdop$. Then it is shown that $K[M]^\Delta$ is generated by the functions $\alpha_{ij}\chi^{u_{ij}}$, where $(u_{ij})_j$ is a finite set of generators for  $\check{\sigma}_i \cap M$ with $\sigma_i$ equal to the local cone ${\rm LC}_{\omega_i}(\Delta)$ and where $\alpha_{ij} \in K$ with $v(\alpha_{ij})+\langle u_{ij}, \omega_i \rangle =0$. \qed

\begin{art} \label{weights} \rm
For $\omega \in N_\rdop$, we will  use the {\it $\omega$-weight}
$$v_\omega(\sum_u \alpha_u \chi^u):=\min_u v(\alpha_u) + \langle u, \omega \rangle$$
on $K[M]$ which extends obviously to a valuation on the field $K(T)$. We may view it as a weighted Gauss-valuation similarly as in  \ref{Tate algebra}. 

For a $\Gamma$-rational polyhedron $\Delta$ in $N_\rdop$, we have $K[M]^\Delta \otimes_{K^\circ} K =K[\check{\sigma} \cap M]$ (see Proposition \ref{properties of polyhedral algebras}). 
It is clear that $v_\Delta := \min_{\omega \in \Delta} v_\omega$ is not necessarily a valuation on $K[\check{\sigma} \cap M]$. However, $\| \phantom{a}\|_\Delta := \exp(-v_\Delta)$ is a {\it $K$-algebra norm} on $K[\check{\sigma} \cap M]$, i.e. we have $\|f \cdot g \|_\Delta \leq \|f\|_\Delta \cdot \| g \|_\Delta$ and $\|\lambda \cdot f\|_\Delta =|\lambda| \cdot \|f\|_\Delta$ for $\lambda \in K$ and $f,g  \in K[\check{\sigma} \cap M]$. 

If $\Delta$ is pointed, then any affine form on $\Delta$ which is bounded below takes its minimum in a vertex and so we have 
\begin{equation} \label{minimum in a vertex}
v_\Delta(f)= \inf_{\omega \in \Delta} v_\omega(f) = \min\{v_\omega(f) \mid \text{$\omega$ vertex of $\Delta$}\}
\end{equation}
for every $f \in K[\check{\sigma} \cap M]$.
\end{art}

If the value group $\Gamma$ is neither discrete nor divisible in $\rdop$, then the $K^{\circ}$-algebra $K[M]^\Delta$ is not necessarily of finite presentation over $K^\circ$ as one can verify in the example $K[x]^\omega$ for any $\omega \in \rdop \setminus \Gamma$ with a non-zero $n \in \ndop$ such that $n \omega \in \Gamma$. The referee has pointed out (and proved) that the following precise equivalence is valid.

\begin{prop} \label{finite presentation in non-discrete case}
Let us assume that the value group $\Gamma$ is not discrete in $\rdop$ and let $\Delta$ be a pointed  $\Gamma$-rational polyhedron in $N_\rdop$. Then the following are equivalent:
\begin{itemize}
 \item[(a)] the $K^\circ$-algebra $K[M]^\Delta$ is of finite presentation;
 \item[(b)] the $K^\circ$-algebra $K[M]^\Delta$ is finitely generated;
 \item[(c)] the vertices of $\Delta$ are in $N_\Gamma$;
 \item[(d)] $v_\Delta(K[M]^\Delta \setminus \{0\}) \subset \Gamma$.
\end{itemize}
\end{prop}

\proof 
The equivalence of (a) and (b) follows from \cite{RG}, Corollaire 3.4.7. 
Since $\Delta$ is a pointed $\Gamma$-rational polyhedron, it is clear that $M \cap \check{\sigma}$ generates $M$ as a group. This and \eqref{minimum in a vertex} prove immediately the equivalence of (c) and (d). The same arguments as in the proof  of Proposition \ref{finite generation} show that (c) implies (b).  

Finally, we prove that (b) yields (d). Let $A$ be the completion of $K[M]^\Delta$ with respect to the algebra norm $\| \phantom{a}\|_\Delta$ from \ref{weights}. Since $A$ has no $(K_v)^\circ$-torsion for the completion $K_v$ of $K$, we have $A \subset \Acal:= A \otimes_{(K_v)^\circ} K_v$ by using Lemma \ref{flat module}. There is a canonical algebra norm on $\Acal$ which extends the norm of the completion $A$. We denote this non-archimedean norm again by $\| \phantom{a}\|_\Delta$. By construction, $A$ is the closed unit ball in $\Acal$ with respect to $\| \phantom{a}\|_\Delta$. 

For a nonzero $\rho$ in the maximal ideal $K^{\circ \circ}$, it is easy to see that the $\rho$-adic completion of $K[M]^\Delta$ is equal to $A$. Since $A$ is of finite presentation over $K^\circ$, we see that $A$ is a topologically finitely generated $(K_v)^\circ$-algebra and hence $\Acal$ is an affinoid $K$-algebra (see \ref{comparision of reduction}). By construction, $\| \phantom{a}\|_\Delta$ is a complete power-multiplicative algebra norm on $\Acal$. By \cite{BGR}, Theorem 6.2.4/1, we easily deduce that $\| \phantom{a}\|_\Delta$ is the spectral radius of $\Acal$ (i.e. the supremum norm) and hence $A =\Acal^\circ$ in the notation of \ref{spectral radius}. Since $A$ is topologically finitely generated and using also that $v$ is not discrete,  we deduce from Corollary 6.4.3/6 of \cite{BGR} that the affinoid algebra $\Acal$ is distinguished. This means that $\Acal$ is a quotient of a Tate algebra such that the spectral radius on $\Acal$ agrees with the quotient norm. In particular, we get $\| \Acal\|_\Delta=|K|$ and hence (d) holds. \qed

\begin{prop} \label{polyhedral algebra is integrally closed}
For a $\Gamma$-rational polyhedron $\Delta$ in $N_\rdop$, the algebra $K[M]^\Delta$ is integrally closed.
\end{prop}

\proof Using the procedure described in \ref{polyhedral scheme} below, we may assume that $\Delta$ is pointed. Since every affine  form on $\Delta$ which is bounded below takes its minimum in a vertex, we deduce that $K[M]^\Delta = \bigcap_\omega K[M]^\omega$ with $\omega$ ranging over the vertices of $\Delta$. Hence it is enough to show that $K[M]^\omega$ is integrally closed in $K[M]$. Since the $\omega$-weight $v_\omega$ from \ref{weights} is a valuation and $K[M]^\omega=\{f \in K[M] \mid v_\omega(f) \geq 0\}$, the same argument as in the case of valuation rings proves the claim. Indeed, let $f^m + a_{m-1}f^{m-1} + \dots + a_0=0$ with  $f \in K[M]$ and all $a_i \in K[M]^\omega$. Then the ultrametric triangle inequality applied to $f^m= - a_{m-1}f^{m-1} - \dots - a_0$ and $v_\omega(a_i) \geq 0$ yield $v_\omega(f) \geq 0$. \qed

\begin{art} \rm \label{polyhedral scheme} \label{torus action on polyhedral scheme}
Let $\Delta$ be a $\Gamma$-rational polyhedron in $N_\rdop$ with recession cone $\sigma$. We call $\Ucal_\Delta:=\Spec(K[M]^\Delta)$ the {\it polyhedral scheme over $K^\circ$ associated to $\Delta$}. By Proposition \ref{properties of polyhedral algebras} and Proposition \ref{polyhedral algebra is integrally closed}, $\Ucal_\Delta$ is a normal scheme which is flat over $K^\circ$. If the value group is discrete or divisible then Proposition \ref{finite generation} shows that $\Ucal_\Delta$ is of finite type over $K^\circ$. 

The $K^\circ$-algebra $K[M]^\Delta$ is $M$-graded and hence $\tdop$ acts on $\Ucal_\Delta$. 
It follows from Proposition \ref{properties of polyhedral algebras} that $\Ucal_\Delta$ is a $\tdop$-toric scheme over $K^\circ$ if and only if $\Delta$ is a pointed polyhedron in the sense of \ref{pointed polyhedron}. In this case, the generic fibre is the affine $T$-toric variety $U_\sigma$ associated to  $\sigma$ (see Proposition \ref{properties of polyhedral algebras}). 

In general, we consider the smallest linear subspace $\rho=\sigma \cap (-\sigma)$ contained in  $\sigma$. Then $\Ucal_\Delta$ is a toric scheme over $K^\circ$ with respect to the split torus $\Spec(K^\circ[M(\rho)])$. Here, we have used the lattice $M(\rho):=M \cap \rho^\bot$ with dual lattice $N(\rho)=N/N_\rho$ where $N_\rho$ is the group $N \cap \rho$. The image of $\Delta$ in $N(\rho)_\rdop$ is a pointed polyhedron whose associated polyhedral scheme is $\Ucal_\Delta$. This procedure is often used to reduce to the case of pointed polyhedra. 
\end{art}

\begin{prop}  \label{open immersions of polyhedral schemes}
Let $\Delta'$ be a closed face of the $\Gamma$-rational polyhedron $ \Delta$  in $N_\rdop$. Then the canonical $\tdop$-equivariant morphism $\Ucal_{\Delta'} \rightarrow \Ucal_\Delta$ is a distinguished open immersion.
\end{prop}

\proof There is a halfspace $\{\omega \in N_\rdop \mid \langle u, \omega \rangle + v(\alpha) \geq 0\}$ containing $\Delta$ such that
the face $\Delta'$ is cut out from $\Delta$ by the hyperplane $\{\omega \in N_\rdop \mid \langle u, \omega \rangle + v(\alpha) = 0\}$ for suitable $u \in M$ and $\alpha \in K$ . We claim that $\Ucal_{\Delta'}$ is the complement  of the closed subscheme of $\Ucal_\Delta$ given by the equation $\alpha \chi^u =0$. To see this, we will show that $K[M]^{\Delta'}$ is the localization $(K[M]^{\Delta})_f$ for $f:= \alpha \chi^u$. Using \ref{polyhedral scheme}, we may assume that $\Delta$ is pointed. By construction, $f$ is in $K[M]^\Delta \subset K[M]^{\Delta'}$ and $f$ is invertible in $K[M]^{\Delta'}$. This yields $(K[M]^{\Delta})_f \subset K[M]^{\Delta'}$. To prove the reverse inclusion, it is enough to show that any homogeneous element $f'=\alpha' \chi^{u'} \in K[M]^{\Delta'}$ is contained in $(K[M]^{\Delta})_f$. From our assumptions, we deduce that there is a sufficiently large $m \in \ndop$ such that $$m\left(v(\alpha)+\langle u, \omega \rangle\right) + v(\alpha') + \langle u' , \omega \rangle \geq 0 $$
for all vertices $\omega$ of $\Delta$. By \eqref{minimum in a vertex},  we get $f'f^{m} \in K[M]^{\Delta}$ proving the claim. \qed

\vspace{2mm}

Let $\Delta$ be any $\Gamma$-rational polyhedron in $N_\rdop$. It follows from \ref{polyhedral scheme} that the split torus $\tdop_{\tilde{K}}$ acts on the special fibre of $\Ucal_\Delta$. Our goal is the description of the orbits of this action and hence only the induced reduced structure $((\Ucal_\Delta)_s)_{\rm red}$ is relevant. 

\begin{lem} \label{reduced special fibre and omega-weight}
The reduced induced structure on the special fibre is given by $$((\Ucal_\Delta)_s)_{\rm red}=\Spec(K[M]^\Delta/\{f \in K[M]^\Delta \mid v_\Delta(f)>0\}  ).$$
The special fibre $(\Ucal_\Delta)_s$ is always of finite type over $\tilde{K}$. If the valuation $v$ is not discrete or if $\Delta$ is pointed with all vertices contained in $N_\Gamma$, then $(\Ucal_\Delta)_s$ is reduced.
\end{lem}

\proof If $v$ is the trivial valuation, then the special fibre is equal to the generic fibre and the claims are obvious. So we may assume that $v$ is non-trivial. The special fibre of $\Ucal_\Delta$ is a closed subscheme of $\Ucal_\Delta$ given by the ideal $I=K^{\circ \circ}K[M]^\Delta$ in $K[M]^\Delta$. Since $v_\Delta$ is power-multiplicative, it is clear that the radical ideal $\sqrt{I}$ of $I$ is contained in the ideal $J=\{f \in K[M]^\Delta \mid v_\Delta(f)>0\} $. On the other hand, $J$ is an $M$-homogeneous ideal in $K[M]^\Delta$ and so it is enough to show that every $f=\alpha \chi^u \in J$ is contained in $\sqrt{I}$. Since the valuation $v$ on $K$ is non-trivial, there is $\beta$ in the maximal ideal $K^{\circ \circ}$ of $K$ and $v(\beta) \leq v_\Delta(f^m)$ for $m \in \ndop$ sufficiently large. We conclude that $f^m \in I$ proving $I \subset \sqrt{J}$ and the first claim. 

We handle now the remaining claims first in the case of a discrete valuation. Then  Proposition \ref{finite generation} yields that the special fibre is of finite type. If $\Delta$ is pointed with all vertices contained in $N_\Gamma$, then \eqref{minimum in a vertex} shows that we may choose $\beta$ as the uniformizing parameter and $m=1$ in the above argument. This proves $I=J$ and hence $(\Ucal_\Delta)_s$ is reduced.

It remains to handle the case of a value group $\Gamma$ which is not discrete in $\rdop$. Then $\Gamma$ is dense in $\rdop$ and hence we may choose $m=1$ in the above argument proving again that $I=J$. We conclude that $(\Ucal_\Delta)_s$ is reduced. To prove that the special fibre is of finite type over ${\tilde{K}}$, we may assume that $\Delta$ is pointed using the procedure described in  \ref{torus action on polyhedral scheme}, . 
Since $\Delta$ is $\Gamma$-rational, we have a non-zero $m \in \ndop$ such that $m \omega \in N_\Gamma$ for every vertex $\omega$ of $\Delta$. We conclude that $M_\omega := \{u \in M \mid \langle u, \omega \rangle \in \Gamma \}$ is a subgroup of $M$ of finite index. 
Let $\sigma_\omega$ be the local cone of $\Delta$ at $\omega$ and let $\Delta_\omega:=\omega + \sigma_\omega$. By Proposition \ref{finite presentation in non-discrete case}, 
the $K^\circ$-algebra $K[M_\omega]^{\Delta_\omega}$ is generated by a finite set $S_\omega$. 

To see that $(\Ucal_\Delta)_s$ is of finite type, it is enough to show that $K[M]^\Delta /K^{\circ \circ} K[M]^\Delta$ is generated by the reductions of $S := \bigcup_\omega S_\omega$, where $\omega$ ranges over the  vertices of $\Delta$. Using that $\Gamma$ is not discrete in $\rdop$, we have seen above that $K[M]^\Delta /K^{\circ \circ} K[M]^\Delta= K[M]^\Delta/\{f \in K[M]^\Delta \mid v_\Delta(f)>0\}$. Hence it is enough to show that any $f=\alpha \chi^u$ with $\alpha \in K$, $u \in M$ and $v_\Delta(f)=0$ is in the algebraic span of $S$ over $K^\circ$ modulo the ideal $J=\{f \in K[M]^\Delta \mid v_\Delta(f)>0\}$. By \eqref{minimum in a vertex}, there is a vertex $\omega$ of $\Delta$ such that $v(\alpha)+\langle u, \omega \rangle = v_\omega(f)=v_\Delta(f)=0$. This yields $u \in M_\omega$ and $f \in K[M_\omega]^{\Delta_\omega}$. We conclude that $f$ is in the algebraic span of $S_\omega$ over $K^\circ$ proving the claim.
\qed

\begin{prop} \label{irreducible components of the special fibre}
Let $\Delta$ be a pointed $\Gamma$-rational polyhedron in $N_\rdop$. Then there is a bijection between the vertices  of  $\Delta$ and the irreducible components of $(\Ucal_\Delta)_s$. The irreducible component corresponding to the vertex $\omega$ is the closed subscheme $Y_\omega$ of $\Ucal_\Delta$ given by the prime ideal $\{f \in K[M]^\Delta \mid v_\omega(f)>0\}$ of $K[M]^\Delta$. 
\end{prop}

\proof 
Since $v_\omega$ is a valuation on $K[M]^\Delta$ for any $\omega \in \Delta$, it is clear that $I_\omega:=\{f \in K[M]^\Delta \mid v_\omega(f)>0\}$ is a prime ideal in $K[M]^\Delta$. Since $\Delta$ is a pointed polyhedron, the restriction  of any affine form $v_\omega(\alpha \chi^u)$ to $\Delta$ with $\alpha \chi^u \in K[M]^\Delta$ takes its minimum in a vertex $\omega$ and for every vertex, there is such an affine  form which has its minimum precisely in this vertex. This means that the set of prime ideals $I_\omega$, with $\omega$ ranging over the vertices of $\Delta$, is a minimal primary decomposition of the ideal $\{f \in K[M]^\Delta \mid v_\Delta(f)>0\}$. We have seen in Lemma \ref{reduced special fibre and omega-weight} that the latter is the ideal of the reduced scheme underlying the special fibre $(\Ucal_\Delta)_s$ in $K[M]^\Delta$ and hence the $I_\omega$ are the ideals of the irreducible components of $(\Ucal_\Delta)_s$. \qed

\begin{cor} \label{irreducible component is toric}
The irreducible component $Y_\omega$ of $(\Ucal_\Delta)_s$ is naturally $\tdop_{\tilde{K}}$-equiva\-ri\-antly isomorphic to $((\Ucal_{\Delta(\omega)})_s)_{\rm red}$ where $\Delta(\omega)= \omega + {\rm LC}_\omega(\Delta)$. Moreover, $M_\omega :=\{u \in M \mid \langle u , \omega \rangle \in \Gamma \}$ is a sublattice of finite index in $M$ and $Y_\omega$ is equivariantly (but non-canonically) isomorphic to the $\Spec(\tilde{K}[M_\omega])$-toric variety over $\tilde{K}$ associated to the local cone ${\rm LC}_\omega(\Delta)$.
\end{cor} 

\proof Since $\Delta \subset \Delta(\omega)$, we have a canonical injective homomorphism 
$$\varphi:  K[M]^{\Delta(\omega)}/\{f \in K[M]^{\Delta(\omega)} \mid v_\omega(f)>0\}      \rightarrow             K[M]^\Delta/\{f \in K[M]^\Delta \mid v_\omega(f)>0\}         .$$
To show surjectivity, it is enough to show that the residue class of $f=\alpha \chi^u \in K[M]^{\Delta}$ is in the image of $\varphi$. We may assume that $v_\omega(f)=0$ otherwise this is trivial. Then the affine form $\Delta \rightarrow \rdop, \, \nu \mapsto v_\nu(f)$ takes its minimum in the vertex $\omega$. This even holds if we extend the affine  form to $\Delta(\omega)$ by definition of the local cone ${\rm LC}_\omega(\Delta)$. We conclude that $f \in K[M]^{\Delta(\omega)}$ proving that $\varphi$ is an isomorphism. By Lemma \ref{reduced special fibre and omega-weight} and Proposition \ref{irreducible components of the special fibre}, we deduce $Y_\omega \cong ((\Ucal_{\Delta(\omega)})_s)_{\rm red}$. Equivariance of this isomorphism follows from the fact that $\varphi$ is an $M$-graded isomorphism.

Since $\Delta$ is $\Gamma$-rational, there is a non-zero $m \in \ndop$ with $m \omega \in N_\Gamma$ and hence $M_\omega$ is a sublattice of finite index in $M$. It is trivial to show that the canonical homomorphism 
from $K[M_\omega]^{\Delta(\omega)}/\{f \in K[M_\omega]^{\Delta(\omega)} \mid v_\omega(f)>0\} $  to $ K[M]^{\Delta(\omega)}/\{f \in K[M]^{\Delta(\omega)} \mid v_\omega(f)>0\}       $
is an isomorphism. We conclude that we may replace $M$ by $M_\omega$ and so we may assume $M=M_\omega$. Then there is $t \in T(K)$ with $\trop_v(t)=\omega$. We may replace $\Delta$ by $\Delta - \omega$  which means geometrically that we use translation by $t^{-1}$ on $T$. Then $\omega = 0$ is the given vertex of $\Delta$. By the first claim, the irreducible component $Y_\omega$ is equivariantly isomorphic to $$\Spec(K[M]^{\Delta(\omega)}/\{f \in K[M]^{\Delta(\omega)} \mid v_\omega(f)>0\})= \Spec(\tilde{K}[M]^{{\rm LC}_\omega(\Delta)})$$
which is the $\tdop_{\tilde{K}}$-toric variety associated to ${\rm LC}_\omega(\Delta)$. 
\qed

\begin{art} \rm \label{completion of an affine toric scheme}
Next, we describe the reduction map with respect to the $\tdop$-toric scheme $\Ucal_\Delta$ over $K^\circ$ associated to the pointed $\Gamma$-rational polyhedron $\Delta$ in $N_\rdop$. Recall from \ref{polyhedral scheme} that the $T$-toric variety $U_\sigma$ is the generic fibre of $\Ucal_\Delta$  where $\sigma$ is the recession cone of $\Delta$. We have seen in \ref{reduction for Xan} that the reduction is a map to the special fibre $(\Ucal_\Delta)_s$ which is defined on the set $U_\sigma^\circ:=\{p \in U_\sigma^{\rm an} \mid p(f) \leq 1 \; \forall f \in K[M]^\Delta\}$. The points of  $U_\sigma^\circ$ are induced by the potentially integral points of $\Ucal_\Delta$. 
\end{art}

We will describe the analytic structure of $U_\sigma^\circ$ using the following result of Joe Rabinoff.

\begin{prop} \label{polyhedral subdomain}
We assume that the valuation $v$ on $K$ is  complete. Let $\Delta$ be a pointed $\Gamma$-rational polyhedron in $N_\rdop$ with recession cone $\sigma$ and let $\|\phantom{a}\|$ be any norm on $M_\rdop$. Then the set of Laurent series
$$\Acal_\Delta:=\left \{ \sum_{u \in \check{\sigma} \cap M} a_u \chi^u \mid \lim_{\|u\| \to \infty} v(a_u) + \langle u , \omega\rangle  =\infty \, \; \forall \omega \in \Delta \right\}$$
is a strictly affinoid algebra  with spectral radius
\begin{equation} \label{polytopal supnorm}
\rho(\sum_{u \in \check{\sigma} \cap M} a_u \chi^u )= \sup_{\omega \in \Delta, \, u \in \check{\sigma} \cap M} |a_u|e^{-\langle u, \omega \rangle} = \max_{\text{$\omega$ \rm vertex, $u \in \check{\sigma} \cap M$}}|a_u|e^{-\langle u, \omega \rangle}.
\end{equation}
\end{prop}

\proof In the case of a non-trivial valuation, we  use \cite{Rab}, Proposition 6.9. If $v$ is trivial, then the sums in the definition of $\Acal_\Delta$ are finite and hence $\Acal_\Delta=K[M]^\Delta$ which immediately yields the claims.  \qed

\begin{rem} \rm \label{remarks to polytopes and regularity}
The special case of polytopal domains  was studied in \cite{Gu3}.  Using Hochster's theorem for toric varieties,   Rabinoff has shown  that $\Acal_\Delta$ is Cohen-Macauley for any $\Gamma$-rational polyhedron $\Delta$ (see \cite{Rab}, Proposition 6.9). If the valuation is discrete or $K$ algebraically closed, then Wilke \cite{Wil} has shown that $\Acal_\Delta$ is a factorial ring for $\Gamma$-rational polytopes $\Delta$. 

\end{rem}

\begin{prop} \label{analytic structure of integral points}
Using the notation from Proposition \ref{polyhedral subdomain}, the Berkovich spectrum $\Mcal(\Acal_\Delta)$ is an affinoid subdomain of $U_\sigma^{\rm an}$ which is equal to $U_\sigma^\circ$. Moreover, the special fibres of $\Mcal(\Acal_\Delta)$ and $\Ucal_\Delta$ agree which means
$ \Spec(\Acal_\Delta^\circ /(K^{\circ \circ} \Acal_\Delta^\circ)) = (\Ucal_\Delta)_s$.
\end{prop}

\proof If $v$ is trivial, then $\Delta$ is a rational pointed cone and we have seen $\Acal_\Delta^\circ=\Acal_\Delta=K[M]^\Delta$ which makes the claims obvious. So we may assume that $v$ is non-trivial. By Proposition \ref{polyhedral subdomain}, $\Ocal(U_\sigma)=K[ \check{\sigma} \cap M]$ is dense in $\Acal_\Delta$ and hence $\Mcal(\Acal_\Delta)$ may be seen as a subset of $U_\sigma^{\rm an}$. In fact, it is shown in \cite{Rab}, Proposition \ref{polyhedral subdomain}, that  $\Mcal(\Acal_\Delta)$ is an affinoid subdomain of $U_\sigma^{\rm an}$. 
Moreover, we deduce from Rabinoff's result that $K[M]^\Delta$ is a subset of $\Acal_\Delta^\circ$ and hence $\Mcal(\Acal_\Delta) \subset U_\sigma^\circ$. We will prove the reverse inclusion (which would follow immediately from \ref{comparision of reduction} if $K[M]^\Delta$ is finitely generated) and so we choose $p \in U_\sigma^\circ$. 

We claim first that $p(\chi^u) \leq \rho(\chi^u)$ for any $u \in \check{\sigma} \cap M$ where $\rho(\chi^u)$ is the spectral radius in $\Acal_\Delta$. There is a vertex $\omega_0$ of $\Delta$ such that the halfspace $H:=\{\omega \in N_\rdop \mid \langle u, \omega \rangle \geq 0\} + \omega_0$ contains $\Delta$. By $\Gamma$-rationality of $\Delta$, there is a non-zero $m\in  \ndop$ such that $m \omega_0 \in N_\Gamma$. We conclude that there is a non-zero $\alpha \in K$ such that $v(\alpha)+\langle mu, \omega_0 \rangle=0$. Using $H= \{\omega \in N_\rdop \mid v(\alpha) + \langle mu, \omega \rangle \geq 0\}$, we get $\alpha \chi^{mu} \in K[M]^\Delta$ and $\rho(\alpha \chi^{mu})=|\alpha|e^{-\langle mu,  \omega \rangle} =1$ follows from Proposition \ref{polyhedral subdomain}. Using power multiplicativity of both  $p \in U_\sigma^\circ$ and $\rho$, we get
$$|\alpha| p(\chi^u)^m=p(\alpha \chi^{mu})\leq 1 = \rho(\alpha \chi^{mu})=|\alpha| \rho(\chi^u)^m.$$
This proves $p(\chi^u) \leq \rho(\chi^u)$ for any $u \in \check{\sigma} \cap M$. 

Next, we will prove $p(f) \leq \rho(f)$ for any $f \in \Ocal(U_\sigma)$. Note that $f=\sum_u \alpha_u \chi^u$ where $u$ ranges over a finite subset of $ \check{\sigma} \cap M$. Using the above and Proposition \ref{polyhedral subdomain}, we get
$$p(f) \leq \max_u |\alpha_u| p(\chi^u )\leq \max_u |\alpha_u| \rho( \chi^u) = \max_u \rho(\alpha_u \chi^u)=\rho(f)$$
as desired. Now $p \leq \rho$ yields that $p$ extends uniquely to a multiplicative seminorm of $\Mcal(\Acal_\Delta)$. This proves $\Mcal(\Acal_\Delta) = U_\sigma^\circ$.  The claim about the special fibres  follows immediately from Proposition \ref{polyhedral subdomain}. \qed

\vspace{2mm}

In the following, we do not necessarily assume that the valuation $v$ on $K$ is complete as the analytifications are defined on the completion of $K$ anyway.

\begin{cor} \label{surjectivity of reduction}
Let $\Delta$ be any pointed $\Gamma$-rational polyhedron in $N_\rdop$ with recession cone $\sigma$. Then 
the reduction map from \ref{reduction for Xan}  maps $U_\sigma^\circ \cap \Tan$ surjectively onto $(\Ucal_\Delta)_s$. 
\end{cor}

\proof Using the procedure described in \ref{polyhedral scheme}, we may assume that $\Delta$ is a  pointed polyhedron. Passing to the completion does not change the special fibre and so we may assume $K$ complete. By Proposition \ref{analytic structure of integral points}, the special fibre of $\Ucal_\Delta$ agrees with the special fibre of $U_\sigma^\circ = \Mcal(\Acal_\Delta)$. Since $T$ is the dense orbit in the generic fibre $U_\sigma$, it is clear that $\Tan \cap U_\sigma^\circ$ is Zariski open and dense in the affinoid subdomain $U_\sigma^\circ$. Now the claim follows from Lemma \ref{lifting lemma}.  \qed

\begin{lem} \label{integral points and tropicalization}
We have $U_\sigma^\circ \cap \Tan = \trop_v^{-1}(\Delta)$ for a pointed $\Gamma$-rational polyhedron $\Delta$.
\end{lem}

\proof By definition, $U_\sigma^\circ \cap \Tan$ is the set of multiplicative seminorms $p$ on $K[M]$ extending $|\phantom{a}|$ with the additional condition that $p(f) \leq 1$ for every $f \in K[M]^\Delta$. The latter is equivalent to the condition that if there are $u \in M$, $\alpha \in K$ with  $v(\alpha) + \langle u , \omega \rangle \geq 0$ for every $\omega \in \Delta$, then $p(\alpha \chi^u) \leq 1$. Since $\Delta$ is defined by such inequalities and since $- \log p(\chi^u)= \langle u,\trop_v(p) \rangle$, this is also equivalent to $\trop_v(p) \in \Delta$. This proves the claim.

If $v$ is non-trivial, then we could  use the above Proposition \ref{analytic structure of integral points} and Proposition 6.9 from \cite{Rab} to get an alternative proof.
\qed

\vspace{2mm}

As we have defined the tropicalization map  only on $\Tan$, we restrict the reduction map in the following proposition to $U_\sigma^\circ \cap \Tan$. By abuse of notation, we denote this restriction $U_\sigma^\circ \cap \Tan \rightarrow (\Ucal_\Delta)_s$ also by $\pi$. In the following result, we use the partial order on the set of orbits (resp. open faces) given by inclusion of closures.

\begin{prop} \label{affine orbit correspondence}
Let $\Delta$ be a pointed $\Gamma$-rational polyhedron in $N_\rdop$ and let $\Ucal_\Delta$ be the associated $\tdop$-toric scheme over $K^\circ$. Then there is a bijective  order reversing correspondence between $\tdop$-orbits $Z$ of $(\Ucal_\Delta)_s$ and open faces $\tau$ of $\Delta$ given by
$$Z=\pi(\trop_v^{-1}(\tau)), \quad \tau = \trop_v(\pi^{-1}(Z)).$$
Moreover, we have $\dim(Z)+ \dim(\tau)=n$. \end{prop}

\proof  Let $\tau$ be an open face of $\Delta$. The affinoid torus $\{p \in \Tan \mid p(\chi^u)=1 \; \forall u \in M\}$ operates on $\trop_v^{-1}(\tau)$. By passing to the reductions, we see that $Z:=\pi(\trop_v^{-1}(\tau))$ is $\tdop_{\tilde{K}}$-invariant. Note that $Z$ is well-defined by Lemma \ref{integral points and tropicalization}. It is clear that distinguished open faces give rise to distinguished $\tdop_{\tilde{K}}$-invariant subsets of $(\Ucal_\Delta)_s$. It remains to show that $Z$ is an orbit. We are allowed to pass to a $\tdop$-invariant open subset and hence we may assume that $\tau = \relint(\Delta)$ by using Proposition \ref{open immersions of polyhedral schemes}. 
 Since $\Delta$ is pointed, it has a vertex $\omega$. 
We have seen in Corollary \ref{irreducible component is toric} that the irreducible component $Y_\omega$ is non-canonically isomorphic to the  $\Spec(\tilde{K}[M_\omega])$-toric variety over $\tilde{K}$ associated to the local cone ${\rm LC}_\omega(\Delta)$. 
We claim that $Z$ is the unique closed orbit $Z'$ of $Y_\omega$. 
Since $Z$ is invariant and $Z'$ is an orbit, it is enough to show that  $\pi(p) \in Z'$ for every $p \in \trop_v^{-1}(\tau)$. 

We will first prove that $\pi(p) \in Y_\omega$. By Proposition \ref{irreducible components of the special fibre}, the latter is given by the $M$-graded prime ideal $\{f \in K[M]^\Delta \mid v_\omega(f)>0\}$. So let us choose $f=\alpha \chi^u \in K[M]^\Delta$ with $v_\omega(f)>0$. Then $v_\nu(f)>0$ for all $\nu \in \tau = \relint(\Delta)$. In particular, this holds for $\nu = \trop_v(p)$ and hence $-\log p(f)=v_\nu(f)>0$. We conclude $p(f)<1$ which means that $f$ is contained in the prime ideal of $\pi(p)$ in $K[M]^\Delta$. This proves $\pi(p) \in Y_\omega$.
 
The well-known orbit-cone correspondence for toric varities over a field shows that the closed orbit $Z'$ of $Y_\omega$ is given by the ideal in $\tilde{K}[Y_\omega]$ generated by $\{\chi^u \mid u \in M_\omega, \, \langle u , \omega' \rangle > 0 \; \forall \omega' \in \tau' \}$ where $\tau':={\rm LC}_\omega(\tau)$. Taking into account how $Y_\omega$ is defined as a $\Spec(\tilde{K}[M_\omega])$-toric variety, we conclude that $Z'$ is given as a closed subscheme of $\Ucal_\Delta$ by the ideal generated by $\{f=\beta \chi^u \mid \beta \in K, \, u \in M_\omega, \, v_{\omega'}(f)>0 \; \forall \omega' \in \omega+\tau'\}$. For such an $f$, we conclude   $- \log p(f) =v_\nu(f)>0$ as above and hence $p(f)<1$. Again this means  $\pi(p) \in Z'$ proving that $Z=Z'$. We conclude that $Z=\pi(\trop_v^{-1}(\tau))$ is a $\tdop$-orbit.

Since $\pi$ maps $U_\sigma^\circ  \cap \Tan = \trop_v^{-1}(\Delta)$ onto $(\Ucal_\Delta)_s$ by Corollary \ref{surjectivity of reduction}, we get a bijective  correspondence between open faces of $\Delta$ and torus orbits of $(\Ucal_\Delta)_s$. Since $Z$ is the torus orbit of $Y_\omega$ corresponding to the open cone $\tau'$, we get $\dim(\tau)+\dim(Z)=\dim(\tau')+\dim(Z)=n$ from the theory of toric varieties over a field. Moreover, we conclude from the special case $Y_\omega$ that the correspondence is order reversing.

Conversely, let $Z$ be any torus orbit of $(\Ucal_\Delta)_s$. By the above, we have $Z=\pi(\trop_v^{-1}(\tau))$ for a unique open face $\tau$ of $\Delta$. This yields $\tau \subset \trop_v(\pi^{-1}(Z))$. Equality follows from the fact that the torus orbits (resp. open faces) form a partition of $(\Ucal_\Delta)_s$ (resp. $\Delta$). 
\qed 

\begin{rem} \rm \label{remarks to the affine orbit correspondence}
The bijective correspondence between open faces and orbits holds more generally for the polyhedral scheme $\Ucal_\Delta$ associated to any $\Gamma$-rational polyhedron $\Delta$ in $N_\rdop$. This follows from the reduction to the case of pointed polyhedra described in \ref{polyhedral scheme}. 

If $\Delta'$ is any $\Gamma$-rational polyhedron contained in $\Delta$, then the canonical equivariant morphism $\Ucal_{\Delta'} \rightarrow \Ucal_\Delta$ is an open immersion if and only if $\Delta'$ is a closed face of $\Delta$. We have seen one direction in Proposition \ref{open immersions of polyhedral schemes} and the converse follows easily from the orbit-face correspondence.

If $v$ is trivial, then $\Delta$ is a pointed rational cone and the  arguments in the proof of Proposition \ref{affine orbit correspondence} show that we get the classical orbit-cone correspondence for toric varieties over a field  from \cite{Fu2}, \S 3.1, or from \cite{CLS}, \S3.2. 
\end{rem}

We consider a field extension $L/K$ and an arbitrary valuation $w$ on $L$ (not necessarily of rank $1$) extending $v$ with valuation ring $L^\circ$ and value group $\Gamma_L$. For $P \in T(L)$, we define $\trop_w(P) \in N_{\Gamma_L}=\Hom(M, \Gamma_L)$ by $u \mapsto w(\chi^u(P))$ similarly as in the case of a valuation of rank $1$. We recall from \ref{compatibility properties}
 that a $\Gamma$-rational polyhedron $\Delta$ in $N_\rdop$ induces a canonical polyhedron $\Delta(\Gamma_L)$ in $N_{\Gamma_L}$. Then we have the following analogue of Lemma \ref{integral points and tropicalization} for $L^\circ$-integral points:

\begin{prop} \label{integral points and tropicalization revisted}
Under the assumptions above, let $\Delta$ be a pointed $\Gamma$-rational polyhedron in $N_\rdop$ and let $P \in T(L)$. Then $P$ is an $L^\circ$-integral point of $\Ucal_\Delta$ if and only if $\trop_w(P) \in \Delta(\Gamma_L)$.
\end{prop}

\proof The point $P$ is $L^\circ$-integral if and only if $ v(\alpha) + w(\chi^u(P)) \geq 0$ for all $\alpha \chi^u \in K[M]^\Delta$. 
We note that $w(\chi^u(P))=\langle u , \trop_w(P) \rangle$. Since $\Delta$ is a $\Gamma$-rational polyhedron, $\Delta$ is the intersection of the half-spaces $H_{u, \alpha}^+:=\{\omega \in N_\rdop \mid v(\alpha) + \langle u, \omega \rangle \geq 0\}$ with $\alpha \chi^u$ ranging over $K[M]^\Delta$. Moreover, $\Delta$ is the intersection of finitely many such half-spaces. It follows from \ref{polyhedra for totally ordered groups} that $\Delta(\Gamma_L)$ is the intersection of the sets $H^+_{u, \alpha}(\Gamma_L)=\{\omega \in N_{\Gamma_L} \mid v(\alpha)+\langle u, \omega \rangle \geq 0\}$ and hence we get the claim. \qed


\section{Toric schemes over a valuation ring} 

In this section, $K$ is a field endowed with a  non-archimedean valuation $v$ and  $\Gamma$ is the valuation group of $v$. We extend the theory of toric schemes over a discrete valuation ring from \cite{KKMS} to this more general situation. More precisely, we will use the affine toric schemes associated to pointed polyhedra from the previous section to define toric schemes. For the glueing process, it is necessary to work with fans in $N_\rdop \times \rdop_+$ rather than polyhedral complexes in $N_\rdop$. 
Recall that $\tdop = \Spec(\Tor)$ is the split multiplicative torus over $K^\circ$ with generic fibre $T$. The character group of $T$ is $M$ and $N$ is the dual lattice.
Further references for the special case of a discrete valuation are \cite{BPS} (with a lot of arithmetic applications) and  \cite{Smi} (from the projective point of view).

\begin{art} \rm \label{affine toric variety and cone}
As a building block, we will use the affine $\tdop$-toric scheme $\Ucal_\Delta$ over $K^\circ$ from \ref{polyhedral scheme} for any pointed $\Gamma$-rational polyhedron $\Delta$ in $N_\rdop$. 
For glueing,  it  is better to replace $\Delta$ by the closed cone $\sigma=c(\Delta)$ in $N_\rdop \times \rdop_+$ generated by $\Delta \times \{1\}$.  For $s \in \rdop_+$, let $\sigma_s := \{\omega \in N_\rdop \mid (\omega,s) \in \sigma\}$. For $s >0$, we have $\sigma_s=s \Delta$ and $\sigma_0$ is the recession cone of $\Delta$. This follows easily from $c(\Delta)=\{(\omega,s) \in N_\rdop \times \rdop_+\mid  \langle u_i, \omega \rangle + s c_i \geq 0 \; \forall i\}$ using that  $\Delta$ is the intersection of finitely many halfspaces $\{\omega \in W \mid \langle u_i , \omega \rangle +  c_i \geq 0)\}$ with $u_i \in M$ and  $c_i \in \Gamma$. 

\end{art}


\begin{art} \rm \label{admissible cones}
A cone $\sigma$ in $N_\rdop \times \rdop_+$ is called {\it $\Gamma$-admissible} if it may be written as
$$\sigma = \bigcap_{i=1}^N \{ (\omega, s) \mid \langle u_i, \omega \rangle + s c_i \geq 0\}$$
for $u_1, \dots , u_N \in M$ and $c_1, \dots, c_N \in \Gamma$ and if $\sigma$ does not contain a line.  For $s \in \rdop_+$, we define $\sigma_s$ as above. 

Note that $\Delta \mapsto c(\Delta)$ gives a bijection between the set of non-empty pointed $\Gamma$-rational polyhedra in $N_\rdop$ and the set of $\Gamma$-admissible cones in $N_\rdop \times \rdop_+$ which are not contained in $N_\rdop \times \{0\}$. The inverse map is $\sigma \mapsto \sigma_1$. 
\end{art}

\begin{Def} \rm \label{algebra of admissible cone}
For a $\Gamma$-admissible cone $\sigma$ in $N_\rdop \times \rdop_+$, we define
$$K[M]^\sigma:= \{\sum_{u \in  M} \alpha_u \chi^u \in K[M] \mid  c v(\alpha_u) + \langle u, \omega \rangle \geq 0 \; \forall (\omega,c) \in \sigma\}$$
and $\Vcal_\sigma := \Spec( K[M]^\sigma)$.
\end{Def}

\begin{prop}   \label{torus orbits on toric scheme}
$\Vcal_\sigma$ is an affine normal $\tdop$-toric scheme over $K^\circ$ with generic fibre equal to the affine toric variety $U_{\sigma_0}$ associated to $\sigma_0$. If the value group $\Gamma$ of $K$ is discrete or divisible in $\rdop$, then $\Vcal_\sigma$ is an affine normal $\tdop$-toric variety over $K^\circ$. 
\end{prop}

\proof If $\sigma$ is contained in the hyperplane $N_\rdop \times \{0\}$,  then  $\Vcal_\sigma$ is the normal toric variety $U_{\sigma_0}$ over  $K$ associated to $\sigma_0$. Since $K$ is of finite type over the valuation ring $K^\circ$, it is also a normal toric variety over $K^\circ$.  

If $\sigma$ is not contained in  $N_\rdop \times \{0\}$, then $\sigma_1$ is a non-empty $\Gamma$-rational polyhedron $\Delta$ in $N_\rdop$ with $\Vcal_\sigma = \Ucal_{\Delta}$ and the claim follows from \ref{polyhedral scheme}. \qed


\begin{Def} \rm \label{admissible fans}
A {\it $\Gamma$-admissible fan} $\Sigma$ in $N_\rdop \times \rdop_+$ is a fan of $\Gamma$-admissible cones. For $s \geq 0$, let $\Sigma_s$ be the polyhedral complex  $\{\sigma_s \mid \sigma \in \Sigma\}$ in $N_\rdop$.  
\end{Def}

\begin{rem} \rm \label{Burgos-Sombra}
It was noticed by Burgos and Sombra \cite{BS} that if $\Ccal$ is a $\Gamma$-rational polyhedral complex in $N_\rdop$, then $c(\Ccal):=\{c(\Delta) \mid \Delta \in \Ccal \}$ is not necessarily a fan in $N_\rdop \times \rdop_+$. However if the support of $\Ccal$ is convex, then they prove that $c(\Ccal)$ is a fan. This gives a bijective correspondence between complete $\Gamma$-rational pointed polyhedral complexes of $N_\rdop$ and complete $\Gamma$-admissible fans of  $N_\rdop \times \rdop_+$.
\end{rem}

\begin{art} \rm \label{toric scheme associated to an admissible fan} 
Let $\Sigma$ be a $\Gamma$-admissible fan in $N_\rdop \times \rdop_+$. Then the affine $\tdop$-toric schemes $\Vcal_\sigma$, $\sigma \in \Sigma$, may be glued along the open subschemes $\Vcal_\rho$ from common subfaces $\rho$ to get a normal $\tdop$-toric scheme $\Ycal_\Sigma$ over $K^\circ$.  The generic fibre of $\Ycal_\Sigma$ is the normal $T$-toric variety $Y_{\Sigma_0}$ over $K$ associated to the fan $\Sigma_0$ in $N_\rdop$. This follows all from the affine case except separatedness which we shall prove next:
\end{art}

\begin{lem} \label{separatedness}
The scheme $\Ycal_\Sigma$ is separated over $K^\circ$.
\end{lem}

\proof Let $\sigma:= \sigma' \cap \sigma''$ for $\sigma', \sigma'' \in \Sigma$. We have to show that the canonical morphism  $\Vcal_{\sigma} \rightarrow \Vcal_{\sigma'} \times_{K^\circ} \Vcal_{\sigma''}$ is a closed embedding. To prove that we may assume that $\sigma',\sigma''$ are not contained in $N_\rdop \times \{0\}$ (as the claim is well-known for toric varieties over a field). Then we have pointed $\Gamma$-rational polyhedra $\Delta':=\sigma'_1$,  $\Delta'':=\sigma''_1$ and $\Delta:=\Delta' \cap \Delta''=\sigma_1$ in $N_\rdop$ with $\Ucal_{\Delta'} =\Vcal_{\sigma'}$, $\Ucal_{\Delta''} =\Vcal_{\sigma''}$ and $\Ucal_\Delta =\Vcal_\sigma$. We have to show that $K[M]^\Delta$ is generated by $K[M]^{\Delta'}$ and $K[M]^{\Delta''}$ as a $K^\circ$-algebra. 

Let $G:=\{ \gamma \in \rdop \mid \exists k \in \ndop \setminus \{0\}\, \; k\gamma  \in \Gamma\}$; then the affine subspace of $N_\rdop$ generated by the face of a $\Gamma$-rational polyhedron is an $N_G$-translate of a rational linear subspace. We conclude that there is $u_0 \in M$ and $\omega_0 \in N_G$ such that $\Delta = \Delta' \cap (\omega_0 + \{u_0\}^\bot)$, $\Delta' \subset \omega_0 + \{\omega \in N_\rdop \mid \langle u_0, \omega \rangle  \geq 0\}$ and $\Delta'' \subset \omega_0 + \{\omega \in N_\rdop \mid \langle u_0, \omega \rangle  \leq 0\}$. There is $\alpha_0 \in K$ and a non-zero $k \in \ndop$ with $v(\alpha_0)=k \langle u_0, \omega_0 \rangle$. We consider a vertex $\omega'$ of $\Delta'$. By construction, 
$$v_{\omega'}( (\chi^{ku_0}/\alpha_0)^{m}f)= k m \langle u_0, \omega' - \omega_0 \rangle +v_{\omega'}(f)$$ 
is non-negative for $m \gg 0$. We conclude that  $v_{\Delta'}(g) \geq 0$ for $g:=(\chi^{ku_0}/\alpha_0)^{m}f$. 

Since $$v_{\omega''}(\alpha_0 \chi^{-ku_0})=k\langle u_0, \omega_0 - \omega''\rangle \geq 0$$ 
for every $\omega'' \in \Delta''$, 
we conclude that $\alpha_0 \chi^{-ku_0} \in K[M]^{\Delta''}$. Using $f=(\alpha_0 \chi^{-ku_0})^m g$, we get the claim proving that $\Ycal_\Sigma$ is separated.  \qed

\begin{art} \rm     \label{tropicalization and orbits} \label{consequences of the orbit correspondence}
We have a bijective correspondence between torus orbits  of $\Ycal_\Sigma$ and open faces of $\Sigma$. The torus orbits in the generic fibre correspond to the open faces contained in $N_\rdop \times \{0\}$ using the theory of toric varieties over a field. The torus orbits in the special fibre correspond to the remaining open faces of $\Sigma$ using the fact that the latter are the open faces of the polyhedral complex $\Sigma_1$ in $N_\rdop$. Indeed, the orbits are contained in an affine $\tdop$-toric scheme $\Vcal_\sigma$ for some $\sigma \in \Sigma$ and so we may use Proposition \ref{affine orbit correspondence}. 
We will later describe the orbit correspondence for $\Ycal_\Sigma$ in a neat way (see Proposition \ref{uniform orbit correspondence}).

Note that $\Ycal_\Sigma$ is a noetherian topological space  which follows from the fact that both the generic and the special fibre are noetherian. If $\Xcal$ is a closed irreducible subset of $\Ycal_\Sigma$ with non-empty generic fibre $\Xcal_\eta$, then it follows from Proposition \ref{purity of special fibre} below that $\Xcal_s$ is of pure dimension equal to $\dim(X_\eta)$. This means that the dimension and the codimension of $\Xcal$ can be computed using any maximal chain of closed subset and hence we get $\dim(\Xcal)=\dim(\Xcal_\eta)+1$. Moreover, if $Z$ is a closed subsets of $\Xcal_s$, then we get $\codim(Z,\Xcal)=\codim(Z,\Xcal_s)+1$ (see \cite{OP}, Lemma 4.2.3 for similar arguments in the case of any irreducible flat scheme of finite type over $K^\circ$).

If $Z_\tau$ is the torus orbit corresponding to the  open face $\tau$ of $\Sigma$, then the above and Proposition \ref{affine orbit correspondence} show that 
$$\dim(\tau)=\codim(Z_\tau, \Ycal_\Sigma).$$
In particular, the  $\tdop$-invariant prime divisors on $\Ycal_\Sigma$ correspond to the halflines in $\Sigma$. 
The irreducible components of the special fibre of $\Ycal_\Sigma$ correspond to the halflines of $\Sigma$ not contained in $N_\rdop \times \{0\}$ or in other words to the vertices of $\Sigma_1$. 
\end{art}

\begin{lem} \label{multiplicities in special fibre}
Suppose that the valuation on $K$ has value group $\Gamma = \zdop$ and suppose $v(\pi)=1$ for $\pi \in K$. Let $Y_\sigma$ be an irreducible component of the special fibre of $\Ycal_\Sigma$ corresponding to the halfline $\sigma$ of $\Sigma$. Then the multiplicity of the divisor $Y_\sigma$ in $\Ycal_\Sigma$ is equal to $k$, where $(\omega,k)$ is the primitive generator of the monoid $\sigma \cap(N\times \zdop)$. 
\end{lem}

\proof See \cite{KKMS}, \S 4.3. \qed

\begin{prop} \label{reduced special fibre}
If the valuation on $K$ is discrete, then the following conditions are equivalent for a $\Gamma$-admissible fan $\Sigma$ in $N_\rdop \times \rdop_+$:
\begin{itemize}
\item[(a)] The vertices of $\Sigma_1$ are contained in $N_\Gamma$. 
\item[(b)] The special fibre $(\Ycal_\Sigma)_s$ is reduced.
\item[(c)] $(\Ycal_\Sigma)_s$ is geometrically reduced.
\item[(d)] For all valued fields  $(L,w)$ extending $(K,v)$, the formation of $\Ycal_\Sigma$ is compatible with  base change to $L^\circ$. 
\item[(e)] For all $\Delta \in \Sigma_1$, the  canonical map $K[M]^\Delta \otimes_{K^\circ} \tilde{K} \rightarrow \Acal_\Delta^\circ/\Acal_\Delta^{\circ \circ}$ is an isomorphism, where we refer to Proposition \ref{polyhedral subdomain} for the definition of $\Acal_\Delta$.
\end{itemize}
\end{prop}

\proof The equivalence of (a) and (b) follows from Lemma \ref{multiplicities in special fibre}. Clearly, (c) implies (b).

Now let $\sigma \in \Sigma$ and $\Delta:=\sigma_1$. Suppose that the vertices $\omega_j$ of $\Delta$  are contained in $N_\Gamma$.  In this case, we may use also the last part of the proof of Proposition \ref{finite generation} to get a set of generators of   $K[M]^\Delta$ which depends only on the combinatorics of $\Delta$ and hence it generates also $L[M]^\Delta$. 
This proves 
$$L[M]^\Delta = K[M]^\Delta \otimes_{K^\circ} L^\circ$$
and hence (a) implies (d). 

Now suppose that (d) holds. There is a finite extension $L/K$ such that (a) holds for the value group of $L$. By the equivalence of (a) and (b), we conclude that the special fibre of $(\Ycal_\Sigma)_L$ is reduced and hence  the special fibre of $\Ycal_\Sigma$ is also reduced. We may repeat this for any finite extension of $K$ and hence (d) yields (c). Since the residue algebra $\Acal_\Delta^\circ/\Acal_\Delta^{\circ \circ}$ of a strictly affinoid algebra is always reduced, we see that (e) implies (b).

Finally we show that (a) implies (e). Since the vertices are in $N_\Gamma$, it is easy to see that the kernel of the quotient homomorphism  $K[M]^\Delta\rightarrow K[M]^\Delta/\langle K^{\circ \circ}  \rangle  =K[M]^\Delta \otimes_{K^\circ}\tilde{K}$ is equal to $\{\sum_{u \in M} a_u \chi^u \in K[M] \mid v(a_u)+ \langle u, \omega \rangle >0  \; \forall \omega \in \Delta\}$. By density of $K[M]^\Delta$ in $\Acal_\Delta^\circ$, we deduce (e). \qed

\begin{prop}  \label{divisible case}
If $v$ is not a discrete valuation, then (b), (c) and (e) hold always. Moreover, properties (a) and (d) are equivalent, but are not always true. In particular, if $\Gamma$ is a divisible group in $\rdop$, then (a)--(e) hold.
\end{prop}

\proof 
We first prove that (c) always holds. Let $L$ be an algebraic closure of $K$. We choose a valuation $w$ on $L$ extending $K$. Then the residue field $\tilde{L}$ is also algebraically closed (see \cite{BGR}, Lemma 3.4.1/4). We have to show that $(\Ucal_\Delta)_s \otimes_{\tilde{K}} \tilde{L}$ is reduced for any $\Delta \in \Sigma_1$. We consider the $L^\circ$-subalgebra $B:=K[M]^\Delta \otimes_{K^\circ} L^\circ$ of $L[M]^\Delta$ and we have to prove that $I:=L^{\circ \circ}B=L^{\circ \circ}K[M]^\Delta$ is a radical ideal in $B$.  
Similarly as in the proof of Lemma \ref{reduced special fibre and omega-weight}, $I$ is contained in the radical ideal $J:=\{f \in B \mid w_\Delta(f)>0\}$ of $B$. To prove the reverse inclusion, it is enough to show that  $f$ is contained in $I$ for any $f=\gamma \alpha \chi^u \in J$ with $\gamma \in L^{\circ}$ and $ \alpha \chi^u \in K[M]^\Delta$. If $\gamma \in L^{\circ \circ}$, then the claim is trivial and so we may assume that $\gamma=1$. 
Then $f \in K[M]^\Delta$ and the claim follows from Lemma \ref{reduced special fibre and omega-weight}. This proves (c).

The above shows that  $\{\sum_{u \in M} a_u \chi^u \in K[M] \mid v(a_u)+ \langle u, \omega \rangle >0  \; \forall \omega \in \Delta\}$ is equal to the kernel of  $K[M]^\Delta\rightarrow K[M]^\Delta/\langle K^{\circ \circ}  \rangle$. 
As in the  proof of Proposition \ref{reduced special fibre}, we conclude that (e) holds. Moreover, this proof shows that (a) yields (d). 

It remains to prove that (d) implies (a). We choose $(L,w)$ as above and we note that the value group of $w$ is divisible. Using Proposition \ref{finite generation} and (d), we see that $(\Ycal_\Sigma)_{L^\circ}$ is of finite type over  $L^\circ$. By \cite{EGA IV}, Proposition 2.7.1, the toric scheme $\Ycal_\Sigma$ is of finite type over $K^\circ$ and (a) follows from Proposition \ref{finite presentation in non-discrete case}.

 If $\Gamma$ is a divisible group in $\rdop$, then the vertices of $\Sigma_1$ are always in $N_\Gamma$ proving  the last claim. It is also clear that (a) has not always to be true if $\Gamma$ is not divisible. \qed

\begin{art} \label{lattices for polyhedra} \rm 
For a $\Gamma$-rational polyhedron $\Delta$ in $N_\rdop$, we introduce the following notation: The affine space generated by $\Delta$ is a translate of $(\ldop_\Delta)_\rdop$ for a rational linear subspace $\ldop_\Delta$ of $N_\qdop$. Then  $N_\Delta := N \cap \ldop_\Delta$ and $N(\Delta):=N/N_\Delta$ are free abelian groups of finite rank with quotient homomorphism $\pi_\Delta:N \rightarrow N(\Delta)$. Dually, we have $M(\Delta):=\ldop_\Delta^\bot \cap M=\Hom(N(\Delta),\zdop)$. 
\end{art}

We return to an arbitrary valued field $(K,v)$. Let $\Sigma$ be a $\Gamma$-admissible cone in $N_\rdop \times \rdop_+$ and let $Z$ be an orbit of $\Ycal_\Sigma$ contained in the generic fibre. By \ref{tropicalization and orbits}, $Z$ corresponds to the relative interior of a rational cone $\sigma \in \Sigma_0$

\begin{prop} \label{generic orbit closure}
Under the hypothesis above,  the closure $\overline{Z}$ of $Z$ in $\Ycal_\Sigma$ is isomorphic to the $\Spec( K^\circ[M(\sigma)] )$-toric scheme over $K^\circ$ associated to the $\Gamma$-admissible fan $\Sigma_\sigma:=\{(\pi_\sigma \times \id_{\rdop_+})(\nu) \mid \nu \in \Sigma, \, \nu \supset \sigma \}$ in $N(\sigma)_\rdop \times \rdop_+$.
\end{prop}

\proof Let $\nu \in \Sigma$ with $\nu \supset \sigma$ and let $\nu_\sigma := (\pi_\sigma \times \id_{\rdop_+})(\nu)$. Then there is a canonical surjective $K^\circ$-algebra homomorphism
\begin{equation*}
K[M]^\nu \rightarrow K[M(\sigma)]^{\nu_\sigma}, \quad  \alpha \chi^u \mapsto 
\begin{cases}
\alpha \chi^u, & \text{if $u \in M(\sigma)$,}\\
0, &  \text{if $u \in M \setminus M(\sigma)$}.
\end{cases}
\end{equation*}
We conclude that the $\Spec( K^\circ[M(\sigma)] )$-toric scheme over $K^\circ$ associated to the $\Gamma$-admissible fan $\Sigma_\sigma$ in $N(\sigma)_\rdop \times \rdop_+$ is a closed subscheme of $\Ycal_\Sigma$. By \cite{Fu2}, \S 3.1, its generic fibre is the closure of $Z$ in the generic fibre of $\Ycal_\Sigma$. By Proposition \ref{characterization of affine closure}, we get the claim. \qed

\vspace{2mm}
Now we assume that the orbit $Z$ of $\Ycal_\Sigma$ is contained in the special fibre. By \ref{tropicalization and orbits}, $Z$ corresponds to the relative interior $\tau$ of $\Delta \in \Sigma_1$. Similarly as in Proposition \ref{irreducible component is toric}, $\Gamma$-rationality of $\Delta$ yields that $M(\Delta)_\tau := \{u \in M(\Delta) \mid \langle u, \omega \rangle \in \Gamma \; \forall \omega \in \tau \}$ is a lattice of finite index in $M(\Delta)$. For $\nu \in \Sigma_1$ with face $\Delta$,  we define ${\rm LC}_\tau(\nu) := {\rm LC}_\omega(\nu)$ which is independent of the choice of $\omega \in \tau$ and where we use the local cones from  \ref{local cone}.

\begin{prop} \label{special orbit closure}
Under the hypothesis above,  the closure $\overline{Z}$ of $Z$ in $\Ycal_\Sigma$ is equivariantly isomorphic to the $\Spec( \tilde{K}[M(\Delta)_\tau] )$-toric variety  over $\tilde{K}$ associated to the rational fan $\{\pi_\Delta ({\rm LC}_\tau(\nu)) \mid \nu \in \Sigma_1, \, \nu \supset \Delta \}$ in $N(\Delta)_\rdop$.
\end{prop}

\proof If $\omega$ is a vertex of $\Sigma_1$, this follows immediately from Corollary \ref{irreducible component is toric}. The general case follows from the corresponding generalization of  Corollary \ref{irreducible component is toric} which can be proved completely analogous. We leave the details to the reader. \qed

\section{Tropical cone of a variety}   \label{TropCone}

In this section, $K$ denotes a field with a non-trivial non-archimedean absolute value $|\phantom{a}|_v$, corresponding valuation $v:=-\log|\phantom{a}|_v$ and valuation ring $K^\circ$. We consider $W:=\{\ve v \mid \ve \geq 0\}$ which is induced by all valuations equivalent to $v$ together with the trivial absolute which we denote by $0$. We may identify $W$ with $\rdop_+$ using $\ve v \leftrightarrow \ve$. The value group of $w \in W$ is denoted by $\Gamma_w$ and the residue field by $k(w)$. Obviously, we have $k(w)=k(v)$ for $w \neq 0$ and $k(0)=K$. At the end of this section, we show how to adjust the notation so that everything works also for the trivial valuation.

We have seen  the advantage of using fans $\Sigma$ in $N_\rdop \times \rdop_+$ rather than polyhedral complexes in $N_\rdop$ to define an associated toric scheme $\Ycal_\Sigma$ over $K^\circ$. It is not surprising that the consideration of the closed cone in $N_\rdop \times \rdop_+$ generated by $\Trop_v(X) \times \{1\}$ is useful to describe information about the closure of the closed subscheme $X$ of $T$ in $\Ycal_\Sigma$ in a uniform way. Moreover, we will see that the tropical variety of $X$ with respect to the trivial valuation is just the intersection of this tropical cone with $N_\rdop \times \{0\}$. 

\begin{art} \rm \label{W-analytification}
For an algebraic scheme $X$ over $K$, we have defined in Section \ref{Section:Analytification} the analytification with respect to the valuation $w$ which we denote here by $X_w^{\rm an}$. In fact, we can define the {\it analytification $X_W^{\rm an}$ of $X$  with respect to $W$} by the same process allowing all multiplicative seminorms $p$ with restriction $w_p:=p|_K \in W$. This gives again a locally compact Hausdorff space which is as a set equal to the disjoint union of all $X_w^{\rm an}$ with $w$ ranging over $W$. The proof follows from Tychonoff's theorem similarly as in the case of a single valuation.
\end{art}

\begin{art} \rm \label{W-tropicalization}
For $w \in W$, let $\trop_w: T_w^{\rm an} \rightarrow N_\rdop$ be the tropicalization map. Proceeding fibrewise, we get the {\it $W$-tropicalization map}
$$\trop_W: T^{\rm an}_W \rightarrow N_\rdop \times W, \quad t \mapsto (\trop_{w_t}(t), w_t).$$
It is clear that $\trop_W$ is continuous. 
\end{art}

\begin{Def} \rm \label{tropical cone of X}
Let $X$ be a closed subscheme of $T$. Then we define the {\it tropical cone} associated to $X$ as $\trop_W(X_W^{\rm an})$ and we denote it by $\Trop_W(X)$. 
\end{Def}

\begin{figure}[h]
\begin{center}$
\begin{array}{cc}
\setlength{\unitlength}{0.14cm}
\ifpdf
\begin{picture}(25,40)
\put(-38,-12){
\includegraphics[width=10.6cm]{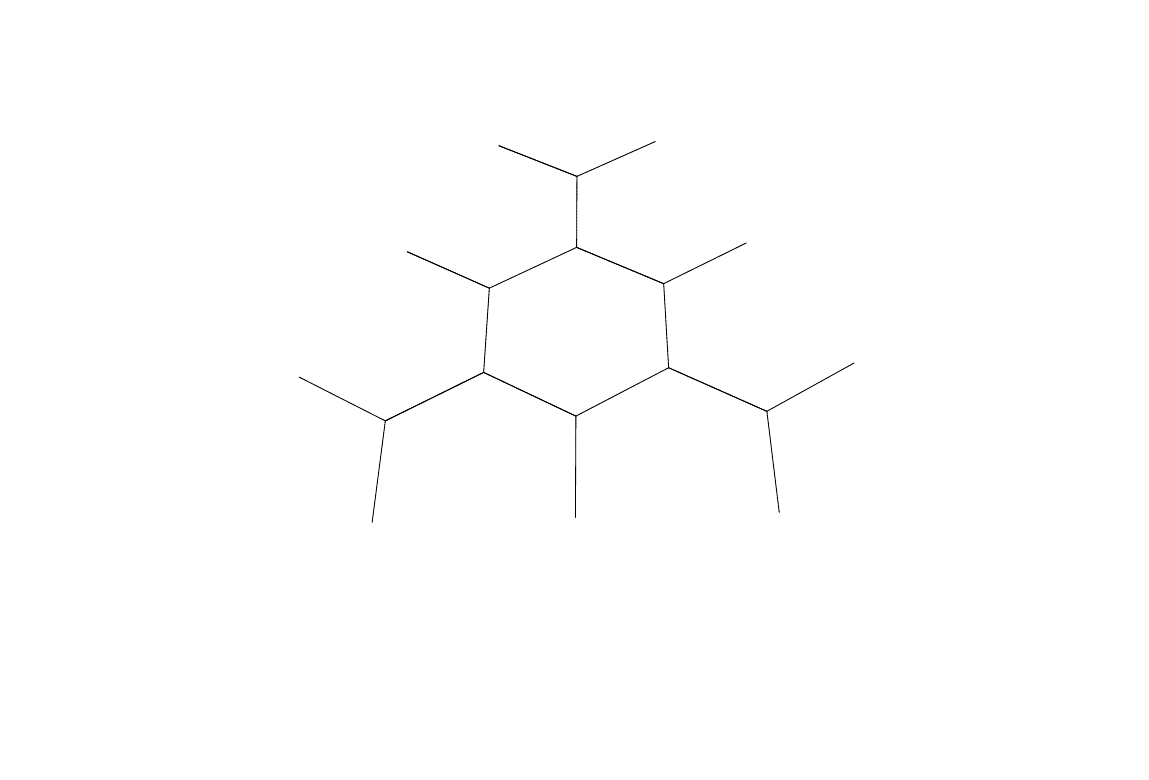}}
\end{picture} 
\else
\begin{picture}(25,40)
\put(-38,-12){
 \includegraphics[width=10.6cm]{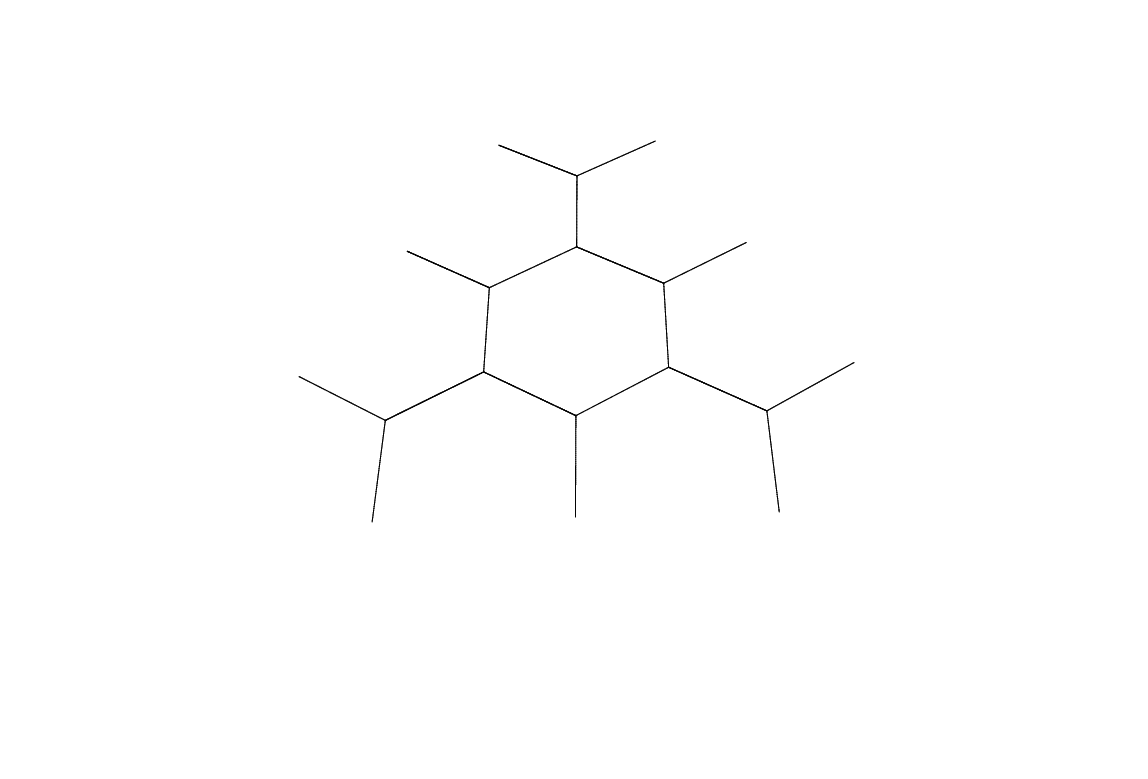}}
\end{picture}  
\fi
\begin{picture}(25,40)
\ifpdf
\put(-20,-8.0){
\includegraphics[totalheight=6.2cm]{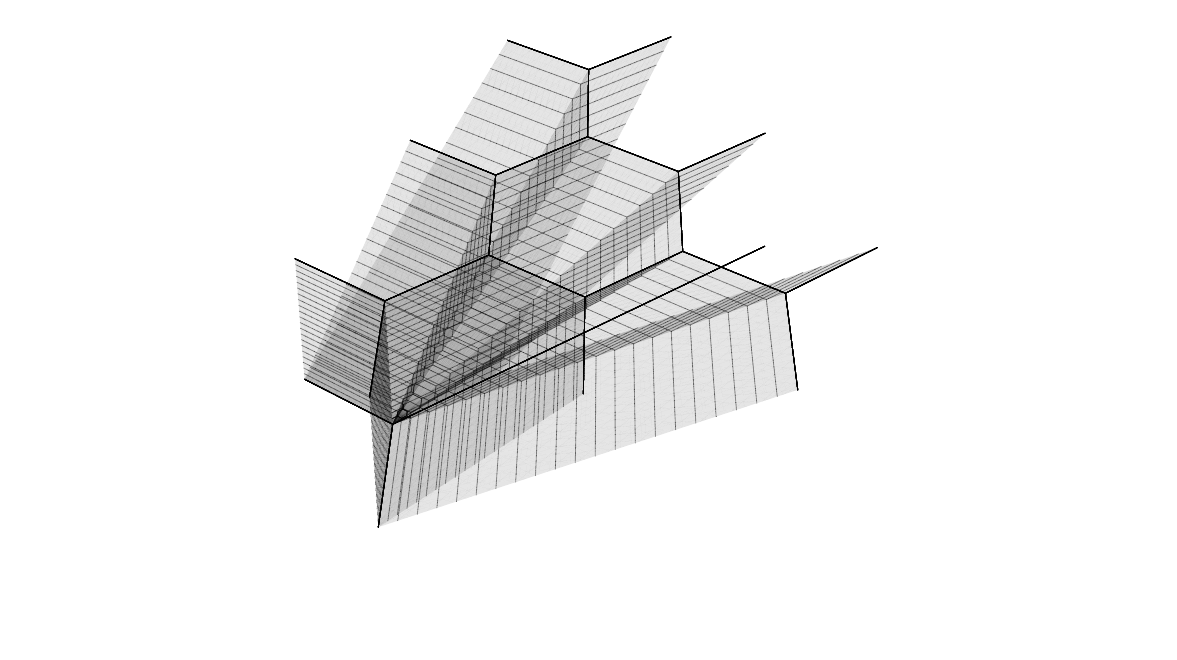}}
\else
\put(-20,-8){
\includegraphics[totalheight=6.2cm]{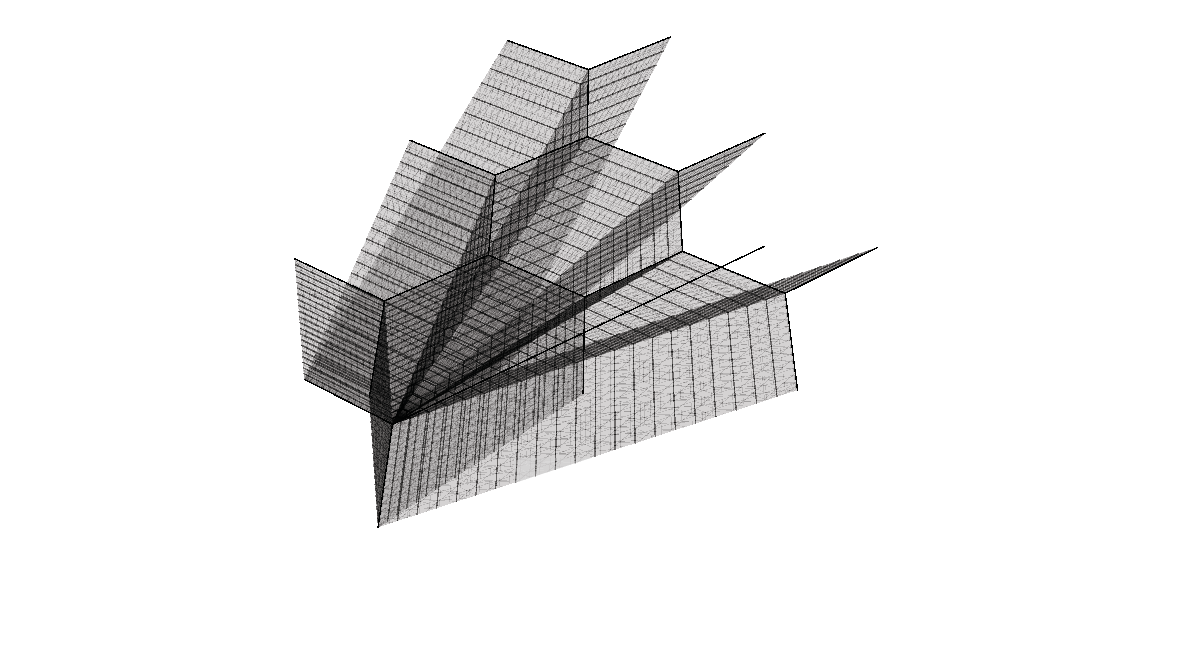}}
\fi
\end{picture}
\end{array}$
\end{center}
\caption{Tropical curve and its tropical cone}
\end{figure}

\begin{rem} \label{cone property} \rm
We will see in Corollary \ref{cone property proved} that $\Trop_W(X)$ is equal to the closure $C$ of the cone in $N_\rdop \times \rdop_+$ generated by $\Trop_v(X)\times \{1\}$. It follows from the Bieri--Groves theorem that $\Trop_W(X)$ is a finite union of $\Gamma$-admissible cones in $N_\rdop \times \rdop_+$. If $X$ is of pure dimension $d$, then we may choose these cones $d+1$-dimensional. 

We will not use these facts until we have proved them in Corollary \ref{cone property proved}. At the moment, it is only clear from the definitions that $C$ agrees with $\Trop_W(X)$ on $N_\rdop \times \rdop_+$. 
\end{rem}

\begin{prop} \label{tropical cone properties}
The tropical cone $\Trop_W(X)$ of $X$ is closed in $N_\rdop \times \rdop_+$.
\end{prop}

\proof We have seen that $\trop_W$ is a continuous map. Moreover, it is easy to see that $\trop_W$ is a proper map, i.e. the preimage of a compact subset is a compact subset. This shows immediately that $\Trop_W(X)=\trop_W(X_W^{\rm an})$ is closed in $N_\rdop \times \rdop_+$. \qed




\vspace{2mm}

In the remaining part of this section, let $\Sigma$ be a $\Gamma$-admissible fan in $N_\rdop \times \rdop_+$ and let $\Ycal_\Sigma^v$ be the associated toric scheme over $K^\circ$. We may identify $T$ with the dense open orbit in $\Ycal_\Sigma^v$ and this  orbit is contained in the generic fibre $Y_{\Sigma_0}$. 

In the following proposition, we consider a field extension $L/K$ and a valuation $u$ on $L$ (not necessarily of rank $1$) extending $v$ with valuation ring $L^\circ$ and value group $\Gamma_L$. Since $\Gamma_L$ is a totally ordered abelian group, any $\Gamma$-rational polyhedron $\Delta$ in $N_\rdop$ induces a polyhedron $\Delta(\Gamma_L)$ in $N_{\Gamma_L}$. We denote by $|\Sigma_1(\Gamma_L)|$ the union of all $\Delta(\Gamma_L)$ with $\Delta$ ranging over the $\Gamma$-rational polyhedral complex $\Sigma_1$ defined in \ref{admissible cones}.  

\begin{prop} \label{integral points of toric variety}
Under the assumptions above, $P \in T(L)$ is an $L^\circ$-integral point of $\Ycal_\Sigma
^v$ if and only if $\trop_u(P)$ is contained in  $|\Sigma_1(\Gamma_L)|$. 
\end{prop}

\proof 
This follows immediately from Proposition \ref{integral points and tropicalization revisted}. \qed

\begin{art} \rm \label{reduction map over W}
We conclude from Proposition \ref{integral points of toric variety} that we have a well-defined reduction map $\pi_W:\trop^{-1}_W(|\Sigma|) \rightarrow \Ycal_\Sigma^v$. Indeed,  we have $\Ycal_{\ve \Sigma}^{\ve v} = \Ycal_\Sigma^v$ for all $\ve >0$ and so we may use the reduction map $\pi_w:\trop_w^{-1}(|\Sigma_\ve|) \rightarrow  (\Ycal_{\ve\Sigma}^w)_s = (\Ycal_\Sigma^v)_s$ in the fibre over $w=\ve v$.  For  $w=0$, the special fibre  agrees with the generic fibre $Y_{\Sigma_0}$ and we use  the reduction  $\pi_0:\trop_0^{-1}(|\Sigma_0|) \rightarrow Y_{\Sigma_0}$. 
Note that  we may use Proposition \ref{integral points of toric variety} also for the trivial valuation $v=0$. 
\end{art}

Then we can describe the orbit-face correspondence in the following uniform way.

\begin{prop} \label{uniform orbit correspondence}
There is a bijective order reversing correspondence between $\tdop$-orbits $Z$ of $\Ycal_\Sigma^v$ and open faces $\tau$ of $\Sigma$ given by
$$Z=\pi_W(\trop_W^{-1}(\tau)), \quad \tau = \trop_W(\pi_W^{-1}(Z)).$$
\end{prop}

\proof We easily reduce to the case of an invariant open subset $\Vcal_\sigma$ of $\Ycal_\Sigma^v$ for $\sigma \in \Sigma$. Then the claim follows from Proposition \ref{affine orbit correspondence} applied to every $w \in W$. To prove that the correspondence is order reversing we use also Proposition \ref{generic orbit closure} to handle an orbit in the generic fibre whose closure contains orbits in the special fibre. \qed

\begin{rem} \label{tropical cone and trivial valuation} \rm
If $v$ is the trivial valuation, then we have to adjust the notation of this section by using the set $\rdop_+$ rather than $W=\{0\}$. We define $X_W^{\rm an}:=X_0^{\rm an} \times \rdop_+$ which is  a locally compact Hausdorff space. Then everything works as above.
\end{rem}

\section{Projectively embedded toric varieties}

In this section, $K$ denotes a  field with a non-archimedean absolute value $|\phantom{a}|$, corresponding valuation $v:=-\log|\phantom{a}|$ and value group $\Gamma := v(K^\times)$. We have defined toric varieties in Definition \ref{toric variety}. Here, we consider projective toric varieties over $K^\circ$ with an equivariant embedding into projective space. These toric varieties are not necessarily normal. This section is inspired by the introductory article of E. Katz (\cite{Ka}, section 4) and we will generalize his results. Further references:  \cite{CLS}, \S 2.1, \S 3.A; \cite{GKZ}, Chapter 5.

Recall that $\tdop = \Spec(\Tor)$ is a split multiplicative torus over $K^\circ$ with generic fibre $T$. The character group of $T$ is $M$ and the character corresponding to $u \in M$ is denoted by $\chi^u$. For convenience, we always choose  coordinates on the projective space $\pdop^N_{K^\circ}$ defined over the valuation ring $K^\circ$.

\begin{art} \rm \label{A-example in the field case}
We first recall the following well-known way to construct a not necessarily normal toric subvariety $Y$ from $A = (u_0, \dots, u_N) \in M^{N+1}$ and $\yb=(y_0: \dots: y_N) \in \pdop^N(K)$ (see \cite{GKZ}, Chapter 5). The torus $T$ acts on $\pdop^N_K$ by 
$$t \cdot \xb:= (\chi^{u_0}(t)x_0: \dots :\chi^{u_N}(t)x_N)$$
and we define $Y$ as the closure of the orbit $T\yb$. Then there is a bijective correspondence between $T$-orbits of $Y$ and faces of the weight polytope $\Wt(\yb)$ which is defined as the convex hull of $A(\yb):=\{u_j \mid y_j \neq 0\}.$ If $Q$ is a face of $\Wt(\yb)$, then the corresponding orbit is given by 
$$Z:=\{\zb \in Y \mid z_j \neq 0 \iff u_j\in Q\}.$$
Duality gives also a bijective correspondence to the normal fan $\Sigma$ of $\Wt(\yb)$. The cone $\sigma$ corresponding to the face $Q$ is the set of $\omega \in N_\rdop$ such that the linear functional $\langle \cdot, \omega \rangle$ achieves its minimum on $\Wt(\yb)$ precisely in the face $Q$. The torus corresponding to the orbit $Z$ has character group $\{\sum m_j u_j \mid \sum m_j =0 \}$, where $j$ ranges over $0, \dots, N$ with $z_j \neq 0$ for any $\zb \in Z$. This character group is 
of finite index in $\zdop (\sigma^\bot \cap  M)$ and hence $\dim(Z)=\dim(Q)=n-\dim(\sigma)$ (see \cite{CLS}, Section 3.A). 
\end{art}

\begin{art} \rm \label{A-example}
The goal of this section is to perform a similar construction over the valuation ring $K^\circ$. Let $A= (u_0, \dots, u_N) \in M^{N+1}$ and let $\yb=(y_0: \dots: y_N) \in \pdop^N(K)$. We define the {\it height function} of $\yb$ by  
$$a:\{0, \dots,N\}  \rightarrow \Gamma \cup \{\infty\}, \quad j \mapsto v(y_j).$$
The torus $\tdop$ operates on $\pdop^N_{K^\circ}$ by 
$$\tdop \times_{K^\circ} \pdop^N_{K^\circ} \rightarrow \pdop^N_{K^\circ},
\quad (t, \xb) \mapsto (\chi^{u_0}(t)x_0: \dots :\chi^{u_N}(t)x_N).$$
The closure of the orbit $T\yb$ in $\pdop^N_{K^\circ}$ is a projective toric variety with respect to the split torus over $K^\circ$ with generic fibre $T/\Stab(\yb)$. We denote this projective toric variety by $\Ycal_{A,a}$ and its generic fibre by $Y_{A,a}$. 
Using the base point $\yb$, the torus $T/\Stab(\yb)$ may be seen as an open dense subset of $\Ycal_{A,a}$. If we apply the following result to the special case $M$ equal to the lattice in $M'$ generated by $A$ and $F$ equal to a translation, then we see that $\Ycal_{A,a}$ depends only on the affine geometry of $(A,a)$.
\end{art}

\begin{prop} \label{affine A-relation}
Suppose that $\tdop'$ is another split multiplicative torus over $K^\circ$ with character lattice $M'$ and that there is an injective affine transformation $F: M \rightarrow M'$ of lattices. 
Let $A=(u_0, \dots, u_N) \in M^{N+1}, A'=(u'_0, \dots, u'_N) \in (M')^{N+1}$ and let $\yb, \yb' \in \pdop^N(K)$ with height functions $a$ (resp. $a'$). Let $\Ycal_{A,a}$ (resp. $\Ycal_{A',a'}$) be the projective toric variety with respect to $A$, $\yb$ (resp. $A'$, $\yb'$). We assume that $F(u_j)=u_j'$ for every $j$ with $y_j\neq 0$. If there is $\lambda \in \Gamma$ such that $a'= a+ \lambda$, then $\Ycal_{A,a}$ is canonically isomorphic to $\Ycal_{A',a'}$. 
\end{prop}

\proof The injective linear map corresponding to $F$  induces a surjective homomorphism $\tdop' \rightarrow \tdop$ of multiplicative tori. If $\yb = \yb'$, then we deduce that $\Ycal_{A,a}=\Ycal_{A',a}$. In general, we have $\yb' = g \yb$ for some $g=(g_0, \dots, g_N) \in K^{N+1}$ with $|g_0|= \dots = |g_N|\neq 0$. Then $g$ induces a linear automorphism of $\pdop^N_{K^\circ}$ mapping $\Ycal_{A,a}$ onto $\Ycal_{A',a'}$. If $y_j \neq 0$, then $g_j$ is uniquely determined and hence we have constructed a canonical isomorphism. \qed

\begin{cor} \label{dimension}
The open dense orbit $T/\Stab(\yb)$ of $\Ycal_{A,a}$ is a torus with character lattice $\{\sum m_j u_j \mid \sum m_j =0 \}$, where $j$ ranges over $0, \dots, N$ with $y_j \neq 0$. This orbit has dimension equal to $\dim(\adop)$, where $\adop$  is the smallest affine subspace of $M_\rdop$ containing $A(\yb):=\{u_j \mid y_j \neq 0\}$.
\end{cor}

\proof Using Proposition \ref{affine A-relation}, we may assume that $u_0=0$ and that $u_1, \dots, u_N$ form a basis of $M$. Then it is easy to see that the stabilizer of $\yb$ is trivial. The corollary is also a special case of the result mentioned at the end of \ref{A-example in the field case}. \qed

\vspace{2mm}

The following result is well-known for fields or discrete valuation rings. We need it for arbitrary valuation rings of rank $1$ which are not noetherian in general and hence we may not use algebraic intersection theory. However, there is an intersection theory with Cartier divisors in this situation (see \cite{Gu1}) which together with the result for the generic fibre will easily imply the claim.

\begin{prop} \label{Picard groups}
The restriction map gives an isomorphism $\Pic(\pdop^N_{K^\circ}) \rightarrow \Pic(\pdop^N_{K})$ and pull-back with respect to the second projection gives an isomorphism  $\Pic(\pdop^N_{K^\circ}) \rightarrow \Pic(\tdop \times_{K^\circ} \pdop^N_{K^\circ})$
\end{prop}

\proof To prove the first claim, we have to show that every line bundle $\Lcal$ on $\pdop^N_{K^\circ}$ is isomorphic to $O_{\pdop^N}(m)$ over $K^\circ$ for some $m \in \zdop$. We consider a Cartier divisor $\Dcal = \{U_i, \gamma_i\}_{i \in I}$ which is trivial on the generic fibre $\pdop^N_{K}$ and we have to prove that $\Dcal$ is trivial. We may assume that $\Ucal_i = \Spec(K^\circ[x_1, \dots, x_N]_{h_i})$ for a  polynomial $h_i \in  K^\circ[x_1, \dots, x_N]$. 	Obviously, we may skip all charts with empty special fibre. This means that the prime factors $p_1, \dots, p_r$ of $h_i$ in $K^\circ[x_1, \dots, x_N]$ are non-constant.    Using unique factorization, we get $\Ocal(U_i)^\times = K^\times p_1^\zdop \cdots p_r^\zdop$ for the generic fibre $U_i$ of $\Ucal_i$. By triviality of $\Dcal$ on the generic fibre $\pdop^N_K$, we get $\gamma_i=\lambda_i h_i'$ for some $\lambda_i \in K^\times$ and $h_i' \in p_1^\zdop \cdots p_r^\zdop$. We want to show that these factorizations fit on an  overlapping $\Ucal_i \cap \Ucal_j$. As $\Ucal_i \cap \Ucal_j$ intersects the special fibre, there is a valued field $(L,w)$ extending $(K,v)$ and an $L^\circ$-integral point $P$ of $\Ucal_i \cap \Ucal_j$. Using $h_i \in \Ocal(\Ucal_i)^\times$, we get $|h_i'(P)|_w=1$. The multiplicity $m(\Dcal, \pdop^N_{\tilde{K}})$ of $\Dcal$ along the special fibre $ \pdop^N_{\tilde{K}}$ was defined in \cite{Gu1}, Section 3. Since the special fibre is irreducible and smooth, it is shown in Proposition 7.6 of \cite{Gu2} that 
$$ m(\Dcal, \pdop^n_{\tilde{K}})=-\log |\gamma_i(P)|=v(\lambda_i).$$
We conclude that $v(\lambda_i)=v(\lambda_j)$. Dividing the equations of $\Dcal$ by a fixed $\lambda_i$, we deduce that $\Dcal$ is trivial on $\pdop^N_{K^\circ}$ proving the first claim.

Similarly, we prove the second claim. The claim holds on the generic fibre and hence it is enough to show that a Cartier divisor $\Dcal$ on $\tdop \times_{K^\circ} \pdop^N_{K^\circ}$ which is trivial on the generic fibre $T \times_K \pdop^N_{K}$ is trivial on $\tdop \times_{K^\circ} \pdop^N_{K^\circ}$. This is done as above replacing $K$ by the unique factorization domain $\Ocal(T)=K[M]$ and using $K[M]^\times=\{\lambda \chi^u \mid u \in M, \, \lambda \in K^\times\}$. \qed

\begin{rem} \label{Picard and product in general} \rm
It was pointed out to the author by Qing Liu and C. P\'epin that the second claim holds more generally for any integral normal scheme $\Xcal$ instead of $\pdop^N_{K^\circ}$. Their argument is as follows:  Injectivity follows from the existence of a section for $p_2$. It remains to prove surjectivity. Let $\xi$ be the generic point of $\Xcal$ with residue field $\kappa(\xi)$. Then the fibre $(\tdop \times_{{K^\circ}} \Xcal)_\xi$ is isomorphic to the split torus $T_{\kappa(\xi)}$ and hence the restriction of any line bundle on $\tdop \times_{{K^\circ}} \Xcal$ to this fibre is trivial. We conclude that it is enough to show that pull-back with respect to the second projection $p_2$  gives an isomorphism from the group of Cartier divisors of $\Xcal$ onto the group of those Cartier divisors  of $\tdop \times_{{K^\circ}} \Xcal$ whose restriction to $(\tdop \times_{{K^\circ}} \Xcal)_\xi$ is zero. 
Injectivity follows again from the existence of a section for $p_2$. Moreover, this shows that surjectivity is a local question on $\Xcal$ and so we may assume $\Xcal$ affine. By the descent argument in \cite{EGA IV}, Proposition 8.9.1, we may assume that $\Xcal$ is noetherian. Then surjectivity follows from Proposition 21.4.11 in the list of Errata and Addendum in \cite{EGA IV}.

The referee gave the following alternative argument for the first claim in Proposition \ref{Picard groups}. It is obvious that the restriction $\Pic(\pdop^N_{K^\circ}) \rightarrow \Pic(\pdop^N_{K})$ is surjective. He noted that injectivity of $\Pic(\Xcal) \rightarrow \Pic(\Xcal_K)$ 
holds more generally for any projective and flat scheme $\Xcal$ over $K^\circ$ with irreducible and reduced geometric fibres if the structure morphism $f:\Xcal \rightarrow S:=\Spec(K^\circ)$ admits a section and if  $f_*(\Ocal_\Xcal)=\Ocal_{S}$. Note that these assumptions hold for $\pdop^N_{K^\circ}$. His argument uses that $\Pic_{\Xcal/S}(S)=\Pic(\Xcal)$ and $\Pic_{\Xcal/S}(\Spec(K))=\Pic(\Xcal_K)$ for the relative Picard functor $\Pic_{\Xcal/S}$ (see \cite{BLR}, Proposition 8.1.4). By a result of Grothendieck, $\Pic_{\Xcal/S}$ is represented by a separated $S$-scheme. Note that Grothendieck's original result (\cite{Gro}, no. 232, Theorem 3.1) was written for locally noetherian schemes, but using the technique of noetherian approximation introduced in \cite{EGA IV}, this holds in general as stated in \cite{BLR}, Theorem 8.2.1. Now the valuative criterion of separatedness (\cite{EGA II}, Proposition 7.2.3) shows that the map 
$\Pic(\Xcal) \rightarrow \Pic(\Xcal_K)$  is injective.

\end{rem}                                                                       

\begin{lem} \label{lift of projective linear action}
Suppose that the torus $\tdop$ acts linearly on $\pdop^N_{K^\circ}$. Then this action lifts to a linear representation of $\tdop$ on $\adop^{N+1}_{K^\circ}$.
\end{lem}

\proof The  action $\sigma: \tdop \times_{K^\circ} \pdop^N_{K^\circ} \rightarrow \pdop^N_{K^\circ}$ is linear which means that it  is given by a homomorphism $\tdop \rightarrow PGL(N+1)$ defined over $K^\circ$. We are looking for a lift to a homomorphism $\tdop \rightarrow GL(N+1)$. This is equivalent to the existence of a $\tdop$-linearization on the line bundle $L:=O(1)$ of $\pdop^N_{K^\circ}$, i.e. an action of $\tdop$ on $L$ which is compatible with the given group action $\sigma$. Here, we use the language of \cite{Mu}, \S 1.3, which is written for schemes over a base field. However, the argument for the existence of a $\tdop$-linearization in \cite{Mu}, Proposition 1.5, extends to the case of a valuation ring. Indeed,  the essential point is the existence of an isomorphism $\sigma^*(L) \cong p_2^*(L)$ which follows from Proposition \ref{Picard groups} and then we may conclude as at the end of the proof of  \cite{Mu}, Proposition 1.5, to prove that $L=O(1)$ has a $\tdop$-linearization. \qed

\begin{prop} \label{converse}
Let $\Ycal$ be a closed irreducible subvariety of $\pdop^N_{K^\circ}$. Suppose that the torus $\tdop$ operates linearly on $\pdop^N_{K^\circ}$ and leaves $\Ycal$ invariant. We assume that $\Ycal$ has an open dense orbit containing a $K$-rational point $\yb$. Then after a suitable change of coordinates on $\pdop^N_{K^\circ}$,  there is $A \in \check{T}^{N+1}$  such that $\Ycal = \Ycal_{A,a}$ for the height function $a:\{0, \dots, N\} \rightarrow \Gamma$ of $\yb$.
\end{prop}

\proof By Lemma \ref{lift of projective linear action}, the projective representation of $\tdop$ on $\pdop^N_{K^\circ}$ lifts to a representation $S$ of $\tdop$ on $\adop^{N+1}_{K^\circ}$.  
Since the multiplicative torus $T$ is split over $K$, it follows  that the vector space $V:=K^{N+1}$ has a simultaneous eigenbasis $v_0, \dots ,v_N$ for the $T$-action  (\cite{Bor}, Proposition  III.8.2). For $j=0, \dots, N$, we have $S_t(v_j)=\chi^{u_j}(t)v_j$ for all $t \in T(K)$ and some $u_j \in M$. 

We endow $V$ with the norm $\|\xb\|:= \max_{j=0,\dots,N} |x_j|$ for $\xb \in K^{N+1}$. By definition, this gives a $K$-cartesian space and hence every subspace $U$ of $V$ is also $K$-cartesian meaning that there is a basis $u_1,\dots,u_r$ of $U$ such that $\|\sum_{j=1}^r \alpha_j u_j\|=\max_{j=1,\dots,r} |\alpha_j|$ for all $\alpha_1, \dots, \alpha_r \in K$ (see \cite{BGR}, Proposition 2.4.1/5). In this non-archimedean situation, such a basis is called orthonormal. We apply this to the simultaneous eigenspaces $V_{u_j}$ for the $T$-action and so we may choose the simultaneous eigenbasis $v_0, \dots ,v_N$ above in such a way that a suitable subset is an orthonormal basis  of $V_{u_j}$ for every $j=0, \dots, N$.

We consider the subgroup $U:=\tdop(K^\circ)=\{t \in T(K) \mid v(t_1)= \dots = v(t_n)=0\}$ of $T(K)$. For $t \in U$, we have $S_t \in GL(N+1, K^\circ)$ and hence the eigenvalues $\chi^{u_j}(t)$ have absolute value $1$. If we use reduction modulo $K^{\circ \circ}$, then the $U$-action becomes a $(\Tor)_{\tilde{K}}$-operation on ${\tilde{K}}^{N+1}$. We note that the reduction of an orthonormal basis in a subspace of $V$ is linearly independent in ${\tilde{K}}^{N+1}$. Using that eigenvectors for distinguished eigenvalues are linearly independent, we conclude that the reduction of $v_0, \dots, v_N$ is a a simultaneous eigenbasis for the  $(\Tor)_{\tilde{K}}$-action. By Nakayama's Lemma, it follows that $v_0, \dots, v_N$  is a $K^\circ$-basis for $(K^\circ)^{N+1}$. 
We choose the coordinates of $\pdop^N_{K^\circ}$ according to this basis and let $a$ be the corresponding height function of $\yb$. For $A=(u_0, \dots, u_N)$, we get $\Ycal = \Ycal_{A,a}$. \qed

\begin{rem} \rm \label{other projective toric varieties}
Every projective normal toric variety over a field can be equivariantly embedded into some projective space endowed with a linear torus action (see \cite{Mu}, \S 1.3).   There are projective non-normal toric varieties over a field for which this is not true (see \cite{GKZ}, Remark 5.1.6).
\end{rem}

\begin{art} \label{weight complex} \rm
In the following, we consider the toric variety $\Ycal_{A,a}$ for some $A \in M^{N+1}$ and $\yb \in \pdop^{N+1}_{K^\circ}$ with height function $a$. 

The {\it weight polytope} $\Wt(\yb)$ is the convex hull of $A(\yb):=\{u_j \mid a(j) < \infty\}$ in $M_\rdop$. The {\it induced subdivision} of $\Wt(\yb)$ is given by projection of the faces of the convex  hull of $\{(u_j,\lambda) \in M_\rdop \times \rdop \mid  j=0, \dots, N, \, \lambda \geq a(j)\}$. The weight subdivision is a polytopal complex denoted by $\Wt(\yb,a)$. The vertices of $\Wt(\yb,a)$ are contained in $A(\yb)$. 
\end{art}

\begin{art} \rm \label{dual complex for Y}
In the following, we need some additional notions from convex geometry which we have introduced in the appendix. By construction, there is a unique proper  polyhedral function $f$ on $M_\rdop$ such that the epigraph of $f$ is equal to the convex hull of $\{(u_j,\lambda) \in M_\rdop \times \rdop \mid  j=0, \dots, N, \, \lambda \geq a(j)\}$. The domain of $f$ is equal to  $\Wt(\yb)$ and $f(u_j)=a(j)$ for all vertices $u_j$ of the weight subdivision $\Wt(\yb,a)$.

We define the {\it dual complex} $\Ccal(A,a)$  of $\Wt(\yb,a)$ as the complete polyhedral complex in $N_\rdop$  characterized by the fact that the $n$-dimensional polyhedra in $\Ccal$ are the domains of linearity of the affine function
$$g(\omega):= \min_{j=0,\dots, N} a(j)+ \langle u_j , \omega \rangle.$$
Obviously, all polyhedra in $\Ccal(A,a)$ are $\Gamma$-rational. There is a bijective order reversing correspondence between the faces  of $\Wt(\yb,a)$ and polyhedra in $\Ccal(A,a)$. The polyhedron $\widehat{Q} \in \Ccal(A,a)$ corresponding to the face $Q$ of $\Wt(\yb,a)$ is given by
\begin{equation*}
\begin{split}
\widehat{Q} &= \{\omega \in N_\rdop \mid g(\omega)= \langle u , \omega \rangle +f (u) \; \forall u \in Q \} \\
&=\{\omega \in N_\rdop \mid g(\omega)= \langle u_j , \omega \rangle +a(j) \; \forall u_j \in A(\yb) \cap Q \}. 
\end{split}
\end{equation*}
Conversely, the face $\widehat{\sigma}$ of $\Wt(\yb,a)$ corresponding to $\sigma \in \Ccal(A,a)$ is given by
$$\widehat{\sigma} = \{ u \in M_\rdop \mid  g(\omega)= \langle u, \omega \rangle + f(u) \; \forall \omega \in \sigma\}$$
and it is also the convex hull of $\{ u_j \in A \mid  g(\omega)= \langle u_j, \omega \rangle + a(j) \; \forall \omega \in \sigma\}$. 
All this can be seen using the dual complex $\Wt(\yb,a)^f$ from \ref{dual complex} and the conjugate polyhedral function $f^*$ of $f$ from \ref{conjugate}. Indeed, we have $f^*(\omega)=-g(-\omega)$ and hence $\Wt(\yb,a)^f=-\Ccal(A,a)$. 

\end{art}
In the next results, we  will also use the tropicalization map ${\trop_v}: \Tan \rightarrow N_\rdop$ and the reduction map $\pi:Y_{A,a}^{\rm an} \rightarrow (\Ycal_{A,a})_s$.

\begin{prop} \label{duality}
There are bijective order reversing correspondences between

\begin{itemize}
 \item[(a)] faces $Q$ of the weight subdivision $\Wt(\yb,a)$;
 \item[(b)] polyhedra $\sigma$ of the dual complex $\Ccal(A,a)$;
 \item[(c)] $\tdop$-orbits $Z$ of the special fibre of $\Ycal_{A,a}$.
\end{itemize}

The  correspondences are given as follows: The face $Q= \widehat{\sigma}$ is the face of $\Wt(\yb,a)$ spanned by those $u_j$ with $x_j \neq 0$ for $x \in Z$. The polyhederon $\sigma$ is given by $\sigma = \widehat{Q}$ and $\relint(\sigma) = {\trop_v}(\{t \in \Tan \mid \pi(t \yb) \in Z\})$. The orbit $Z$ is equal to
$$ \{\xb \in (\Ycal_{A,a})_s \mid x_j \neq 0 \iff u_j \in A(\yb) \cap Q\}=\{\pi(t\yb)\mid t \in \Tan \cap {\trop_v}^{-1}( \relint(\sigma))\}.$$
The correspondence $Q \leftrightarrow Z$ is preserving the natural orders and the other correspondences are order reversing. Moreover,  we have $\dim(Q)= \dim(Z)=n - \dim(\sigma)$.
\end{prop}

\proof We have discussed the  correspondence  $Q \leftrightarrow \sigma$ in \ref{dual complex for Y}. Next, we note that every point $z$ of $(\Ycal_{A,a})_s$ is the reduction of a point in $\Tan\yb$. Since $T\yb$ is an open dense subset of the generic fibre of $\Ycal_{A,a}$, this follows from  Lemma \ref{lifting lemma}.

Now let $\sigma$ be a polyhedron from $\Ccal(A,a)$. We will show next that $Z:=\{\pi(t\yb)\mid t \in \Tan \cap {\trop_v}^{-1}( \relint(\sigma))\}$ is a $\tdop$-invariant subset of  $(\Ycal_{A,a})_s$. Let us consider the formal affinoid torus $T^\circ$ which is the affinoid subdomain of $T^{\rm an}$ given by the equations $|x_1|= \dots = |x_n|=1$. The reduction map  induces a surjective group homomorphism $T^\circ \rightarrow \tdop_s$ and  $\pi:Y_{A,a}^{\rm an} \rightarrow (\Ycal_{A,a})_s$ is equivariant with respect to this homomorphism. Since $T^\circ$ leaves  ${\trop_v}^{-1}(\relint(\sigma))$ invariant, we conclude that $Z$ is invariant under the $\tdop_s$-action. 

For $z \in (\Ycal_{A,a})_s$, we have seen above that there is  $t \in \Tan$ with $z=\pi(t \yb)$. It follows from \ref{dual complex for Y} that $\omega \in \relint(\sigma)$ if and only if 
$$A(\yb) \cap Q = \{u_j \in A \mid g(\omega)= a(j) + \langle u_j , \omega \rangle\},$$ 
i.e. precisely the functions $a(j) +\langle u_j, \omega \rangle$ with $u_j \in A(\yb) \cap Q $ are minimal in $\omega$.  If we apply this with $\omega:= {\trop_v}(t)$, then we have $a(j) +\langle u_j, \omega \rangle=v(\chi^{u_j}(t) \cdot y_j)$ and we deduce 
\begin{equation} \label{orbit equality}
Z =\{\xb \in (\Ycal_{A,a})_s \mid x_j \neq 0 \iff u_j \in A(\yb) \cap Q\}. 
\end{equation}

Next, we prove that $Z$ is a $\tdop_s$-orbit. We have already seen that $Z$ is $\tdop_s$-invariant. It remains to show that the action is transitive and so we consider $z_1,z_2 \in Z$. There is a complete valued field $(F,u)$ extending $(K,v)$ such that $z_1,z_2$ are $\tilde{F}$-rational. Let $L=F((\rdop))$ be the  the Mal'cev-Neumann ring. Note that $L$ is a complete field consisting of certain power series in the variable $x$ and with real exponents (see \cite{Poo} for details).  The advantage is that we have a canonical homomorphism $\rho:\rdop \rightarrow L^*$ with $v \circ \rho = \id$. Using suitable coordinates, we get a homomorphism $N_\rdop \rightarrow T(L)$ which is a section of ${\trop_v}$ and which we also denote by $\rho$. 

For $i=1,2$, there is $t_i\in \Tan$ with $z_i=\pi(t_i \yb)$ and ${\trop_v}(t_i) \in \relint(\sigma)$. Choosing $F$ sufficiently large, we may assume that $t_i$ is induced by an $F$-rational point in $T$ which we also denote by $t_i$. 
For $t \in T(L)$, we set $t^\circ:=t \cdot \rho(-{\trop_v}(t))$. This is an element of the formal affinoid torus $T^\circ(L)$ and hence reduces to an element $\tilde{t^\circ} \in \tdop(\tilde{L})$. The map $t \mapsto t^\circ$ is a homomorphism as well as the reduction. We will use this construction for $t_1$, $t_2$ and $t:=t_2/t_1$. We claim that $\tilde{t^\circ}z_1=z_2$. To see this, we note for $j \in  A(\yb) \cap Q$ that
$$(t^\circ t_1 \yb)_j= \chi^{u_j}(t^\circ t_1) y_j=\chi^{u_j}(\rho\circ {\trop_v}(t_1/t_2)t_2)y_j=\lambda_j (t_2 \yb)_j$$
with factor 
$$\lambda_j:=\chi^{u_j}(\rho \circ {\trop_v}(t_1/t_2))=\rho(\langle u_j, {\trop_v}(t_1)-{\trop_v}(t_2)\rangle).$$
From the above considerations and using that $\trop_v(t_i) \in \relint(\sigma) $, we conclude that $\langle u_j,  {\trop_v}(t_i) \rangle + a(j) $ does not depend on the choice of $u_j \in A(\yb) \cap Q$ and hence the factor $\lambda_j$ does not depend on $j \in  A(\yb) \cap Q$ as well. For $i=1,2$, let $x_i \in \pdop^N$ be the point with coordinates $(t_i \yb)_j$ for $j \in A(\yb) \cap Q$ and with all other coordinates $0$. Then the above shows $t^\circ x_1=x_2$ and hence $\tilde{t^\circ}z_1=z_2$ by the equivariance of the reduction maps. This proves transitivity.



Conversely, if the orbit $Z$ is given, then we may recover $A(\yb) \cap Q$ by \eqref{orbit equality} and this set generates the face $Q$ of $\Wt(\yb,a)$. Let $\sigma = \widehat{Q}$ be the corresponding polyhedron in the dual complex $\Ccal(A,a)$; then 
$${\trop_v}(\{t \in \Tan \mid \pi(t \yb) \in Z\}) \subset \relint(\sigma)$$
is also clear from what we have proven in \eqref{orbit equality}. Then we get immediately equality as the left hand side forms a partition of $N_\rdop$ for varying $Z$. 

The character group of the torus corresponding to the orbit $Z$ has finite index in $\zdop(\sigma^\bot \cap M)$. This is clear as we may choose a base point $\yb'$ in $Z$ and then apply \ref{A-example in the field case}  with $A, \yb'$ and with $\tilde{K}$ replacing $K$. This and \ref{dual complex} prove immediately the identities relating the dimensions. Finally, the claims about the orders are evident.
\qed

\begin{rem} \rm \label{normal fan}
If $v$ is the trivial valuation, then the dual complex $\Ccal(A,a)$ is just the normal fan of the weight polytope $\Wt(\yb)$.
\end{rem}

\begin{cor} \label{duality for generic fibre}
There are bijective order  correspondences between

\begin{itemize}
 \item[(a)] faces $Q$ of the weight polytope $\Wt(\yb)$;
 \item[(b)] polyhedra $\sigma$ of the normal fan of  $\Wt(\yb)$;
 \item[(c)] $T$-orbits $Z$ of the generic fibre of $\Ycal_{A,a}$.
\end{itemize}
\end{cor}

\proof The generic fibre $Y_{A,a}$ is the closure of the $T$-orbit of $\yb$ in  $\pdop^N_K$ with respect to the $T$-action induced by the character set $A$ and so we have seen the claim already in \ref{A-example in the field case}. We note here that the equivalence is also a special case of Proposition   \ref{duality}. To see this, we  replace $v$ by the trivial valuation and then the special fibre is equal to the generic fibre. \qed

\begin{art} \rm \label{weight polytope of orbit}
Let $Z$ be an orbit of $\Ycal_{A,a}$ corresponding to a face $Q$ of the weight subdivision $\Wt(\yb,a)$ (resp. the weight polytope $\Wt(\yb)$). We choose a base point $\zb \in Z(K)$. Then the closure of $Z$ in $\pdop^N$ is the projective toric variety $\Ycal_{A,a}(\zb)$ in $\pdop^N_{\tilde{K}}$ (resp. in $\pdop^N_K$) constructed from $\zb$ and $A(\yb) \cap Q$ as in \ref{A-example in the field case}. We conclude that $Q$ is its weight polytope. 
\end{art}

\begin{rem} \rm \label{preparation for normalization}
The polyhedra of $\Ccal(A,a)$ are pointed if and only if $\Wt(\yb)$ has dimension $n$. In other words, this means that the smallest affine space containing $A(\yb)$ is $N_\rdop$ and this is equivalent to $\dim(\Stab(\yb))=0$ (see Corollary \ref{dimension}).

By passing to a sublattice of $M$, we may always achieve this situation and we may even assume that $M= \zdop A(\yb)$ (see Proposition \ref{affine A-relation}). Since $\Ccal(A,a)$  is a complete complex, it follows from \ref{Burgos-Sombra} that $\Ccal(A,a)= \Sigma_1$ for a complete $\Gamma$-admissible fan $\Sigma$ in $N_\rdop \times \rdop_+$. 
\end{rem}




\section{The Gr\"obner complex} \label{Section: The Groebner complex}

In this section, $K$ denotes a field with a  non-archimedean absolute value $|\phantom{a}|$, corresponding valuation $v:=-\log|\phantom{a}|$, valuation ring $K^\circ$, residue field $\tilde{K}$ and value group $\Gamma = v(K^\times)$. Then we consider a closed subscheme $X$ of the split multiplicative torus $T$ over $K$. We will introduce its Gr\"obner complex on $N_\rdop$ which is related to the natural orbit of $X$ in the Hilbert scheme of a projective compactification. This is a certain complete $\Gamma$-rational complex which has a subcomplex with support equal to $\Trop_v(X)$. At the end, we relate the Gr\"obner complex to the initial degenerations of $X$. This section is inspired by \cite{Ka}, Section 5, which in turn was influenced by Tevelev. We work here with more general base fields, but the ideas are the same. For an elementary approach using Gr\"obner bases and for examples, we refer to \cite{MS}, Section 2.4.


\begin{art} \rm \label{Hilbert scheme}
First, we recall the  following property of the {\it Hilbert scheme} $\Hilb_p(\pdop^m_S)$ for the projective space $\pdop^m_S$ over a locally noetherian 
scheme $S$ and for a Hilbert polynomial $p(x) \in \qdop[x]$. This property  characterizes the Hilbert scheme up to unique isomorphism: 

There is a projective scheme $\Hilb_p(\pdop^m_S)$ over the base scheme $S$ and a closed subscheme $\Univ_p(\pdop^m_S)$ of $\pdop^m_S\times_S \Hilb_p(\pdop^m_S)$ which is flat over $\Hilb_p(\pdop^m_S)$ and which has Hilbert polynomial $p$ such that for every scheme $Z$ over $S$, the map from the set of morphisms $Z \rightarrow \Hilb_p(\pdop^m_S)$ to the set of closed subschemes of $\pdop^m_Z$ with Hilbert polynomial $p$ and flat over $Z$, given by mapping $f$ to the inverse image scheme 
$$(\id \times f)^{-1}(\Univ_p(\pdop^m_S))= \Univ_p(\pdop^m_S) \times_{\Hilb_p(\pdop^m_S)} Z,$$
is a bijection. In other words, there is a bijective correspondence $Y \mapsto [Y]$ between the set of closed subschemes of $\pdop^m_Z$ which are flat over $Z$ and which have Hilbert polynomial $p$ and the set of $Z$-valued points of $\Hilb_p(\pdop^m_S)$. For a proof, we refer to \cite{Kol}, Section 1.1. 
Note that the Hilbert polynomial of a closed subscheme $Y$ of $\pdop^m_Z$ is defined for every fibre over a  point $z$ of $Z$. If $Y$ is flat over $Z$ and if $Z$ is connected, then the Hilbert polynomial does not depend on the choice of $z$. 

The valuation ring $K^\circ$ has not to be noetherian and so we cannot directly apply the above. However, the Hilbert scheme exists also for an arbitrary base scheme $S$ if we require that it represents the functor mapping $Z$ to the set of closed subschemes of $\pdop^m_Z$ with Hilbert polynomial $p$ and which are flat and {\it finitely presented} over $Z$ (see \cite{AK}, Corollary 2.8). Note that if $S=\Spec(K^\circ)$, then every closed subscheme  of $\pdop^m_S$ which is flat over $S$ is  of finite presentation (use \cite{RG}, Corollaire 3.4.7), hence it is defined over a noetherian subring of $K^\circ$  and so we can construct the Hilbert scheme $\Hilb_p(\pdop^m_S)$ from the noetherian case by base change.
\end{art}

\begin{art} \rm \label{construction of Hilbert scheme}
We briefly sketch the construction of the Hilbert scheme as far as we need it later. For simplicity, we restrict to the case $S=\Spec(F)$ for a field $F$. The general case follows similary using graded ideal sheaves instead of graded ideals. For details, we refer to \cite{Kol}, Section 1.1. 

Let $I_Y(k)$ be the $k$-th graded piece of the graded ideal $I_Y$ in $F[x_0, \dots, x_m]$ of a closed subscheme $Y$ of $\pdop^m_F$ with Hilbert polynomial $p$. For sufficiently large $k$ depending only on $p$, we have $\dim(I_Y(k))=q(k)-p(k)$ and the map $Y \mapsto I_Y(k)$ is an injective map from the set of closed subschemes of $\pdop^m_F$ to the Grassmannian $G(q(k)-p(k), q(k))$, where $q$ is the Hilbert polynomial of $\pdop^m$. The image is $\Hilb(\pdop^m_F)$ which we may endow with a suitable  structure as a closed subscheme of the Grassmannian and with a family $\Univ_p(\pdop^m_F)$ which satisfies the required universal property. Using the Grassmann coordinates $L \mapsto \bigwedge^{q(k)-p(k)}(L)$, we get $G(q(k)-p(k),q(k))$ as a closed subscheme of $\pdop^N_F$ for $N:= \binom {q(k)} {p(k)} -1 $ and hence $\Hilb(\pdop^m_F)$ may be seen as a closed subscheme of $\pdop^N_F$ as well. 
\end{art}

\begin{art} \rm \label{torus action on Hilbert scheme}
We consider a linear action of the torus $\tdop$ on $\pdop^m_{K^\circ}$. It follows easily from the universal property of the Hilbert scheme that $\tdop$ operates also on $\Hilb_p(\pdop^m_{K^\circ})$ such that for any scheme $Z$ over $K^\circ$ and any closed subscheme $Y$ of $\pdop^m_Z$ with Hilbert polynomial $p$ which is finitely presented and flat over $Z$, we have $g \cdot [Y]=[g^{-1}Y]$. It makes the following formulas more natural if  the action is by pull-back with respect to multiplication by $g$ rather than push-forward. If we use the closed embedding of $\Hilb_p(\pdop^m_{K^\circ})$ into $\pdop^N_{K^\circ}$ similarly as in \ref{construction of Hilbert scheme}, then the $\tdop$-action on $\Hilb_p(\pdop^m_{K^\circ})$ extends to a linear action of $\tdop$ on $\pdop^N_{K^\circ}$. Indeed, if $A_t$ is the $(m+1)\times (m+1)$-matrix representing the action of $t$ on $\pdop^m_{K^\circ}$ similarly as in the proof of Lemma \ref{lift of projective linear action}, then $(A_t \xb)^\mb$ is a linear combination of monomials of degree $|\mb|$ and this shows easily the claim using the Grassmann coordinates. 
\end{art}

\begin{prop} \label{torus orbit on Hilbert scheme}
Let $Y$ be a closed subscheme of $\pdop^m_K$ with Hilbert polynomial $p$. Then the closure of the $T$-orbit of $[Y]$ in $\Hilb_p(\pdop^m_{K^\circ})$ is equal to $\Ycal_{A,a}$ for  suitable $A \in M^{N+1}$ and height function $a:\{0, \dots , N\} \rightarrow \Gamma \cup \{\infty\}$ using suitable coordinates on $\pdop^N_\kcirc$.
\end{prop}

\proof This follows from Proposition \ref{converse}. \qed

\begin{Def} \rm \label{Groebner complex for Y}
The dual complex $\Ccal(A,a)$  from \ref{dual complex for Y} is called the {\it Gr\"obner complex of $Y$}.  
\end{Def}

\begin{Def} \rm \label{compactified initial degeneration}
Let $(L,w)$ be a valued field extension of $(K,v)$. For $t \in T(L)$,  the special fibre of the closure of $t^{-1}Y_L$ in $\pdop^m_{L^\circ}$ is called the {\it initial degeneration of $Y$ at $t$}. This is a closed subscheme of $\pdop^m_{\tilde{L}}$ defined over the residue field $\tilde{L}$ which we denote  by $\In_t(Y)$. 

For $\omega = \trop_w(t)$, we set $\In_\omega(Y) = \In_t(Y)$. Similarly as in Proposition \ref{dependence on tropicalization}, this is independent of the choice of $t$ up to multiplication by an element from $\tdop$ defined over a suitable field extension of $\tilde{K}$. Since ${\trop_v}$ is surjective,  $\In_\omega(Y)$ is defined for every $\omega \in N_\rdop$. 
\end{Def}

\begin{art} \label{initial degeneration and Hilbert scheme} \rm
In the situation above,  $[t^{-1}Y_L]=t \cdot [Y_L]$ is an  $L$-rational point of $\Hilb(\pdop^m)$. By projectivity of the Hilbert scheme, we conclude that $[t^{-1}Y_L]$ extends uniquely to an  $L^\circ$-valued point $h_t$ of $\Hilb(\pdop^m)$ and hence corresponds to a closed subscheme of $\pdop^m_{L^\circ}$ which is flat over $L^\circ$ and has generic fibre $t^{-1}Y_L$. By Proposition \ref{characterization of affine closure} and Remark \ref{closure in general}, this has to be the closure of $t^{-1}Y_L$ and hence the special fibre is $\In_t(Y)$. In other words, $[\In_t(Y)]$ is equal to the reduction of $h_t$ in $\Hilb(\pdop^m)(\tilde{L})$. 
\end{art}

\begin{prop} \label{Groebner complex and initial degeneration}
Suppose that $\tdop$ acts linearly on $\pdop^m_{K^\circ}$. Let $Y$ be a closed subscheme of $\pdop^m_K$ and let $(L,w)$ be a valued field extending $(K,v)$. For $t_1,t_2 \in T(L)$, the following conditions are equivalent:
\begin{itemize}
\item[(a)] There is a polyhedron $\sigma$ of the Gr\"obner complex $\Ccal(A,a)$ of $Y$ such that for $i=1,2$, we have  $\trop_w(t_i) \in \relint(\sigma)$.
\item[(b)] There is $g \in \tdop(\tilde{L})$ with $\In_{t_2}(Y)=g \cdot  \In_{t_1}(Y)$.
\end{itemize}
\end{prop}

\proof This follows from Proposition \ref{duality}, Proposition \ref{torus orbit on Hilbert scheme} and \ref{initial degeneration and Hilbert scheme}. \qed

\begin{prop} \label{small perturbation of initial degeneration}
Let $\omega_1=\omega_0 + \Delta \omega$ in $N_\rdop$ and suppose that there is a polyhedron $\sigma$ from the Gr\"obner complex $\Ccal(A,a)$ with $\omega_0 \in \sigma$ and $\omega_1 \in \relint(\sigma)$. Then we have
\begin{equation} \label{perturbation}
\In_{\omega_1}(Y)=\In_{\Delta \omega}(\In_{\omega_0}(Y)),
\end{equation}
where  we consider $\In_{\omega_0}(Y)$ as a closed subscheme of $\pdop^m$ over a trivially valued extension of the residue field $\tilde{K}$ and then we take its initial degeneration with respect to $\Delta \omega$. In particular, we have \eqref{perturbation} for all $\omega_1 \in N_\rdop$ in a sufficiently small neighbourhood of $\omega_0$. 
\end{prop}

\proof It follows from  Proposition \ref{Groebner complex and initial degeneration} that $z_1:=[\In_{\omega_1}(Y)]$ is in the orbit $Z_\sigma$ of the special fibre of $\Ycal_{A,a}$ corresponding to $\sigma$. If $\rho$ is the closed face of $\sigma$ with $\omega_0 \in \relint(\rho)$, then $z_0:=[\In_{\omega_0}(Y)]$ is in the orbit $Z:=Z_\rho$ corresponding to $\rho$. 

Now we repeat the procedure taking the  closure of the orbit $Z$ with respect to the base point $z_0$ in $\Hilb(\pdop^m_{\tilde{K}})$.  We have seen in \ref{weight polytope of orbit} that the dual polytope $\hat{\rho}$ is the weight polytope of the projective toric variety $\overline{Z}$. Since we use the trivial valuation on $\tilde{K}$, the dual complex of $\hat{\rho}$ is the complete fan formed by the local cones ${\rm LC}_{\omega_0}(\nu)$ with $\nu$ ranging over all  polyhedra from $\Ccal(A,a)$ containing $\rho$. 
Then $z:=[\In_{\Delta \omega}(\In_{\omega_0}(Y))]$ is in the orbit of $\overline{Z}$ corresponding to the fan ${\rm LC}_{\omega_0}(\nu)$ containing   $\Delta \omega$ in its relative interior. Obviously, this holds for $\nu=\sigma$.

Recall that $\Hilb(\pdop^m_{\tilde{K}})$ is the special fibre of $\Hilb(\pdop^m_{{K}^\circ})$  and  we have $Z=Z_\rho$. Moreover, $\overline{Z}$ is contained in the special fibre of $\Ycal_{A,a}$.  We note that every orbit of $\overline{Z}$  is an orbit of $(\Ycal_{A,a})_s$ and the corresponding fan ${\rm LC}_{\omega_0}(\nu)$ transforms to $\nu$ taking into account that the base point has changed from $[Y]$ to $z_0$.  We conclude that  $z$ and $z_1$ are in the same orbit. This proves \eqref{perturbation}. Finally, the last claim is obvious from the fact that the above local fan in $\omega_1$  is complete. \qed

\begin{art} \rm  \label{cone situation}
In the remaining part of this section, we 
consider the following important special case, where we can compare the definitions in \ref{compactified initial degeneration} and in \ref{initial degeneration at t}: We consider a projective toric variety $\Ycal_{B,0}$ over $K^\circ$ given by $B \in M^{m+1}$ and height function identically zero, i.e. the base point $\zb \in \pdop^m(K)$ in the open dense orbit satisfies $v(z_j)=0$ for $j=0,\dots ,m$. Recall that $\Ycal_{B,0}$ is a closed subvariety of $\pdop^m_{K^\circ}$ and the torus action extends to a linear action on $\pdop^m_{K^\circ}$ (see \ref{A-example}). 
We assume further that the stabilizer of $\zb$ is trivial and so we may identify $T$ with the open dense orbit $T\zb$. By Corollary  \ref{dimension},  the affine span of $B$ is $M_\rdop$. For example, the standard embedding of $\tdop$ in $\pdop^n_{K^\circ}$ fulfills all these requirements.

The triviality of the height function implies that the weight polytope is equal to the weight subdivision and the dual complex is just the normal fan of $\Wt(\zb)$. Moreover, we may identify $\tdop$ with the $\tdop$-invariant open subset of $\Ycal_{B,0}$ whose generic fibre is the open dense orbit and whose special fibre is the orbit corresponding to the vertex $0$ of the cones.  
\end{art}

\begin{art} \label{comparision of initial degenerations} \rm
We consider a closed subscheme $X$ of $T$ and we denote by $Y$ its closure in $\pdop^m_K$. For a valued field $(L,w)$ extending $(K,v)$ and $t \in T(L)$, it follows immediately from comparing Definitions \ref{initial degeneration at t} and \ref{compactified initial degeneration} that $$\In_{t}(X_L)=\In_t(Y_L) \cap \tdop_{\tilde{L}}.$$
\end{art}

\begin{cor} \label{toric perturbation}
For $\omega_0 \in N_\rdop$, there is a neighbourhood $\Omega$ of $\omega_0$ in $N_\rdop$ such that
\begin{equation*} 
\In_{\omega_1}(X)=\In_{\Delta \omega}(\In_{\omega_0}(X)),
\end{equation*}
for every $\omega_1 \in \Omega$ and $\Delta \omega := \omega_1 - \omega_0$. On the right hand side, the initial degeneration of $\In_{\omega_0}(X)$ at $\Delta \omega$ is with respect to a trivially valued field of definition for  $\In_{\omega_0}(X)$.
\end{cor}

\proof This follows from Proposition \ref{small perturbation of initial degeneration} and \ref{comparision of initial degenerations}. \qed

\begin{art} \label{Groebner complex for X} \rm 
We apply the above for $Y=\overline{X}$ leading to a  polyhedral complex  $\Ccal(A,a)$ in $N_\rdop$ which we call the  {\it Gr\"obner complex for $X$}. It depends on the choices from \ref{cone situation}. 
\end{art}

\begin{thm} \label{support of Groebner}
The Gr\"obner complex $\Ccal(A,a)$ of $X$ is a complete $\Gamma$-rational complex in $N_\rdop$ and $\{\sigma \in \Ccal(A,a) \mid \sigma \subset \Trop_v(X)\}$ is a subcomplex $\Ccal_X$ of $\Ccal(A,a)$ with support equal to $\Trop_v(X)$.
\end{thm}

\proof All statements are evident by construction except the claim about the support. Let $\omega \in \Trop_v(X)$. By completeness of the Gr\"obner complex, there is $\sigma \in \Ccal(A,a)$ with $\omega \in \relint(\sigma)$. We have to prove that $\sigma \subset \Trop_v(X)$. Since $\Trop_v(X)$ is closed in $N_\rdop$, it is enough to show that every $\omega' \in \relint(\sigma)$ is contained in $\Trop_v(X)$. There is a valued field  $(L,w)$ extending $(K,v)$ and $t, t' \in T(L)$ with $\trop_w(t)= \omega$ and $\trop_w(t')=\omega'$. By Proposition \ref{Groebner complex and initial degeneration}, there is $g \in \tdop(\tilde{L})$ with $\In_{t'}(Y)=g \cdot \In_{t}(Y)$. By \ref{comparision of initial degenerations} , we conclude that $\In_{\omega'}(X)=  \In_\omega(X) $. Using $\omega \in \Trop_v(X)$, Theorem \ref{fundamental theorem} implies that $\In_\omega(X)$ is non-empty and hence the same is true for $\In_{\omega'}(X)$. Using this equivalence the other way round, we deduce that 
$\omega' \in \Trop_v(X)$ proving the claim. \qed

\vspace{2mm}
The following result is very useful for reducing local statements about the tropical variety to the case of trivial valuations. We will see in Proposition \ref{tropical multiplicities are local} that this is also compatible with tropical multiplicities.

\begin{prop} \label{trop trop}
Let $X$ be a closed subscheme of $T$ and let $\omega \in N_\rdop$. Using the local cone at $\omega$ from Appendix \ref{local cone},  we have
$$\Trop_0(\In_\omega(X))={\rm LC}_\omega(\Trop_v(X)).$$
\end{prop}

\proof 
The fundamental theorem of tropical algebraic geometry (Theorem \ref{fundamental theorem}) says that $ \Delta\omega\in N_\rdop$ is in $\Trop_0(\In_\omega(X))$ if and only if $\In_{\Delta \omega}(\In_\omega(X))$ is non-empty. If we choose $\Delta \omega$ sufficiently small, then we deduce from Corollary \ref{toric perturbation} that these conditions are also equivalent to $\In_{\omega + \Delta \omega}(X) \neq \emptyset$. Theorem \ref{fundamental theorem} again shows that this is equivalent to $\omega + \Delta \omega \in \Trop_v(X)$. 
As we are working in a sufficiently small neighbourhood of $\omega$, this is equivalent to $\Delta \omega \in {\rm LC}_\omega(\Trop_v(X))$ proving the claim. \qed

\begin{art}  \rm \label{Groebner fan}
For a polyhedron $\Delta$ in $N_\rdop$, let us recall that $c(\Delta)$ denotes  the closed cone in $N_\rdop \times \rdop_+$ generated by $\Delta \times \{1\}$. 
We call  $\Sigma(A,a):=c(\Ccal(A,a))$ the {\it  Gr\"obner fan of $X$ in $N_\rdop \times \rdop_+$}. 
\end{art}

The following result uses that the tropical cone $\Trop_W(X)$ is the closure of the cone in $N_\rdop \times \rdop_+$ generated by $\Trop_v(X) \times \{1\}$ which is proved in the next section (see Corollary \ref{cone property proved}). We will not use Corollary \ref{support of Groebner fan} and the following consequence in \ref{pointed subdivision} in the next section.

\begin{cor}  \label{support of Groebner fan}
The Gr\"obner fan $\Sigma(A,a)$ of $X$ in $N_\rdop \times \rdop_+$ is a complete $\Gamma$-rational fan  and $\Sigma_X:=\{\sigma \in \Sigma(A,a) \mid \sigma \subset \Trop_W(X)\}$ is a subfan  of $\Sigma(A,a)$ with support equal to the tropical cone $\Trop_W(X)$ from \ref{tropical cone of X}.
\end{cor}

\proof Since $\Ccal(A,a)$ is a complete $\Gamma$-rational polyhedral complex, it follows from Remark \ref{Burgos-Sombra} that $\Sigma(A,a)$ is a complete $\Gamma$-rational fan in $N_\rdop \times \rdop_+$. Then the claim  follows from Theorem \ref{support of Groebner} and Corollary \ref{cone property proved} below. \qed

\begin{art} \rm \label{pointed subdivision}
By \ref{preparation for normalization}, $\Ccal(A,a)$ is a pointed polyhedral complex if and only if $\Stab(\yb)$ is zero-dimensional. By definition of the torus action on the Hilbert scheme, we have $\Stab(\yb)=\Stab(Y)=\Stab(X)$, where $Y$ is the closure of $X$ in $\pdop^m_K$. In general, it is clear that $\Ccal(A,a)$ is isomorphic to the product of an affine  space and  the Gr\"obner complex of $X/\Stab(X)$. By the above, the latter is pointed and so it is obvious that $\Ccal(A,a)$ has always a $\Gamma$-rational subdivision $\Ccal$ consisting of pointed polyhedra. By Corollary \ref{support of Groebner fan}, $c(\Ccal)$ is a $\Gamma$-admissible fan in $N_\rdop \times \rdop_+$ with support $\Trop_W(X)$.
\end{art}






\section{Compactifications in toric schemes} \label{Compactifications in toric schemes}

Let  $K$ be a field with a  non-archimedean absolute value $|\phantom{a}|$, corresponding valuation $v:=-\log|\phantom{a}|$, valuation ring $K^\circ$, residue field $\tilde{K}$ and value group $\Gamma = v(K^\times)$. Let $\tdop$ be the split torus over ${K^\circ}$ with generic fibre $T$ associated to the character lattice $M$ of rank $n$ and dual lattice $N$. We keep the usual notation. In this section, we consider a closed subscheme $X$ of $T$ and we study its closure $\Xcal$ in the toric scheme $\Ycal_\Sigma$ associated to a $\Gamma$-admissible fan $\Sigma$ in $N_\rdop \times \rdop_+$ (see \ref{admissible fans}).  First, we prove surjectivity of the reduction map which is called the tropical lifting lemma. Then, we show Tevelev's lemma which is a tropical characterization of the orbits intersecting $\Xcal$. Finally, we give several equivalences for properness of the occurring schemes.

\vspace{2mm}

We start with a lemma due to Draisma.

\begin{lem} \label{Draisma's Lemma}
Let $(L,w)$ be a valued field extending $(K,v)$ and let $r,s \in \ndop$. For $a_{ij}, b_i \in K$ and $\lambda_i \in \rdop$, we consider the following system of equalities
$$a_{i1}x_1+ \dots 
+ a_{it}x_t = b_i \quad (1 \leq i \leq r)$$
and inequalities 
$$w(a_{i1}x_1+ \dots + a_{it}x_t) \geq \lambda_i \quad (r+1 \leq i \leq r+s).$$
If this system has a solution $y \in L^t$, then it has also a solution $z \in K^t$.
\end{lem}

\proof This follows from the same arguments as Lemma 4.3 in \cite{Dr}. \qed

\begin{lem} \label{base change and tropical compactification}
Let $(L,w)$ be a valued field extending $(K,v)$ and let $\Xcal'$ be the closure of $X_L$ in the toric scheme over $L^\circ$ associated to $\Sigma$. Then 
the canonical morphism $\phi:(\Xcal')_s \rightarrow \Xcal_s$ is surjective.
\end{lem}

\proof We will first prove the claim if the value group $\Gamma$ is a divisible subgroup of $\rdop$ and then we will reduce the claim to this special case in several steps.

\vspace{2mm}
\noindent {\it Step 1: If the value group $\Gamma$ is a divisible subgroup of $\rdop$, then $\phi$ is surjective.}
\vspace{2mm}

In this case, we have seen in Proposition \ref{divisible case} that the toric scheme over $L^\circ$ associated to $\Sigma$ is the base change of $\Ycal_\Sigma$ to $L^\circ$. By Corollary \ref{closure and base change}, we have $\Xcal'=\Xcal_{L^\circ}$ and hence $(\Xcal')_s$ is the base change of $\Xcal_s$ to the residue field $\tilde{L}$. This yields surjectivity of $\phi$.

\vspace{2mm}

In particular, this proves the claim for $v$ trivial. We may assume that $v$ is non-trivial and that $\Ycal_\Sigma = \Ucal_\Delta$ for a pointed $\Gamma$-rational polyhedron $\Delta$ in $N_\rdop$.  Let $\sigma$ be the recession cone of $\Delta$. Then $X$ is given by an ideal $I_X$ in $K[M]^{\sigma}$ and its closure $\Xcal$ is given by the ideal $I_X \cap K[M]^\Delta$ in $K[M]^\Delta$. Similarly, $\Xcal'$ is the closed subscheme given by the ideal $(I_X L[M]^\sigma)\cap L[M]^\Delta$ in $L[M]^\Delta$. 


\vspace{2mm}
\noindent {\it Step 2: The morphism $\phi$ is dominant.}
\vspace{2mm}

Let $f \in K[M]^\Delta$ such that the residue class of $f$ in $L[M]^\Delta/ ((I_X L[M]^\sigma)\cap L[M]^\Delta) )\otimes_{\tilde{K}} \tilde{L}$ is zero. 
We have to prove that there is $m \in \ndop$ such that $f^m \in (I_X \cap K[M]^\Delta) + K^{\circ \circ} K[M]^\Delta$. By assumption, we have 
\begin{equation} \label{equation for f}
f= g_1h_1 + \dots + g_r h_r + \lambda f_1 + \dots +\lambda f_s
\end{equation}
with $g_i \in I_X$, $h_i \in L[M]^\sigma$,  $\lambda \in L^{\circ \circ}$ and $f_j \in L[M]^\Delta$. 
We may assume that $h_i = \beta_i \chi^{u_i}$ for some $\beta_i \in L$ and $u_i \in \check{\sigma} \cap M$. Similarly, we may assume that $f_j = \gamma_j \chi^{v_j}$ for some $\gamma_j \in L$ and $v_j \in M$. Since the valuation $v$ is non-trivial, there is $m \in \ndop$ such that $\lambda^m$ is divisible by an element of $K^\circ$. Replacing $f$ by $f^m$, we may assume that $\lambda \in K^{\circ \circ}$. If we compare the coefficients on both sides of equation \eqref{equation for f}, then we get a finite system of linear equations with coefficients in $K$ and unknowns $\beta_1, \dots , \beta_r$ and $\gamma_1, \dots, \gamma_s$. The conditions $f_j \in L[M]^\Delta$ are equivalent to the finitely many inequalities $v(\gamma_j) + \langle v_j, \omega \rangle \geq 0$, where $\omega$ ranges over the vertices of $\Delta$. By assumption, this system of equalities and inequalities has a solution in $L^{r+s}$. By Lemma \ref{Draisma's Lemma}, there is a solution with $\beta_1, \dots, \beta_r,\gamma_1, \dots, \gamma_s \in K$ which means that we find a representation in \eqref{equation for f} with all $h_i \in K[M]^\sigma$ and all $f_j \in K[M]^\Delta$. We conclude that $f \in (I_X \cap K[M]^\Delta) + K^{\circ \circ} K[M]^\Delta$ proving Step 2. 

\vspace{2mm}
\noindent {\it Step 3: If $L$ is an algebraic closure of $K$, then $\phi$ induces a finite surjective map $(\Xcal')_s \rightarrow \Xcal_s \otimes_{\tilde{K}} \tilde{L}$.}
\vspace{2mm}

We use first that the value group $\Gamma_L$ of $w$ is equal to $\{ \lambda \in \rdop \mid \exists m \in \ndop \setminus \{0\}, \, m \lambda \in \Gamma\}$. It follows  that the vertices of $\Delta$ are in $N_{\Gamma_L}$ and there is a non-zero $m \in \ndop$ such that $m \omega \in N_\Gamma$ for every vertex $\omega$ of $\Delta$. For every $u \in \check{\sigma} \cap M$, there is $\beta_u \in L$ with $v_\Delta(\beta_u \chi^u)=0$. 
For each vertex $\omega$ of $\Delta$, we choose a finite generating set of the semigroup $\check{\sigma}_\omega \cap M$, where $\sigma_\omega$ is the local cone of $\Delta$ at $\omega$. We have seen in the proof of Proposition \ref{finite generation} that $L[M]^\Delta$ is generated as an $L^\circ$-algebra by $\beta_u \chi^u$, where  $u$ ranges over the union $S$  of all these generating sets. 

We claim that the finite set $H:=\{\prod_{u \in S}(\beta_u \chi^u)^{k_u} \mid 0 \leq k_u < m\}$ generates $L[M]^\Delta$ as a $K[M]^\Delta \otimes_{K^\circ} L^\circ$-module. Indeed, every $f \in L[M]^\Delta$ has the form $f=\sum_{h,k} \lambda_{hk}h  \prod_{u \in S}(\beta_u \chi^u)^{mk_u} $ where $h$ ranges over $H$, $k$ over $\ndop^{S}$ and only finitely many coefficients $\lambda_{hk} \in L^\circ$ are non-zero. The construction of $m$ yields that $\langle m u, \omega \rangle \in \Gamma$ for every vertex $\omega$ of $\Delta$ and hence there is $\alpha_u \in K$ with $v_\Delta(\alpha_u \chi^{m u})=0$. We conclude that $\beta_u^{m}=\alpha_u \gamma_u$ for some $ \gamma_u \in L^\circ$. Since $\alpha_u \chi^{m u} \in K[M]^\Delta$, this implies that $ \prod_{u \in S}(\beta_u \chi^u)^{mk_u}\in L^\circ K[M]^\Delta$ proving that $H$ generates the module $L[M]^\Delta$.

Since $(\Xcal')_s$ is a closed subscheme of $\Spec(L[M]^\Delta)$ and since $\Xcal_s \otimes_{\tilde{K}} \tilde{L}$ is a closed subscheme of $\Spec(K[M]^\Delta \otimes_{K^\circ} L^\circ)$, we conclude that $(\Xcal')_s \rightarrow \Xcal_s \otimes_{\tilde{K}} \tilde{L}$ is a finite map. It follows from Step 2 that this map is dominant and hence it is surjective proving Step 3.

\vspace{2mm}

We will now deduce the claim from Step 3. We endow an algebraic closure $E$ of $L$ with a valuation $u$ extending $w$. Let $F$ be the algebraic closure of $K$ in $E$ endowed with the restriction of $u$. Let $\Xcal''$ (resp. $\Xcal'''$) be the closure of $X_F$ (resp. $X_E$) in the toric scheme over $F^\circ$ (resp. $E^\circ$) associated to $\Delta$. Then we have a commutative diagram
\begin{equation*} 
\begin{CD}
(\Xcal''')_s @> >> \Xcal_s \otimes_{\tilde{K}} \tilde{F}\\
@VV V    @VV V\\
(\Xcal')_s@>{\phi}>> \Xcal_s
\end{CD}
\end{equation*}
of canonical morphisms. The first row has the factorization $(\Xcal''')_s \rightarrow(\Xcal'')_s \rightarrow \Xcal_s \otimes_{\tilde{K}} \tilde{F}$ and hence it is surjective by Steps 1 and 3. Since the second column is  surjective as well, we  deduce that $\phi$ is surjective. \qed

\begin{prop} \label{purity of special fibre}
The closure $\Xcal$ of $X$ in $\Ycal_\Sigma$ is a separated flat scheme over $\kcirc$. The special fibre $\Xcal_s$ is of finite type over $\tilde{K}$.  If  $X$ is of pure dimension $d$ and $\Xcal_s$ is non-empty, then $\Xcal_s$ is of pure dimension $d$. If the value group $\Gamma$ is divisible, then $\Xcal$ is of finite presentation over $\kcirc$.
\end{prop}

\proof It follows from Lemma \ref{separatedness} that $\Xcal$ is separated over $\kcirc$. Flatness is a consequence of  \ref{polyhedral scheme} and \ref{closure in general}. By Lemma \ref{reduced special fibre and omega-weight},    $\Xcal_s$ is  of finite type over $\ktilde$. The dimensionality claim is clear in case of a divisible value group $\Gamma$ as in this case $\Xcal$ is a flat scheme of finite type over $K^\circ$ (see Proposition \ref{finite generation}). In this case, it follows from \cite{RG}, Corollaire 3.4.7, that $\Xcal$ is of finite presentation over $\kcirc$. In general, we will reduce to the divisible case: We may assume that $v$ is non-trivial and that $\Ycal_\Sigma = \Ucal_\Delta$ for a $\Gamma$-rational polyhedron $\Delta$ in $N_\rdop$. Then the dimensionality claim follows from Step 3 in the proof of Proposition \ref{base change and tropical compactification}.  \qed


\begin{art} \rm \label{recall reduction}
We recall from \S \ref{algebraic and formal models} that the reduction map is defined on an analytic subdomain $(Y_{\Sigma_0})^\circ$ of the generic fibre $Y_{\Sigma_0}$ and maps to the special fibre of the $K^\circ$-model $\Ycal_\Sigma$. The points of $(Y_{\Sigma_0})^\circ$ are induced by potentially integral points and Proposition \ref{integral points of toric variety} shows that $(Y_{\Sigma_0})^\circ \cap \Tan = \trop_v^{-1}(|\Sigma_1|)$. We conclude that the potentially integral points of $X$ with respect to  $\Xcal$  induce an analytic subdomain $X^\circ = \trop_v^{-1}(|\Sigma|) \cap \Xan$ of $\Xan$ where we have a well-defined reduction map $\pi:X^\circ \rightarrow \Xcal_s$. 
\end{art}

We have here the following generalization of  Jan Draisma's tropical lifting lemma (see \cite{Dr}, Lemma 4.4). 

\begin{prop}   \label{tropical lifting lemma}
Using the above notation, we have $\pi(\Uan \cap X^\circ) = \Xcal_s$ for every open dense subset $U$ of $X$. Moreover, if $K$ is algebraically closed and $v$ is non-trivial, then every closed point of $\Xcal_s$ is the reduction of a closed point of $U$. 
\end{prop}

\proof The additional difficulty here in contrast to Draisma's paper is that $\Xcal$ and the ambient toric scheme $\Ycal_\Sigma$ might be not of finite type (see Proposition \ref{finite presentation in non-discrete case}). Let $L$ be an algebraic closure of $K$ and let us choose a valuation $u$ on $L$ extending $v$. Let $\Xcal'$ be the closure of $X_L$ in the toric scheme over $L^\circ$ associated to the fan $\Sigma$. Then  $\Xcal'$ is a flat scheme of finite type over $L^\circ$ by Proposition \ref{purity of special fibre}. By Proposition \ref{finite generated lifting lemma}, the reduction map $\pi_L:(U_L)^{\rm an} \cap (X_L)^\circ \rightarrow (\Xcal')_s$ is surjective. We have a canonical commutative diagram
\begin{equation*} 
\begin{CD}
(U_L)^{\rm an} \cap (X_L)^\circ @>{\pi_L}>> \Xcal_s' \\
@VV V    @VV V\\
\Uan \cap X^\circ @>{\pi}>> \Xcal_s
\end{CD}
\end{equation*}
where  the second column is surjective by Lemma \ref{base change and tropical compactification}. This proves surjectivity of $\pi$. The last claim follows directly from Proposition \ref{finite generated lifting lemma}. \qed

\vspace{2mm}

The following result is called {\it Tevelev's Lemma}. We will use the tropical cone $\Trop_W(X)$ and the notation of the previous section. The bijective correspondence between open faces and orbits from Proposition \ref{uniform orbit correspondence} will be important for the understanding of the following.



\begin{lem} \label{Tevelev's Lemma}
Let $\sigma \in \Sigma$. Then the orbit $Z_\sigma$ corresponding to $\relint(\sigma)$ intersects $\Xcal$ if and only if $\Trop_W(X) \cap \relint(\sigma)$ is non-empty. 
\end{lem}

\proof If $\omega \in \Trop_W(X) \cap \relint(\sigma)$, then there is $x \in X_W^{\rm an}$ with  $\trop_W(x)= \omega$.  Let $\pi_W: \trop_W^{-1}(|\Sigma|) \rightarrow \Ycal_\Sigma$ be the reduction map. We deduce from Proposition \ref{uniform orbit correspondence}  that $\pi_W(x) \in Z_\sigma$. Since $\pi_W(x)$ is also contained in $\Xcal$, we see that $\Xcal \cap Z_\sigma$  is non-empty.

Suppose that $z \in \Xcal \cap Z_\sigma$. By Proposition \ref{tropical lifting lemma}, $\pi_W$ induces a surjective map $X_W^{\rm an} \cap \trop_W^{-1}(|\Sigma|)  \rightarrow \Xcal$ and hence there is $x \in X_W^{\rm an}$ with $z= \pi_W(x)$. Again by Proposition \ref{uniform orbit correspondence}, we see that $\trop_W(x) \in \relint(\sigma)$. 
\qed

\begin{rem} \label{reduction to trivial valuation} \rm
We will give a procedure which can be often used to reduce questions about  $\Xcal$  to the case of the trivial valuation:

Let $\omega$ be a vertex of $\Sigma_1$. By \ref{tropicalization and orbits} or Proposition \ref{irreducible components of the special fibre}, we have a corresponding irreducible component $Y_\omega$ of the toric scheme $\Ycal_\Sigma$. It is the closure of the orbit $Z_\omega$ corresponding to the vertex $\omega$. By Proposition \ref{irreducible component is toric},    $Y_\omega$ may be viewed as a normal toric variety over $\tilde{K}$ associated to the fan ${\rm LC}_\omega(\Sigma_1):=\{{\rm LC}_\omega(\Delta) \mid \Delta \in \Sigma_1\}$. Note that the acting torus is $T_\omega:=\Spec(\tilde{K}[M_\omega])$, where $M_\omega :=\{ u \in M \mid \langle u, \omega \rangle \in \Gamma\}$ is a sublattice of $M$ of finite index. To identify  it with the dense open orbit $Z_\omega$ of $Y_\omega$ involves the  choice of a basepoint in $Z_\omega(\tilde{K})$ which does not influence tropical varieties of closed subschemes of $Z_\omega(\tilde{K})$ as we deal with the trivial valuation on $\tilde{K}$.

We assume now that the vertex $\omega$ is also contained in  $\Trop_v(X)$. In local problems involving $\omega$, the relevant closed subscheme of $T_\omega$ is $X_\omega := \Xcal \cap Z_\omega$. By Tevelev's Lemma \ref{Tevelev's Lemma}, $X_\omega$ is non-empty and its closure $\Xcal_\omega$ is contained in $\Xcal \cap Y_\omega$. We claim that $\Trop_0(X_\omega)$ is the local cone of $\Trop_v(X)$ at $\omega$.  

To prove the claim, we note first that the induced reduced structure of the special fibre is compatible with base change by Lemma \ref{reduced special fibre and omega-weight}. As the tropical variety is also invariant under base change (Proposition \ref{tropical base change}), we may assume that $\omega \in N_\Gamma$. Then there is $t \in T(K)$ with $\omega = \trop_v(t)$ and we may choose the basepoint of $Z_\omega$ equal to $\pi(t)$. Using translation by $t^{-1}$, we conclude easily that $X_\omega$ is isomorphic to $\In_\omega(X)$ and hence the claim follows from Proposition \ref{trop trop}.
\end{rem}

\begin{prop}  \label{properness of toric schemes}
For a $\Gamma$-admissible fan $\Sigma$ in $N_\rdop \times \rdop_+$, the following conditions are equivalent:
\begin{itemize}
\item[(a)] $|\Sigma | = N_\rdop \times \rdop_+$;
\item[(b)] $|\Sigma_1| = N_\rdop$;
\item[(c)] the special fibre of $\Ycal_\Sigma$ is non-empty and proper over $\tilde{K}$;
\item[(d)] $\Ycal_\Sigma$ is universally closed over $K^\circ$.
\end{itemize}
If the equivalent conditions (a)--(d) hold, then the generic fibre of $\Ycal_\Sigma$ is also proper over $K$. If the value group $\Gamma$ is  divisible or discrete in $\rdop$, then (a)--(d) are equivalent to $\Ycal_\Sigma$ proper over $K^\circ$.
\end{prop}

\proof Clearly, (a) and (b) are equivalent. Suppose that (a) holds. Then $\Sigma_0$ is a complete fan  and hence the generic fibre $Y_{\Sigma_0}$ of $\Ycal_\Sigma$ is complete (see \cite{Fu2}, \S 2.4). The special fibre of $\Ycal_\Sigma$ is the union of its finitely many irreducible components corresponding to the vertices $\omega_j$ of $\Sigma_1$ (see \ref{tropicalization and orbits}). Such an irreducible component is a  toric variety with fan ${\rm LC}_\omega(\Sigma_1)$ generated by the local cones in $\omega_j$ (see Proposition  \ref{special orbit closure}). Since $\Sigma_1$ satisfies (b), all these fans are also complete and hence the irreducible components are proper over $\tilde{K}$. Proposition \ref{purity of special fibre} shows that the special fibre is separated and of finite type over $\ktilde$. We conclude that $(\Ycal_\Sigma)_s$ is proper over $\ktilde$ by \cite{EGA IV}, Corollaire 5.4.5. This proves (a) $\Rightarrow$ (c). 

Next, we show that (c) implies (b). If the special fibre  $(\Ycal_\Sigma)_s$ is proper over $\tilde{K}$, then every irreducible component of  $(\Ycal_\Sigma)_s$ is complete. As we have seen above, such an irreducible component $Y$  is associated to a vertex $\omega$ of $\Sigma_1$ and $Y$ is a toric variety over $\tilde{K}$ associated to the fan ${\rm LC}_\omega(\Sigma_1)$. By completeness of $Y$, the fan   ${\rm LC}_\omega(\Sigma_1)$ is also complete (see \cite{Fu2}, \S 2.4) . As this holds for every vertex of $\Sigma_1$, we conclude that $\Sigma_1$ is complete. This proves (c) $\Rightarrow$ (b).


We have seen now that (a)--(c) are equivalent and that the generic fibre of $\Ycal_\Sigma$ is proper over $K$  in this case. Now we will prove that (a) yields (d). So we assume that the equivalent conditions (a)--(c) hold. We note that $\Ycal_\Sigma$ is a quasi-compact separated scheme over $K^\circ$ and so we may apply the  valuative criterion of universal closedness (\cite{EGA II}, Th\'eor\`eme 7.3.8) which holds also in the non-noetherian situation. Let $L^\circ$ be a valuation ring with fraction field $L$ and suppose that we have a commutative diagram
\begin{equation*} 
\begin{CD}
\Spec(L) @>{\psi}>>\Ycal_\Sigma \\
@VV V    @VV V\\
\Spec(L^\circ) @>>> \Spec(K^\circ)
\end{CD}
\end{equation*}
of morphisms. To prove universal closedness of $\Ycal_\Sigma$, the criterion says that it is enough to show that there is a morphism $g: \Spec(L^\circ) \rightarrow \Ycal_\Sigma$ over $K^\circ$ which factors through $\psi$. First, we assume that the homomorphism $K^\circ \rightarrow L^\circ$ is not injective. Then the kernel is the maximal ideal $K^{\circ \circ}$ of $K^\circ$ and $\psi$ factors through the special fibre of $\Ycal_\Sigma$. It follows from  (c) that the special fibre is proper over $\tilde{K}$ and hence the  valuative criterion of properness gives the existence of $g$. So we may assume that $K^\circ \subset L^\circ$. Then the intersection of the maximal ideal $L^{\circ \circ}$ of $L^\circ$ with $K^\circ$ is either $\{0\}$ or $K^{\circ \circ}$. In the first case, we may replace the second column in the diagram by its generic fibre. We have seen at the beginning of the proof that the generic fibre of $\Ycal_\Sigma$ is proper over $K$. By the  valuative criterion of properness again, we get the existence of $g$. 

Finally, we have to consider the case $L^{\circ \circ} \cap K^\circ = K^{\circ \circ}$. In this case, $L^\circ$ is the valuation ring of a valuation $w$ on $L$ extending $v$. Since the value group $\Gamma_L$ of $w$ is a totally ordered abelian group, the polyhedra $\Delta \in \Sigma_1$ induce polyhedra $\Delta(\Gamma_L)$ in $N_{\Gamma_L}$ and (b) yields that they are covering the whole space (see \ref{compatibility properties}). We note that $\psi$ corresponds to an $L$-rational point $P$ of the generic fibre  $(\Ycal_\Sigma)_\eta=Y_{\Sigma_0}$ and we have to show that $P$ is an $L^\circ$-point of $\Ycal_\Sigma$. If $P$ is contained in the dense orbit $T$, then $\trop_w(P) \in \Delta(\Gamma_L)$ for some $\Delta \in \Sigma_1$ and it follows from Proposition \ref{integral points of toric variety} that $P$ extends to an $L^\circ$-integral point of $\Ycal_\Sigma$. If $P$ is contained in another orbit $Z$ of $Y_{\Sigma_0}$, then this extension property holds as well by using Proposition \ref{generic orbit closure} to reduce to the previous case. This proves that $\Ycal_\Sigma$ is universally closed over $K^\circ$ and hence we get (a)--(c) $\Rightarrow$ (d).

Conversely, if (d) holds, then the special fibre is non-empty and universally closed over $\ktilde$. By Proposition \ref{purity of special fibre}, the special fibre is also separated and of finite type over $\ktilde$. We conclude that $(\Ycal_\Sigma)_s$ is proper over $\ktilde$. This proves (d) $\Rightarrow$ (c) and hence all four properties are equivalent. If the value group $\Gamma$ is divisible or discrete in $\rdop$, Proposition \ref{purity of special fibre} shows that $\Ycal_\Sigma$ is separated and of finite type over $K^\circ$. If (d) holds, then it follows that $\Ycal_\Sigma$ is proper over $K^\circ$. \qed

\begin{art} \rm \label{toric morphism}
We consider now another free abelian group $M'$ of finite rank with dual $N'$ and split multiplicative torus $T' = \Spec(K[M'])$ over $K$. Then a $\Gamma$-admissible fan $\Sigma'$ in $N'_\rdop \times \rdop_+$ induces a toric scheme $\Ycal_{\Sigma'}$ over $K^\circ$ with dense orbit $T'$. We assume that $f:N' \rightarrow N$ is a homomorphism such that $f_\rdop \times \id_{\rdop_+}$ maps each cone $\sigma'$ of $\Sigma'$ into  a suitable cone $\sigma$ of $\Sigma$. Then the dual homomorphism of $f$ induces a canonical equivariant morphism $\Vcal_{\sigma'} \rightarrow \Vcal_\sigma$. We can patch these homomorphisms together to get an equivariant morphism $\varphi:\Ycal_{\Sigma'} \rightarrow \Ycal_\Sigma$ of toric schemes over $K^\circ$ which is canonically determined by $f$ through the fact that $\varphi$ restricts to the homomorphism $T' \rightarrow T$ of tori induced by $f$.
\end{art}

\begin{prop} \label{proper morphism}
Under the hypothesis above, the following properties are equivalent:
\begin{itemize}
 \item[(a)] the morphism  $\varphi$ is closed with  generic fibre $\varphi_\eta$ and special fibre $\varphi_s$ both proper;
 \item[(b)] $(f_\rdop \times \id_{\rdop_+})^{-1}(|\Sigma|)=|\Sigma'|$;
 \item[(c)] the morphism $\varphi$ is universally closed.
\end{itemize}
If the value group $\Gamma$ is  divisible or discrete in $\rdop$, then (a)--(c) are also equivalent to $\varphi$ proper over $K^\circ$. 
\end{prop}

\proof We assume that (a) holds. By the criterion of properness for homomorphisms of toric varieties over a field (\cite{Fu2}, \S 2.4), we have $f^{-1}(|\Sigma_0|)=|\Sigma_0'|$. To prove (b), it remains to see that $f^{-1}(|\Sigma_1|)=|\Sigma_1'|$. 
Let $\omega' \in N_\rdop$ with $f(\omega') \in |\Sigma_1|$. There is $t' \in (T')^{\rm an}$ with $\trop_v(t')=\omega'$ and hence $\trop_v(t)=f(\omega') \in |\Sigma_1|$ for $t := \varphi^{\rm an}(t')$. By Proposition \ref{integral points of toric variety}, we have $t \in Y_{\Sigma_0}^\circ$ and hence we have a well-defined reduction $\pi(t)$ in the special fibre of the closure of $\varphi(T')$ in $\Ycal_\Sigma$. Since $\varphi$ is closed, we have $\pi(t)= \varphi_s(z')$ for some $z' \in (\Ycal_{\Sigma'})_s$. By Proposition \ref{tropical lifting lemma}, there is $t_0' \in (T')^{\rm an} \cap Y_{\Sigma_0'}^\circ$ with reduction $\pi(t_0')=z'$. Again Proposition \ref{integral points of toric variety} shows that $\omega_0' := \trop_v(t_0') \in \Sigma_1'$. We have $\pi(\varphi^{\rm an}(t_0')) = \varphi_s(\pi(t_0'))=\pi(t)$. The orbit correspondence in Proposition \ref{uniform orbit correspondence} yields that $f(\omega_0') = \trop_v(\varphi^{\rm an}(t_0'))$ is in the same open face $\tau$ of $\Sigma_1$ as $f(\omega')=\trop_v(t)$.

 Arguing by contradiction, we assume that $\omega' \not \in |\Sigma_1'|$. We consider now the closed segment $[\omega_0',\omega']$ in $N_\rdop'$. Let $\omega_1'$ be the point of $[\omega_0',\omega'] \cap  |\Sigma_1'|$ which is closest to $\omega'$. Then $\omega_1'$ is contained in an open face $\tau'$ of $\Sigma_1'$. Let $\omega_2'$ be a vertex of $\overline{\tau'}$. Using $[\omega_1',\omega'] \cap |\Sigma_1'|=\{\omega_1'\}$ and moving $\omega'$ sufficiently close to $\omega_1'$, we may assume also that 
\begin{equation} \label{point intersection}
[\omega_2',\omega'] \cap |\Sigma_1'|=\{\omega_2'\}.
\end{equation}
Now we use the notation from Proposition \ref{special orbit closure} and we apply this result two times. The irreducible component $Y_{\omega_2'}$ corresponding to the vertex $\omega_2'$ is the 
toric variety over $\tilde{K}$ associated to the fan ${\rm LC}_{\omega_2'}(\Sigma_1')$. The closure $\overline{Z}$ of the orbit $Z$ of $(\Ycal_{\Sigma})_s$ corresponding to $\tau$ is the 
toric variety over $\tilde{K}$ associated to the fan in $N(\overline{\tau})_\rdop$ which is given by the projections of ${\rm LC}_{\tau}(\Sigma_1)= \{ {\rm LC}_\tau(\nu) \mid \nu \in \Sigma_1, \, \nu \supset \tau \}$ to $N(\overline{\tau})_\rdop$.  We have an equivariant morphism $\varphi_{\omega_2'}: Y_{\omega_2'} \rightarrow \overline{Z}$ induced by $\varphi_s$. We deduce from \eqref{point intersection} that $|{\rm LC}_{\omega_2'}(\Sigma_1')|$ is a proper subset of $f_\rdop^{-1}(|{\rm LC}_{\tau}(\Sigma_1)|)$.  By the criterion of properness for homomorphisms of toric varieties over a field (see \cite{Fu2}, \S 2.4), $\varphi_{\omega_2'}$ is not proper. This contradicts properness of $\varphi_s$. We conclude that (a) implies (b).

To prove the converse, we assume that (b) holds. We get $f^{-1}(|\Sigma_0|)=|\Sigma_0'|$ and hence $\varphi_\eta$ is proper again by the criterion in \cite{Fu2}, \S 2.4. 
By \ref{consequences of the orbit correspondence}, the irreducible 
components of  $(\Ycal_{\Sigma'})_s$ correspond to the vertices $\omega'$ of $\Sigma_1$. Moreover, let $Z$ be the orbit of $(\Ycal_{\Sigma})_s$ corresponding to the open face $\tau$ of $\Sigma_1$ containing $f(\omega')$. As above, we get an equivariant morphism $\varphi_{\omega'}:Y_{\omega'} \rightarrow \overline{Z}$ of toric varieties over $\overline{K}$. It follows from (b) that $f_\rdop^{-1}(|{\rm LC}_{\tau}(\Sigma_1)|)=|{\rm LC}_{\omega'}(\Sigma_1')|$. The criterion in \cite{Fu2}, \S 2.4, shows that $\varphi_\omega$ is proper. As this holds for any irreducible component $Y_{\omega'}$, we get properness of $\varphi_s$. Indeed, it follows from Proposition \ref{purity of special fibre} that $\varphi_s$ is separated and of finite type and so we may use 
 \cite{EGA II}, Corollaire 5.4.5, to deduce properness of $\varphi_s$. 

It remains to see that $\varphi(\Xcal')$ is closed for any closed subset $\Xcal'$ of $\Ycal_{\Sigma'}$. We may assume that $\Xcal'$ is irreducible. Since $\varphi_s$ is proper, we may also assume that $\Xcal'$ is the closure of a closed subvariety $X'$ of $(\Ycal_{\Sigma'})_\eta$. Using Proposition \ref{generic orbit closure}, we may reduce the claim to the case $X' \cap T' \neq \emptyset$. 
Since $\varphi_\eta$ is proper, the generic fibre of $\varphi(\Xcal')$ is a closed subvariety $X$ of $\Ycal_\Sigma$. It remains to show that any point $z$ in the special fibre of the closure $\Xcal$ of $\varphi(\Xcal')$ is contained in $\varphi(\Xcal')$. By Proposition \ref{tropical lifting lemma}, the reduction map $\pi: X^\circ \cap \Tan \rightarrow \Xcal_s$ is surjective.  We conclude that $z=\pi(x)$ for some $x \in X^\circ \cap \Tan$ and hence $\trop_v(x) \in |\Sigma_1|$ by the orbit correspondence in Proposition \ref{uniform orbit correspondence}.  There is $x' \in (X')^{\rm an}$ with $\varphi^{\rm an}(x')=x$. By Chevalley's theorem, $\varphi_\eta(X' \cap T')$ is a constructible dense subset of $X$ and hence contains an open dense subset of $X$. We conclude that we may choose $x$ and $x'$ in the above argument such that $x' \in (T')^{\rm an}$ additionally. Using $f(\trop_v(x'))=\trop_v(\varphi^{\rm an}(x'))=\trop_v(x) \in |\Sigma_1|$, our assumption (b) on the fans leads to $\trop_v(x') \in |\Sigma_1'|$. By Proposition \ref{integral points of toric variety}, we have $x' \in (X')^\circ$ and hence its reduction $z':=\pi(x')$ is well-defined in $(\Xcal')_s$. We get $\varphi(z')=\varphi_s\circ \pi(x')=\pi \circ \varphi^{\rm an}(x')=z$. This proves $z \in \varphi(\Xcal')$ and therefore the morphism $\varphi$ is closed. We conclude that (b) implies (a). 

Next we prove that  (a) $\Rightarrow$ (c). We will use  similar arguments as for (a) $\Rightarrow$ (c) in the proof of Proposition \ref{properness of toric schemes}. By the above, we may assume that the equivalent properties (a) and (b) hold. We note that $\varphi$ is a quasi-compact separated morphism (see Proposition \ref{purity of special fibre}) and so we may apply the valuative criterion of universal closedness (\cite{EGA II}, Th\'eor\`eme 7.3.8). 
Let $L^\circ$ be a valuation ring with fraction field $L$ and suppose that we have a commutative diagram
\begin{equation*} 
\begin{CD}
\Spec(L) @>{\psi}>>\Ycal_{\Sigma'} \\
@VV V    @VV{\varphi} V\\
\Spec(L^\circ) @>{h}>> \Ycal_\Sigma
\end{CD}
\end{equation*}
of morphisms.  The criterion says that the morphism $\varphi$ is universally closed, if there is a morphism $g: \Spec(L^\circ) \rightarrow \Ycal_{\Sigma'}$ over $h$ which factors through $\psi$. If $h$ maps the generic point of $\Spec(L^\circ)$ into the special fibre of $\Ycal_\Sigma$, then we may replace the second column in the diagram by the special fibre $\varphi_s$ and the existence of $g$ follows from the valuative criterion of properness for $\varphi_s$. If $h$ maps the closed point of $\Spec(L^\circ)$ into the generic fibre of $\Ycal_{\Sigma'}$, then we may replace the second column in the diagram by the generic fibre $\varphi_\eta$ and the existence of $g$ follows from the valuative criterion of properness for $\varphi_\eta$. So we may assume that $h$ maps the generic point of $\Spec(L^\circ)$ into the generic fibre and the special point of $\Spec(L^\circ)$ into the special fibre. Then $L^{\circ\circ} \cap K^\circ = K^{\circ \circ}$ and it follows that $L^\circ$ is the valuation ring of a valuation $w$ on $L$ extending $v$. The morphism $\psi$ corresponds  to an $L$-rational point $P$ of the generic fibre  $(\Ycal_{\Sigma'})_\eta=Y_{\Sigma_0'}$ and we have to show that $P$ is an $L^\circ$-point of $\Ycal_{\Sigma'}$. It follows from Proposition \ref{generic orbit closure} again that we may assume $P \in T(L)$. Using that $\varphi(P)$ is an $L^\circ$-integral point of $\Ycal_\Sigma$,  Proposition \ref{integral points of toric variety} implies that $f_{\Gamma_L}(\trop_w(P))=\trop_w(\varphi(P)) \in \Delta(\Gamma_L)$ for some $\Delta \in \Sigma_1$. By (b), we have $f^{-1}(|\Sigma_1|)=|\Sigma'_1|$. It follows from \ref{compatibility of affine maps} that $\trop_w(P) \in \Delta'$ for some $\Delta' \in \Sigma_1$. Again Proposition \ref{integral points of toric variety} shows that $P$ is an $L^\circ$-integral point of $\Ycal_{\Sigma'}$. This proves (a) $\Rightarrow$ (c).

Conversely, we prove (c) $\Rightarrow$ (a). If (c) holds, then the generic and the special fibre of $\varphi$ are universally closed. As both of them are separated and of finite type (use Proposition \ref{purity of special fibre}), they are proper and hence we get (a).  If the value group $\Gamma$ is divisible or discrete in $\rdop$, Proposition \ref{purity of special fibre} implies that $\varphi$ is separated and of finite type over $K^\circ$. If (c) holds, then  $\varphi$ is proper over $K^\circ$. \qed

\vspace{2mm}
We have the following generalization of the above proposition to the case of a closed subscheme:

\begin{prop} \label{properness criterion for restriction of toric morphism}
Using the assumptions and notations from \ref{proper morphism}, let $X'$ be a closed subscheme of $T'$ with closure $\Xcal'$ in $\Ycal_{\Sigma'}$ and let $\phi$ be the restriction of the equivariant morphism $\varphi:\Ycal_{\Sigma'} \rightarrow \Ycal_\Sigma$ to $\Xcal'$. Then the following conditions are equivalent:
\begin{itemize}
\item[(a)] $\Trop_W(X') \cap (f_\rdop \times \id_{\rdop_+})^{-1}(|\Sigma|)= \Trop_W(X') \cap |\Sigma'|$;
\item[(b)] $\phi:\Xcal' \rightarrow \Ycal_\Sigma$ is universally closed.
\end{itemize}
If the equivalent conditions (a) and (b) hold, then the generic fibre of $\phi$ is proper. If  $\Gamma$ is divisible or discrete in $\rdop$, then (a) and (b) are equivalent to $\phi$ proper.
\end{prop}

\proof 
We assume  (a) and we will show (b): Similarly as in \cite{BS}, Proposition 3.15, one can construct a $\Gamma$-admissible fan $\Sigma''$ which subdivides $(f \times \id_{\rdop_+})^{-1}(\Sigma)$ and such that $\Sigma''$ has a subfan $\Sigma'''$ which subdivides $\Sigma'$.  We get a commutative diagram
\begin{equation*} 
\begin{CD}
\Ycal_{\Sigma'''} @>>>\Ycal_{\Sigma''} \\
@VV{\varphi'''} V    @VV{\varphi''} V\\
\Ycal_{\Sigma'} @>{\varphi}>> \Ycal_\Sigma
\end{CD}
\end{equation*}
of canonical equivariant morphisms. Since  $\Sigma'''$ is a subfan of $\Sigma''$, it is obvious that $\Ycal_{\Sigma'''}$ is an open subset of $\Ycal_{\Sigma''}$. By Proposition \ref{proper morphism}, the morphisms $\varphi''$ and $\varphi'''$ are universally closed. By Tevelev's Lemma \ref{Tevelev's Lemma}, the assumption (a) yields that the closure $\Xcal''$ of $X'$ in $\Ycal_{\Sigma''}$ is contained in the open subset $\Ycal_{\Sigma'''}$. It follows that $\Xcal''$ is the closure of $X'$ in $\Ycal_{\Sigma'''}$. Using that $\varphi'''$ is closed and restricts to the identity on the dense orbit $T'$, we get $\varphi'''(\Xcal'')=\Xcal'$. Let $\phi'':\Xcal'' \rightarrow \Ycal_\Sigma$ (resp. $\phi''':\Xcal'' \rightarrow \Xcal'$) be the restriction of $\varphi''$ (resp. $\varphi'''$) to the closed subscheme $\Xcal''$. Since $\varphi''$ is universally closed, the same is true for $\phi''$. We have the factorization $\phi''=\phi \circ \phi'''$. Using that the surjective morphism $\varphi'''$ remains surjective after base change (\cite{EGA I}, Proposition 3.5.2(ii)) and that $\phi''$ is universally closed, we deduce universal closedness of $\phi$. This proves (a) $\Rightarrow$ (b). 

We consider now the property
\begin{itemize}
\item[(a')] $\Trop_v(X') \cap f_\rdop^{-1}(|\Sigma_1|)= \Trop_v(X') \cap |\Sigma_1'|$.
\end{itemize}
If we would have assumed (a') instead of (a) in the above argument, then we could still show that the special fibre of $\Xcal''$ is contained in  $\Ycal_{\Sigma'''}$ and the same proof would show that the special fibre of $\phi$ is universally closed. This will be important in the proof of Proposition \ref{properness of closure}.

We assume now that (b) holds  and we will prove first (a'). It is clear that $\Trop_v(X') \cap f_\rdop^{-1}(|\Sigma_1|)\supset \Trop_v(X') \cap |\Sigma_1'|$ and we have to prove the reverse inclusion. Let $\omega' \in \Trop_v(X) \cap f_\rdop^{-1}(|\Sigma_1|)$. By Remark \ref{rationality of seminorms}, there is a valued field $(L,w)$ extending $(K,v)$ and an $L$-rational point $P$ of $X'$ such that $\trop_w(P)=\omega'$. We have $\trop_w(\phi(P))=f(\trop_w(P))=f(\omega') \in |\Sigma_1|$ and hence $\phi(P)$ is an $L^\circ$-integral point of $\Ycal_\Sigma$ (see Proposition \ref{integral points of toric variety}). We note that $\phi$ is a quasi-compact separated morphism (see Proposition \ref{purity of special fibre}) and hence we may apply the valuative criterion of universal closedness (\cite{EGA II}, Th\'eor\`eme 7.3.8). It follows that $P$ is an $L^\circ$-integral point of $\Xcal'$. By Proposition \ref{integral points of toric variety}, this means $\omega'=\trop_w(P) \in |\Sigma'_1|$ proving (a').

Now we show (b) $\Rightarrow$ (a). We have seen that (b) implies (a'). If we apply this with the trivial valuation  instead of $v$, we get 
$$\Trop_0(X') \cap f_\rdop^{-1}(|\Sigma_0|)= \Trop_0(X') \cap |\Sigma_0'|.$$
This together with (a') implies that (a) holds. 

Finally, we assume that $\Gamma$ is divisible or discrete in $\rdop$. It follows from Proposition \ref{purity of special fibre} and \cite{EGA I}, Proposition 5.5.1, that $\phi$ is a separated morphism of finite type. Then (b) yields that $\phi$ is proper. \qed

\begin{prop} \label{properness of closure} \label{proper closure and divisible case}
Let $X$ be a  closed subscheme of $T$. For the closure $\Xcal$ of $X$ in $\Ycal_\Sigma$, the following  are equivalent:
\begin{itemize}
\item[(a)] $\Trop_W(X) \subset |\Sigma|$;
\item[(b)] $\Trop_v(X) \subset |\Sigma_1|$;
\item[(c)] $\Xcal$ is universally closed over $\kcirc$.
\end{itemize}
If  $X$ is geometrically connected, then (a)--(c) are equivalent to:
\begin{itemize}
\item[(d)] the special fibre of  $\Xcal$ is non-empty and proper over $\tilde{K}$.
\end{itemize}
If the value group $\Gamma$ is divisible or discrete in $\rdop$, then (a)--(c) are equivalent to:
\begin{itemize}
\item[(e)] $\Xcal$ is proper over $\kcirc$.
\end{itemize}
If (a) holds, then the generic fibre of $\Xcal$ is proper over $K$. 
\end{prop}

\proof It follows from Proposition \ref{properness criterion for restriction of toric morphism} that (a) and (c) are equivalent and that they imply properness of the generic fibre $\Xcal_\eta$ over $K$. If the value group $\Gamma$ is divisible or discrete in $\rdop$, the same result shows that (a) and (c) are also equivalent to (e). If we apply this with the trivial valuation, then we get the following result:

\vspace{2mm}
\noindent {\it Step 1: Suppose that the valuation $v$ is trivial. Then the closure $\Xcal$ of $X$ in the toric variety $Y_{\Sigma_0}$ associated to the rational fan $\Sigma_0$ in $N_\rdop$ is proper if and only if $\Trop_0(X) \subset |\Sigma_0|$.}
\vspace{2mm}

Obviously, (a) yields (b).  We assume now (b) and we will show (d). We have seen in the proof of  Proposition \ref{properness criterion for restriction of toric morphism} that (b) (which is a special case of (a') there) implies (c). In particular, the special fibre $\Xcal_s$ is non-empty and universally closed over $\ktilde$. Since $\Xcal_s$ is always separated and of finite type over $\ktilde$ (see Proposition \ref{purity of special fibre}), we conclude that $\Xcal_s$ is proper over $\ktilde$. This proves (b) $\Rightarrow$ (d) without assuming irreducibility of $X$ or completeness of $K$.

For the converse, we assume $X$ geometrically connected.  
By the first step, the converse holds  for the trivial valuation and so we may assume that $v$ is non-trivial complete valuation. 
Arguing by contradiction, we assume that (d) holds and that $\Trop_v(X)$ is not a subset of $|\Sigma_1|$. Since the special fibre of $\Xcal$ is non-empty, Tevelev's Lemma \ref{Tevelev's Lemma} yields that $\Trop_v(X)$ intersects $|\Sigma_1|$. Since $\Trop_v(X)$ is a connected finite union of $\Gamma$-rational polyhedra (see Theorem \ref{Bieri-Groves theorem} and Proposition \ref{connected}), there is $\omega \in \Trop_v(X) \cap |\Sigma_1|$ such that $\Omega \cap \Trop_v(X)$ is not contained in $|\Sigma_1|$ for every neighbourhood $\Omega$ of $\omega$. Moreover, we may assume $m \cdot \omega \in N_\Gamma$ for some non-zero $m \in \ndop$. Then there is a $\Gamma$-admissible subdivision $\Sigma'$ of $\Sigma$  such that $\omega$ is a vertex of $\Sigma'_1$. Let $\Xcal'$ be the closure of $X$ in $\Ycal_{\Sigma'}$. By Proposition \ref{proper morphism}, the canonical $\tdop$-equivariant morphism $\varphi:\Ycal_{\Sigma'} \rightarrow \Ycal_{\Sigma}$ is closed and the special fibre $\varphi_s$ is proper. This shows  $\varphi(\Xcal')=\Xcal$  
and that the special fibre of $\Xcal'$ is also non-empty and proper over $\tilde{K}$. To simplify the notation, we may assume that $\Xcal = \Xcal'$. 

By Remark \ref{reduction to trivial valuation}, $X_\omega := \Xcal \cap Z_\omega$ is a closed subscheme of the dense orbit $Z_\omega$ of the toric variety $Y_\omega$ over $\tilde{K}$ with $\Trop_0(X_\omega)$ equal to the local cone of $\Trop_v(X)$ at $\omega$. This means that $\Trop_0(X_\omega)$ is not contained in the fan ${\rm LC}_\omega(\Sigma_1)$ of the toric variety $Y_\omega$. By Step 1, we conclude that the closure $\Xcal_\omega$ of $X_\omega$ in $Y_\omega$ is not proper over $\tilde{K}$. On the other hand, $\Xcal_\omega$ is a closed subscheme of $\Xcal \cap Y_\omega$. Since the special fibre of $\Xcal$ is assumed to be proper over $\tilde{K}$,  this has to be true also for its closed subscheme $\Xcal_\omega$. This is a contradiction and hence (d) implies (b). 

Next, we assume   $X$  geometrically connected and that (d) holds. We will show that the generic fibre $\Xcal_\eta$ is proper over $K$.  
By base extension and Proposition \ref{tropical base change}, we may assume that $K$ is an algebraically closed complete field and hence the value group $\Gamma$ is divisible in $\rdop$. Then $\Xcal$ is flat, separated and of finite presentation over $\kcirc$ (see Proposition \ref{purity of special fibre}). Using that $\Xcal_s$ is non-empty, we deduce that $\Xcal$ is faithfully flat over $K$. We have seen in \ref{reduction for Xan} that the reduction map $\pi:\Xcal_\eta^\circ \rightarrow \Xcal_s$ is defined on the compact analytic subdomain $\Xcal_\eta^\circ$ of $\Xcal_\eta^{\rm an}$. We have already seen that (d) yields (b) and hence $\Xcal_\eta^\circ \cap \Tan = X^\circ  = \trop_v^{-1}(|\Sigma_1|) \cap \Xan=\Xan$ implies that $\Xcal_\eta^\circ = \Xcal_\eta^{\rm an}$. Since $\Xcal_\eta$ is connected and $K$ is complete,  $\Xcal_\eta^{\rm an}$ is connected (see \cite{Berk1}, Theorem 3.4.8 and Theorem 3.5.3). By anticontinuity of the reduction map $\pi$, we deduce that $\Xcal_s$ is connected. Since $\tilde{K}$ is algebraically closed as well (\cite{BGR}, Lemma 3.4.1/4), we conclude that the special fibre $\Xcal_s$ is geometrically connected. The same holds obviously for the generic fibre $\Xcal_\eta$. We conclude that all assumptions of \cite{EGA IV}, Corollaire 15.7.11, are satisfied and this result shows that $\Xcal_\eta$ is proper over $K$.  

Finally, we prove that (b) yields (a). By base change and Proposition \ref{tropical base change}, we may assume that $K$ is  algebraically closed and that $X$ is irreducible. 
We have seen that (b) yields (d) and the above shows that $\Xcal_\eta$ is proper over $K$. By Step 1, we get $\Trop_0(X) \subset |\Sigma_0|$.  Moreover, (b) yields $\Trop_{\ve v}(X) \subset |\Sigma_\ve|$ for every $\ve > 0$ and so we conclude $\Trop_W(X) \subset |\Sigma|$. 
\qed



\begin{cor} \label{cone property proved}
Let $X$ be a closed subscheme of $T$. Then $\Trop_W(X)$ is the closure of the cone in $N_\rdop \times \rdop_+$ which is generated by $\Trop_v(X) \times \{1\}$. 
\end{cor}

\proof By base change and Proposition \ref{tropical base change}, we may assume  $K$  algebraically closed and $X$ irreducible. 
Let $C$ be the closure of the cone in $N_\rdop \times \rdop_+$  generated by $\Trop_v(X) \times \{1\}$. It follows easily from the definitions that $C$ agrees with $\Trop_W(X)$ on $N_\rdop \times (0, \infty)$ and hence $C \subset \Trop_W(X)$ using that $\Trop_W(X)$ is closed in $N_\rdop \times \rdop_+$ (see Proposition \ref{tropical cone properties}). The Gr\"obner complex $\Ccal(A,a)$ from \ref{Groebner complex for X} is a complete $\Gamma$-rational polyhedral complex in $N_\rdop$. We have seen in Theorem \ref{support of Groebner} that $\Ccal(A,a)$ has a subcomplex $\Ccal$ with support $\Trop_v(X)$. It follows from \ref{pointed subdivision} that $\Ccal$ has a pointed $\Gamma$-rational subdivision $\Ccal'$. By Remark \ref{Burgos-Sombra}, $\Sigma:=c(\Ccal')$ is a $\Gamma$-admissible fan in $N_\rdop \times \rdop_+$ with support $C$. Since $\Sigma_1=\Trop_v(X)$, Proposition \ref{properness of closure} yields that $\Trop_W(X) \subset |\Sigma|=C$ and we get the claim. \qed

\section{Tropical compactifications}  \label{tropical compactifications}

We keep the notation from Section \ref{Compactifications in toric schemes}, where we have studied the closure $\Xcal$ of a closed subscheme $X$ of $T$ in the toric scheme $\Ycal_\Sigma$ over $K^\circ$. In this section, we study tropical compactifications $\Xcal$ related to certain fans $\Sigma$ supported  on the tropical cone $\Trop_W(X)$ introduced in Section \ref{TropCone}. This generalizes results of Tevelev who handled the case of an integral $X$ over an algebraically closed field with trivial valuation (see \cite{Tev}) and of Qu  who obtained some results in the case of a discrete valuation (see \cite{Qu}). Their definition of a tropical fan seems simpler, but our definition is better suited to handle the case of a non-reduced $X$ and the definitions agree in the case of reduced closed subschemes.

\vspace{2mm}
Let $(K,v)$ be an arbitrary valued field and let $X$ be  any closed  subscheme  of $T$.

\begin{Def} \rm \label{tropical fan}
A {\it $\Gamma$-admissible tropical fan} for $X$ is a $\Gamma$-admissible fan $\Sigma$ in $N_\rdop \times \rdop_+$ such that $\Trop_v(X) \subset |\Sigma_1|$ and such that there is a closed subscheme $\Fcal$ of $\tdop \times_{K^\circ} \Ycal_\Sigma$ with the following properties:
\begin{itemize}
\item[(a)] The second projection induces a faithfully flat map $f:\Fcal \rightarrow \Ycal_\Sigma$. 
\item[(b)] The map $\Phi: \tdop \times_{K^\circ} \Ycal_\Sigma \rightarrow  \tdop \times_{K^\circ} \Ycal_\Sigma, (t,y) \mapsto (t^{-1}, t \cdot y)$ maps  $T\times_K X$ isomorphically onto  $f^{-1}(T)$.         
\end{itemize}
In this case, we call the closure $\Xcal$ of $X$ in $\Ycal_\Sigma$ a {\it tropical compactification} of $X$.
\end{Def}

\begin{rem}\rm  \label{fibre over T}
Let us consider just multiplication $m: T \times_K X \rightarrow T$. Then this is isomorphic to the trivial fibre bundle $X \times_K T$ over $T$. The isomorphism is given by $(t,x) \mapsto (x, t \cdot x)$. In particular, we see that $m$ is faithfully flat.  A tropical fan asks for extension of faithful flatness for $p_2$ from $\Phi(T \times_K X)$ to a closed subscheme $\Fcal$ of $\tdop \times_{K^\circ} \Ycal_\Sigma$. 

If $\Sigma$ is a $\Gamma$-admissible tropical fan, then it follows from flatness that the open subset $f^{-1}(T)$ of $\Fcal$ is dense. Using the isomorphism $\Phi$, we get $(\Phi^{-1}(\Fcal))_{\rm red}= \tdop \times_{K^\circ} \Xcal_{\rm 
red}$ since the right hand side is reduced by \cite{EGA IV}, Proposition 17.5.7. We conclude that the multiplication map $m:\tdop \times_{K^\circ} \Xcal_{\rm red} \rightarrow \Ycal_\Sigma$ is surjective and has the same topological properties as $f$.

If $X$ is reduced, then the closure $\Xcal$ is reduced. If we assume additionally that $\Gamma$ is divisible or discrete in $\rdop$,  then Proposition \ref{proper closure and divisible case} implies that $\Sigma$ is a  $\Gamma$-admissible tropical fan if and only if $\Xcal$ is proper and the multiplication map induces a faithfully flat map $\tdop \times_{K^\circ} \Xcal \rightarrow \Ycal_\Sigma$. Hence our definition is the same as Tevelev's definition of a tropical fan for varieties over a trivially valued algebraically closed field.

\end{rem}

We can now generalize Tevelev's  result to our framework:

\begin{thm} \label{existence of tropical fans}
Let $\Sigma(A,a)$ be the Gr\"obner fan for $X$ in $N_\rdop \times \rdop_+$ and let $\Sigma_X$ be the subcomplex with support $\Trop_W(X)$ as in Corollary \ref{support of Groebner fan}. Then every $\Gamma$-ad\-mis\-si\-ble fan $\Sigma$ which subdivides $\Sigma_X$ is a $\Gamma$-admissible tropical fan for $X$. In particular,  $\Gamma$-ad\-mis\-si\-ble tropical fans exist for every closed subscheme $X$ of $T$.
\end{thm}

\proof We keep the notation introduced in  Section \ref{Section: The Groebner complex} about Gr\"obner complexes. Let $\Ycal_{A,a}$ be the orbit closure of $\yb :=   [\overline{X}]$ in $\Hilb(\pdop^m_{K^\circ})$. 
Since the $\Gamma$-admissible fan $\Sigma$ subdivides the  subcomplex $\Sigma_X$ of $\Sigma(A,a)$,  the canonical morphism $T \rightarrow T \yb$ between the dense orbits extends to a  $\tdop$-equivariant morphism $\varphi:\Ycal_\Sigma \rightarrow \Ycal_{A,a}$ (see \ref{toric morphism}). Indeed, $\varphi$ is given on $T$ by $t \mapsto (\chi^u(t)y_u)_{u \in A} \in \pdop_K^N$ (see \ref{A-example}) and it is easy to see that this extends to the desired morphism.

We consider the closed subscheme $\Gcal := (\id \times \varphi)^{-1}(\Univ(\pdop^m_{K^\circ}))$ of $\pdop^m_{K^\circ} \times_{K^\circ} \Ycal_\Sigma$ which is flat over $\Ycal_\Sigma$. The fibre $\Gcal_y$ over $y \in T$ is equal to $y^{-1}\overline{X} \subset \pdop^m_K$. This makes it easy to check that
\begin{equation} \label{trivial family}
h: \Gcal|_{ T} \rightarrow \overline{X} \times_K T, \quad (z,y) \mapsto (y \cdot z, y)
\end{equation}
is an isomorphism over $T$.  Let $\Fcal$ be the restriction of $\Gcal$ to $\tdop \times_{K^\circ} \Ycal_\Sigma$. Then the second projection restricts to a flat morphism $f:\Fcal \rightarrow \Ycal_\Sigma$. Moreover,  axiom (b) from the definition of a tropical fan 
follows  from \eqref{trivial family}.  
By Corollary  \ref{support of Groebner fan}, $\Trop_W(X)$ is the support of $\Sigma_X$ and hence also from its subdivision $\Sigma$. 
By Tevelev's Lemma \ref{Tevelev's Lemma}, we conclude that every orbit intersects $\Xcal$ and hence the multiplication map $m:\tdop \times_{K^\circ} \Xcal_{\rm red} \rightarrow \Ycal_\Sigma$ is surjective. By the same argument as in Remark \ref{fibre over T}, we conclude that $f$ is surjective and hence faithfully flat. 
This means that $\Sigma$ is a tropical fan for $X$. 
Finally, we have seen in \ref{pointed subdivision} that a $\Gamma$-admissible fan exists which subdivides $\Sigma_X$. \qed

\begin{prop} \label{subdivision of tropical fan}
Let $\Sigma$ be a {\it $\Gamma$-admissible tropical fan} for $X$ and let $\Sigma'$ be a  $\Gamma$-rational fan which subdivides $\Sigma$. Then $\Sigma'$ is  a $\Gamma$-admissible tropical fan for $X$. 
\end{prop}

\proof Since $\Sigma'$ is a subdivision of $\Sigma$, we have $\Trop_v(X) \subset |\Sigma_1| = |\Sigma_1'|$. By \ref{toric morphism}, we get  a canonical $\tdop$-equivariant morphism $\varphi:\Ycal_{\Sigma'}\rightarrow \Ycal_\Sigma$ which is the identity on the dense open orbits $T$. 
Let us define the closed subscheme $\Gcal$ of $\tdop \times_{K^\circ} \Ycal_{\Sigma'}$ by the following Cartesian diagram:
\begin{equation*} 
\begin{CD}
\Gcal @>{f'}>> \Ycal_{\Sigma'}\\
@VV{\varphi'}V    @VV{\varphi}V\\
\Fcal @>{f}>> \Ycal_\Sigma
\end{CD}
\end{equation*}
Since $f'$ is obtained from $f$ by base change, we conclude that $f'$ is faithfully flat. Since we have $(\varphi')^{-1}(f^{-1}(T))=(f')^{-1}(T)$ and $\varphi$ is the identity on $T$, we deduce easily axiom (b) from the definition of a tropical fan. This proves the claim. \qed


\begin{prop} \label{support of tropical fan}
Every $\Gamma$-admissible tropical fan for $X$ in  $N_\rdop \times \rdop_+$  has support equal to $\Trop_W(X)$. 
\end{prop}

\proof It follows from Proposition \ref{properness of closure} that the support of a $\Gamma$-admissible tropical fan $\Sigma$ contains $\Trop_W(X)$. We have to show that every $\sigma \in \Sigma$ is contained in $\Trop_W(X)$.  We argue by contradiction and so we assume that $\sigma$ is not contained in $\Sigma$. Passing to a subdivision and using Proposition \ref{subdivision of tropical fan}, we may assume that $\sigma$ is disjoint from $\Trop_W(X)$.   It follows from Lemma \ref{Tevelev's Lemma} that the tropical compactification $\Xcal$ is disjoint from the orbit corresponding to $\relint(\sigma)$.  We conclude that the multiplication map $m:\tdop \times_{K^\circ} \Xcal_{\rm red} \rightarrow \Ycal_\Sigma$ is not surjective.  This contradicts Remark \ref{fibre over T}. \qed


\begin{prop} \label{proper intersection}
Let $X$ be a pure dimensional closed subscheme of $T$ and let $\Sigma$ be a $\Gamma$-admissible tropical fan for $X$ with tropical compactification $\Xcal$ of $X$. If $Z$ is any torus orbit in the generic (resp. special) fibre of $\Ycal_\Sigma$, then $Z \cap \Xcal$ is a non-empty pure dimensional scheme over $K$ (resp. $\tilde{K}$) with $$\dim(Z \cap \Xcal)=\dim(X)+\dim(Z)-n.$$ In particular, $Z$ intersects the generic (resp. special) fibre of $\Xcal$ properly.   
\end{prop}

\proof By Proposition \ref{purity of special fibre}, the special fibre of $\Xcal$ is also pure dimensional of the same dimension as $X$. By flatness of $f$ and Remark \ref{fibre over T}, the multiplication map $m:\tdop \times_{K^\circ} \Xcal \rightarrow \Ycal_\Sigma$ has pure dimensional fibres of constant fibre dimension. By Remark \ref{fibre over T}, this fibre dimension is equal to $\dim(X)$. For a closed point $z \in Z$, the fibre $m^{-1}(z)=\{(t,x) \in \tdop \times_{K^\circ} \Xcal \mid t \cdot x = z\}$ 
projects onto $ \Xcal \cap Z$. The fibres of this projection are isomorphic to $\Stab(z)$ and hence they have dimension $n-\dim(Z)$. By the fibre dimension theorem, we get
$$\dim(X)=\dim(m^{-1}(z))=\dim(Z \cap \Xcal) + n - \dim(Z)$$
proving the claim. The fibre dimension theorem yields also that $Z \cap \Xcal$ is pure dimensional. \qed




\vspace{2mm}

For a tropical fan and   $z \in (\Ycal_\Sigma)_s$, Remark \ref{fibre over T} shows $(m^{-1}(z))_{\rm red}\cong f^{-1}(z)_{\rm red}$; this is closely related to a certain initial degeneration, as we will see in the next remark.

\begin{rem}  \rm \label{multiplication fibre and initial degeneration}
Let $\Sigma$ be a $\Gamma$-admissible tropical fan in $N_\rdop \times \rdop_+$ and let $z \in (\Ycal_\Sigma)_s$. 
By the tropical lifting lemma (Proposition \ref{tropical lifting lemma}), there is 
$y \in \Tan \cap  Y_{\Sigma_0}^\circ$ with reduction $\pi(y)=z$. There is a valued field $(L,w)$ extending $(K,v)$ such that $y$ is an $L$-rational point of $T$ in the sense of Remark \ref{rationality of seminorms}. 
Let $f:\Fcal \rightarrow \Ycal_\Sigma$ be the faithfully flat family from Definition \ref{tropical fan}. We claim that $f^{-1}(z)_{\tilde{L}}$ is isomorphic to the special fibre of the closure of $y^{-1}X$ in $\tdop_{L^\circ}$. 
This means that $\In_{\trop_w(y)}(X)$ is  represented by the embedding $f^{-1}(z)\rightarrow \tdop_{\tilde L}$ given by the first projection.

To prove the claim, we may assume $K=L$. By flatness of $f$, the closure of $f^{-1}(y)$ is equal to $f^{-1}(\overline{y})$ (see Corollary \ref{affine closure and flat pull-back} and Remark \ref{closure in general}). We restrict the flat family $f$ to the closure $\overline{y}$ of $y$ in $\Ycal_\Sigma$. The generic fibre of this restriction is $f^{-1}(y)$ which is isomorphic to $y^{-1}X$ using the first projection of $T \times_K Y_{\Sigma_0}$ and axiom (b) in \ref{tropical fan}. 
Note that the first projection also gives a closed embedding of $f^{-1}(\overline{y})$ into $\tdop$ and hence the special fibre  $f^{-1}(z)$ is isomorphic to the special fibre of the closure of $y^{-1}X$ in $\tdop$ as claimed.
\end{rem}  

\begin{prop} \label{isomorphic flat fibres}
Let $\Sigma$ be a $\Gamma$-admissible tropical fan in $N_\rdop \times \rdop_+$ and let $z$ be an $F$-rational point of $\Ycal_\Sigma$ for a field $F$. If we use the first projection to identify the following fibres of $f$ with closed subschemes of $\tdop_F$  as in Remark \ref{multiplication fibre and initial degeneration}, then we have $f^{-1}(sz)=s^{-1} \cdot f^{-1}(z)$ for all $s \in \tdop(F)$.
\end{prop}

\proof Note that $F$ is either an extension of $K$ or of the residue field $\tilde{K}$. It is enough to consider the case $F/\tilde{K}$ as we may deduce the case  $F/K$ from the previous one by using the trivial valuation on $K$. There is a valued field $(L,w)$ extending $(K,v)$ such that the residue field $\tilde{L}$ contains $F$. Since all the objects are defined over $F$, it is enough to show $f^{-1}(sz)=s^{-1} \cdot f^{-1}(z)$ over $\tilde{L}$. Let $t\in\tdop(L^\circ)$ be a lift of $s$, i.e. $\pi(t)=s$. By enlarging $L$, we may assume that $z = \pi(y)$ for some $y \in T(L)$ (see Proposition \ref{tropical lifting lemma}). Using Remark \ref{multiplication fibre and initial degeneration}, we see that the special fibre of the closure of $(ty)^{-1}(X)$ is equal to $f^{-1}(sz)$. On the other hand, multiplication with $t^{-1}$ induces an automorphism of $\tdop$ and hence is compatible with taking closures. This automorphism is given on the special fibre by multiplication with $s^{-1}$ and hence we get the claim. \qed

\begin{cor} \label{tropical fan and initial degeneration}
Let $\Sigma$ be a $\Gamma$-admissible tropical fan in $N_\rdop \times \rdop_+$. Suppose that $\omega, \nu \in \relint(\sigma_1)$ for some $\sigma \in \Sigma$.  Then we have $\In_\omega(X)=\In_{\omega'}(X)$.
\end{cor}

\proof This follows immediately from the orbit correspondence (Proposition \ref{uniform orbit correspondence}), Remark \ref{multiplication fibre and initial degeneration} and Proposition \ref{isomorphic flat fibres}. \qed

\begin{prop} \label{small perturbation}
Let $\sigma$ be a cone of the  $\Gamma$-admissible tropical fan $\Sigma$ in $N_\rdop \times \rdop_+$. For every $\omega \in \sigma_1$ and every $\omega'= \omega +\Delta \omega \in \relint(\sigma_1)$, we have $$\In_{\omega'}(X)=\In_{\Delta \omega}(\In_\omega(X)),$$ where the initial degeneration at $\Delta \omega$ is with respect to the trivial valuation. 
\end{prop}

\proof We have seen in Corollary \ref{toric perturbation} that the identity holds in a neighbourhood of $\omega$. Now the claim follows from Corollary \ref{tropical fan and initial degeneration}. \qed

\section{Tropical multiplicities}  \label{tropical mult}

In this section, $X$ is a closed subscheme of $T$. We will define a tropical multiplicity function on the tropical variety $\Trop_v(X)$. It will be used to define $\Trop_v(X)$ as a tropical cycle, i.e. a weighted polyhedral complex satisfying the balancing condition. This appeared first in Speyer's thesis \cite{Spe}. The balancing condition relies on the description of the Chow cohomology of a toric variety given by Fulton--Sturmfels \cite{FS}. This is very implicit in the presentation here as we reduce the claim to the case of the trivial valuation where the balancing condition of $\Trop_v(X)$ is a result of Sturmfels and Tevelev based on \cite{FS}. Further references: \cite{AR}, \cite{BPR}, \cite{ST}.

\begin{Def}  \rm \label{regular points of tropical variety}
A point $\omega$ of $\Trop_v(X)$  is called {\it regular} if there is a polytope $\sigma \subset \Trop_v(X)$ such that $\relint(\sigma)$  is a neighbourhood of $\omega$ in $\Trop_v(X)$.
\end{Def}

\begin{prop} \label{local regularity}
A point $\omega$ of $\Trop_v(X)$ is regular if and only if $0$ is regular in $\Trop_0(\In_\omega(X))$. 
\end{prop}

\proof This follows immediately from Proposition \ref{trop trop}. \qed

\begin{art} \rm \label{setup for tropical multiplicity}
For $\omega \in N_\rdop$, we have seen that $\In_\omega(X)$ is a closed subscheme of the special fibre of $\tdop$ defined over a field extension of the residue field and it is well-defined up to multiplication with elements $g \in \tdop$ which are rational over a possibly larger field extension. Let $F$ be an algebraically closed field extension over which $\In_\omega(X)$ is defined. Then the irreducible components of $\In_\omega(X)$ over $F$ are also irreducible components over every field extension of $F$ and hence the following definition makes sense.
\end{art}

\begin{Def} \rm \label{tropical multiplicity}
The {\it tropical multiplicity} $m(\omega,X)$ of $\omega \in N_\rdop$ is defined as the sum of the multiplicities of the irreducible components of $\In_\omega(X)$ in  $\In_\omega(X)$ over the algebraically closed field $F$. For a cycle $Z = \sum m_Y Y$ of $T$ with prime components $Y$, we define the {\it tropical multiplicity of $Z$ in $\omega$} by $m(\omega,Z):=\sum_Y m_Y m(\omega,Y)$. 
\end{Def}

We have defined the initial degeneration as an equivalence class of closed subschemes up to multiplication by torus elements over an extension of the residue field (see \ref{initial degeneration at omega}). In the next result, we form the cycle of an initial degeneration. This means that we consider cycles up to the obvious linear extension of the above equivalence relation from prime components to all cycles. Note that the following result is a special case of \cite{OP}, Theorem 4.4.5. Here, we give a different proof using intersection theory with Cartier divisors. We have to deal with the fact that the models are usually non-noetherian and hence we cannot use  algebraic intersection theory, but there is an analytic replacement introduced in \cite{Gu1}. 

\begin{lem} \label{initial degeneration and cycle}
Let $\cyc(X)=\sum_Y m_Y Y$ be the representation of the cycle associated to $X$ as a sum of its irreducible components $Y$ counted with multiplicities. Then we have
$$\cyc(\In_\omega(X))= \sum_Y m_Y  \cyc(\In_\omega(Y)).$$ 
\end{lem}

\proof By base change, we may assume that $v$ is non-trivial and that $K$ is an algebraically closed complete field such that all occurring initial degenerations are  defined over $\tilde{K}$. Moreover, we may suppose that  $\omega = \trop_v(t)$ for some $t \in T(K)$. Replacing $X$ by $t^{-1}X$, we may assume $t=e$ and $\omega = 0$. Then $\In_\omega(X)$ is the special fibre $\Xcal_s$ of the closure $\Xcal$ of $X$ in $\tdop$.

By Proposition \ref{finite generation}, we have $\Xcal = \Spec(A)$ for a flat $K^\circ$-algebra of finite type. Let us choose a non-zero $\nu \in K^{\circ \circ}$. We have seen in \ref{comparision of reduction} that the $\nu$-adic completion $\hat{A}$ of $A$ is a flat $K^\circ$-algebra which is topologically of finite type, i.e. $\hat{\Xcal}:=\Spf(\hat{A})$ is an admissible formal affine scheme over $K^\circ$ in the theory of Raynaud, Bosch and L\"utkebohmert (see \cite{BL3}, \S 1). Its generic fibre is defined as the Berkovich spectrum $\Mcal(\Acal)$ of the strictly affinoid algebra $\Acal := \hat{A} \otimes_{K^\circ} K$ and it is equal to the affinoid subdomain $X^\circ$ of $\Xan$ from \ref{reduction for Xan}. If $T^\circ$ is the formal affinoid torus constructed in the same way from $T$, then we have $X^\circ=\Xan \cap T^\circ$. 

Using that $\Acal$ is a noetherian algebra, we have a theory of cycles and Cartier divisors on $X^\circ$ (see \cite{Gu1}, \S 2, for details). Hence we have a cycle decomposition $\cyc(X^\circ) = \sum_{W \in S} m_W W$ for a finite set $S$ of prime cycles of $X^\circ$. If $Y$ is an irreducible component of $X$, then the GAGA principle shows that $Y^\circ$ is a closed reduced analytic subvariety of $X^\circ$, but $Y^\circ$ is not necessarily irreducible. Hence we have
$\cyc(Y^\circ)=\sum_{W \in S_Y} W$
for a subset $S_Y$ of $S$. It is clear that $S$ is the disjoint union of the sets $S_Y$. By \cite{Gu1}, Proposition 6.3, we have 
$m_Y=m_W$ for all $W \in S_Y$. Moreover, it is obvious that $S_Y \neq \emptyset$ if and only if $Y$ meets the affinoid torus $T^\circ$ and the latter is equivalent to $\In_0(Y) \neq \emptyset$. By \cite{Gu2}, Lemma 4.5, we have
\begin{equation} \label{admissible cycle identity 1}
\cyc(\Xcal_s) = \sum_{W \in S}  m_W \cyc(\overline{W}_s)
\end{equation}
where $\overline{W}_s$ is the special fibre of the closure $\overline{W}$ of $W$ in $\hat{\Xcal}$. Similarly, we get
\begin{equation} \label{admissible cycle identity 2}
\cyc(\overline{Y}_s) = \sum_{W \in S_Y}  \cyc(\overline{W}_s)
\end{equation}
where $\overline{Y}$ is the closure of $Y$ in $\Xcal$. Using \eqref{admissible cycle identity 1}, \eqref{admissible cycle identity 2} and the above facts, we get 
$$\cyc(\In_0(X))=  \sum_Y \sum_{W \in S_Y} m_Y \cyc(\overline{W}_s)= \sum_Y m_Y \cyc(\overline{Y}_s) =       \sum_Y m_Y  \cyc(\In_0(Y))$$
proving the claim. \qed


\begin{prop} \label{properties of tropical multiplicities}
Tropical multiplicities have the following properties:
\begin{itemize}
\item[(a)] They are invariant under base change of $X$ or $Z$ to valued field $(L,w)$ extending $(K,v)$.
\item[(b)] The tropical multiplicity $m(\omega,Z)$ is linear in the cycle $Z$.
\item[(c)] For the cycle $\cyc(X)$ associated to $X$, we have $m(\omega,X)=m(\omega,\cyc(X))$.
\end{itemize}
\end{prop}

\proof Property (a) follows from Proposition \ref{initial degeneration and base change} and (b) is obvious. Finally, (c) follows from Lemma \ref{initial degeneration and cycle}. \qed

\vspace{2mm}
The following result shows that we may compute tropical multiplicities locally over the trivially valued residue field.

\begin{prop} \label{tropical multiplicities are local}
For  $\omega_0 \in N_\rdop$, there is a neighbourhood $\Omega$ of $\omega_0$ in $N_\rdop$ such that $m(\omega, X)=m(\omega - \omega_0, \In_{\omega_0}(X))$ for all $\omega \in \Omega$.
\end{prop}

\proof This follows from Corollary \ref{toric perturbation}. \qed

\vspace{2mm}
We have now the setup to generalize the following result of Sturmfels--Tevelev, which was given in the case of trivial valuations. 

\begin{thm} \label{tropical multiplicity is locally constant}
The restriction of the tropical multiplicity function $m(\cdot, X)$ to the open subset of regular points in $\Trop_v(X)$ is locally constant.
\end{thm}

\proof By Proposition \ref{tropical multiplicities are local}, we reduce to the case of a trivially valued base field. By Proposition \ref{properties of tropical multiplicities}, we may assume that base field is algebraically closed and that $X$ is an irreducible subvariety.  Then the claim follows from \cite{ST}, Corollaries 3.8 and 3.15. \qed

\begin{art} \rm \label{tropical cycles}
A $\Gamma$-rational polyhedral complex $\Ccal$ in $N_\rdop$ is called {\it of pure dimension $d$} if every maximal $\sigma \in \Ccal$ has dimension $d$. Such a complex is called {\it weighted} if it is endowed with a {\it multiplicity function} $m$ which maps every $d$-dimensional $\sigma \in \Ccal$ to a number $m_\sigma \in \zdop$. 

A polyhedron $\sigma \in \Ccal$ generates an affine space in $N_\rdop$ which is a translate of a linear space $\ldop_\sigma$. By $\Gamma$-rationality of $\sigma$, the vector space $\ldop_\sigma$ is defined over $\qdop$ and $N_\sigma := \ldop_\sigma \cap N$ is a lattice in $\ldop_\sigma$. 

We say that a weighted $\Gamma$-rational polyhedral complex $\Ccal$  in $N_\rdop$  of pure dimension $d$ satisfies the { \it balancing condition} if for every $d-1$-dimensional polyhedron $\nu$, we have
$$\sum_{\sigma \supset \nu} m_\sigma n_{\sigma,\nu} \in N_\nu,$$
where $\sigma$ ranges over all $d$-dimensional polyhedra of $\Ccal$ containing $\nu$, and $n_{\sigma,\nu}$ is any representative  of the generator of the $1$--dimensional lattice $N_\sigma/N_\nu$ pointing in the direction of $\sigma$. 

A weighted $\Gamma$-rational polyhedral complex $\Ccal$  in $N_\rdop$  of pure dimension $d$ is called a {\it tropical cycle} if it satisfies the balancing condition. We identify tropical cycles if there is a common $\Gamma$-rational subdivision of both complexes for which the multiplicities coincide. This allows us to add tropical cycles. In general, we define a $\Gamma$-rational tropical cycle $\Ccal$ in $N_\rdop$ as a formal sum $\Ccal=\sum_{j=0}^n \Ccal_j$, where $\Ccal_j$ is a tropical cycle   in $N_\rdop$  of pure dimension $j$. For details about tropical cycles, we refer to \cite{AR}. 
\end{art}

\begin{art} \rm \label{setup for tropical balancing theorem}
We suppose that $X$ is a pure-dimensional closed subscheme of $T$ and we set $d:=\dim(X)$.  Let $\Ccal$ be any $\Gamma$-rational polyhedral complex with support equal to $\Trop_v(X)$. By Theorem \ref{support of Groebner}, we know that such complexes exist. Theorem \ref{Bieri-Groves theorem} shows that $\Ccal$ is of pure dimension $d$. Note that the relative interior of a $d$-dimensional polyhedron $\sigma \in \Ccal$ is contained in the regular part of $\Trop_v(X)$. By Theorem \ref{tropical multiplicity is locally constant}, the multiplicity function $m(\cdot, X)$ is constant on $\relint(\sigma)$ and this constant is denoted by $m_\sigma$. We call $m_\sigma$ the {\it  tropical multiplicity} of $\sigma$.
\end{art}

\begin{thm} \label{balancing theorem}
Under the hypothesis of \ref{setup for tropical balancing theorem},  the complex $\Ccal$ endowed with the tropical multiplicities is a $\Gamma$-rational tropical cycle of pure dimension $d$.
\end{thm}

\proof The balancing condition is a local condition in any $\omega \in  N_\rdop$. By Propositions \ref{trop trop} and  \ref{tropical multiplicities are local}, it is enough to check the balancing condition for $\In_\omega(X)$ in a neighbourhood of $0$. Hence we have reduced the claim to the case of trivial valuation. Again, we may assume that the base field is algebraically closed and that $X$ is an irreducible subvariety. This case is proved in \cite{ST}, Corollary 3.8. \qed

\begin{rem} \rm  \label{tropical cycle for closed subschemes}
It follows from Theorem \ref{tropical multiplicity is locally constant} that the tropical cycle from Theorem \ref{balancing theorem} does not depend on the choice of the complex $\Ccal$ from \ref{setup for tropical balancing theorem}. We conclude that $\Trop_v(X)$ is canonically a tropical cycle which we denote also by $\Trop_v(X)$. 

If $X$ is any closed subscheme of $T$, then we define $\Trop_v(X)$ by linearity in its irreducible components, i.e. we set $\Trop_v(X):=\sum_Y m_Y \Trop_v(Y)$ as a tropical cycle, where $m_Y$ is the multiplicity of $X$ in the irreducible component $Y$. This is a tropical cycle in $N_\rdop$ with support equal to the set-theoretic tropical variety of $X$. By Proposition \ref{properties of tropical multiplicities}, this agrees with the above construction in the pure dimensional case.

If $Z=\sum_Y m_Y Y$ is any cycle on $X$ with prime components $Y$, then we define the {\it tropical cycle associated to $Z$} by $\Trop_v(Z):= \sum_Y m_Y \Trop_v(Y)$, where we use the induced reduced structure on every $Y$.
\end{rem}



\begin{prop} \label{Tropical cycle and localization}
Let $X$ be any closed subscheme of $T$ and let $\omega \in N_\rdop$. Replacing the polyhedra in the tropical cycle $\Trop_v(X)$ by its local cones in $\omega$ and using the same tropical multiplicities, we get a tropical cycle in $N_\rdop$ which is equal to the tropical cycle $\Trop_0(\In_\omega(X))$ with respect to the trivial absolute value $0$.
\end{prop}

\proof This follows from Proposition \ref{tropical multiplicities are local} and Proposition \ref{trop trop}. \qed

\begin{art} \label{tropicalization of a homomorphism} \rm
Let $T'$ be another split torus over $K$ with lattice $N'$ of one--parameter--subgroups. Let $\varphi:T \rightarrow T'$ be a homomorphism of split tori over $K$. This induces a homomorphism $M' \to M$ of character lattices and hence we get a linear map $N \rightarrow N'$. The base change of this map to $\rdop$ is easily seen to be the unique map $\Trop_v(\varphi): N_\rdop \rightarrow N_\rdop'$ such that $\Trop_v(\varphi) \circ \trop_v = \trop_v \circ \varphi$. 
\end{art}

\begin{art} \label{push-forward of cycles} \rm 
The push-forward of a cycle $Z$ on $X$ with respect to the homomorphism $\varphi$ is a cycle $\varphi_*(Z)$ on the closure $X'$ of $\varphi(X)$ defined in the following way: If $Z$ is a prime cycle and $Z'$ is the closure of $\varphi(Z)$, then 
\begin{equation*}
\varphi_*(Z):= \begin{cases}
[K(Z):K(Z')] Z', & \text{if $[K(Z):K(Z')]< \infty$,}\\
0, &  \text{if $[K(Z):K(Z')]= \infty$}.
\end{cases}
\end{equation*}
In general, $\varphi_*(Z)$ is defined by linearity in its prime components. 

Usually, the push-forward of cycles is defined with respect to proper morphisms. This could be easily obtained by using tropical compactifications as in Section \ref{tropical compactifications}, but as we are not interested in compatibility with rational equivalence of cycles, this plays no role here.
\end{art}

\begin{art} \label{push-forward of tropical cycles} \rm
We will explain how the linear map $f:=\Trop_v(\varphi):N_\rdop \rightarrow N_\rdop'$ induces a push--forward map of tropical cycles. For details, we refer to \cite{AR}, \S 7. Let $\Ccal$ be a tropical cycle in $N_\rdop$ of pure dimension $d$. After a subdivision of $\Ccal$, we may assume that 
$$f_*(\Ccal):=\{f(\sigma) \mid \text{$\sigma$ is a face of $\nu \in \Sigma$ with $\dim(f(\nu))=d$}\}$$
is a ($d$-dimensional $\Gamma$-rational) polyhedral complex in $N_\rdop'$. We define the multiplicity of a $d$-dimensional $f(\sigma) \in f_*(\Ccal)$ by
$$m_{f(\sigma)}:=\sum_{\nu \subset f^{-1}(f(\sigma))}[N_\nu:N_{f(\sigma)}] m_\nu ,$$
where $\nu$ ranges over all $d$-dimensional $\nu \in \Ccal$ contained in $ f^{-1}(f(\sigma))$. Endowed with these multiplicities, we get a weighted polyhedral complex which is a tropical cycle in $N_\rdop'$. It might happen that $f_*(\Ccal)$ is empty, then we get the tropical zero-cycle. 
\end{art}

The following result is the Sturmfels--Tevelev multiplicity formula (see \cite{ST}). It was generalized to the case of non-trivial valuations in \cite{BPR}, Corollary 8.4 and Appendix A.

\begin{thm} \label{Sturmfels-Tevelev multiplicity formula}
Let $\varphi:T \rightarrow T'$ be a homomorphism of split tori over $K$ and let $Z$ be cycle on  $T$. Then we have 
$$\Trop_v(\varphi)_*(\Trop_v(Z))=\Trop_v(\varphi_*(Z))$$ 
as an identity of tropical cycles.
\end{thm}

\proof By base change, we may assume that $K$ is an algebraically closed field with a complete non-trivial valuation. Using linearity of the identity in the prime  components of $Z$, we may assume that $X=Z$ is an integral closed subscheme of $T$. If $\dim(X') < \dim(X)$ for $X':=\varphi(X)$, then $\varphi_*(X)=0$ by definition. Since $\Trop_v(X)$ is a polyhedral complex of pure dimension $d:=\dim(X)$ and since $\varphi_*(\Trop_v(X))$ is a tropical cycle supported in $\Trop_v(X')$ which is of lower dimension, we conclude $\varphi_*(\Trop_v(X))= 0$ as well. So we may assume that $\varphi$ induces a generically finite map $X \rightarrow X'$ and then we may deduce the claim from \cite{BPR}, Corollary 8.4. \qed

\section{Proper intersection with orbits}

As usual, $(K,v)$ is a valued field which serves as a ground field. 
Let $\Sigma$ be a $\Gamma$-admissible fan in $N_\rdop \times \rdop_+$ with associated toric scheme $\Ycal_\Sigma$ over $K^\circ$. Let $X$ be a closed subscheme  of the dense torus $T$ with closure $\Xcal$ in $\Ycal_\Sigma$. We have seen in Proposition \ref{proper intersection} that $\Xcal$ intersects the orbits of $\Ycal_\Sigma$ properly in case of a tropical fan. In this section, we will generalize this result and we prove that this property is a purely combinatorial property of the fan $\Sigma$. 
 I am very grateful to Sam Payne for explaining to me some of the arguments for this nice result. 

\begin{art} \label{notation from convex geometry} \rm
Let $\sigma \in \Sigma_1$ and $\tau = \relint(\sigma)$. Then $\sigma$ generates an affine space in $N_\rdop$ which is a translate of a linear space $\ldop_\sigma$. By $\Gamma$-rationality of $\sigma$, the vector space $\ldop_\sigma$ is defined over $\qdop$. Then  $N_\sigma := N \cap \ldop_\sigma$ and $N(\sigma):=N/N_\sigma$ are free abelian groups of finite rank with quotient homomorphism $\pi_\sigma:N \rightarrow N(\sigma)$. Dually, we have $M(\sigma):=\ldop_\sigma^\bot \cap M=\Hom(N(\sigma),\zdop)$. 

For $S \subset N_\rdop$, we define the {\it local cone of $S$ at $\tau$} by 
$${\rm LC}_\tau(S) := \bigcup_{\omega \in \tau}{\rm LC}_\omega(S)$$
using the local cones at points from \ref{local cone}. If $S$ is a polyhedron containing $\tau$, then we have ${\rm LC}_\tau(S)= {\rm LC}_\omega(S)$ for any $\omega \in \tau $. 
\end{art}

\begin{art} \label{torus for orbit} \rm
We recall from \ref{consequences of the orbit correspondence} that $\tau$ corresponds to an orbit $Z=Z_\tau$ of the special fibre of   $\Ycal_\Sigma$. 
By choosing a base point $z_0 \in Z(\ktilde)$, Proposition \ref{special orbit closure} shows that $Z$  may be identified with the torus $\Spec(\tilde{K}[M(\sigma)_\tau])$ for the sublattice $M(\sigma)_\tau := \{u \in M(\sigma) \mid \langle u, \omega \rangle \in \Gamma \; \forall \omega \in \tau \}$ of finite index in $M(\sigma)$. We get tropical varieties of closed subschemes of $Z$ with respect to the trivial valuation which do not depend on the choice of the base point $z_0$. This is used in the following result which generalizes Remark \ref{reduction to trivial valuation}.
\end{art}

\begin{prop} \label{tropical variety for orbit intersection}
Using the notions from above, we have $$\Trop_0(\Xcal \cap Z_\tau) = \pi_\sigma({\rm LC}_\tau(\Trop_v(X))).$$
\end{prop}

\proof By base change and Lemma \ref{base change and tropical compactification}, we may assume that $\Gamma$ is divisible and hence we have $M(\sigma)=M(\sigma)_\tau$. We assume first that $\Sigma$ has a tropical subfan for $X$ (i.e. a subfan $\Sigma'$ which is a tropical fan for $X$, see Definition \ref{tropical fan} and \ref{cone}). If $\tau \cap \Trop_v(X)$ is empty, then Tevelev's Lemma \ref{Tevelev's Lemma} shows that both sides of the claim are empty. By Proposition \ref{support of tropical fan}, we have $\Trop_v(X)=|\Sigma_1'|$ and so we may assume that $\tau \subset \Trop_v(X)$. We choose $\omega \in \tau  \cap N_\Gamma$. By translation, we may assume that $\omega = 0$ and therefore the affine toric scheme $\Ucal_\omega$ from \ref{polyhedral scheme} is just the split torus $\tdop$ over $K^\circ$. We conclude that $\In_\omega(X)$ is the special fibre of the closure of $X$ in $\Ucal_\omega$. 
To identify $Z_\tau$ with $T(\sigma):= \Spec(\tilde{K}[M(\sigma)])$,    we choose the base point $z_0$ of ${Z_\tau}$ as the reduction of the unit element in $T(K)$. Then the canonical quotient homomorphism $q: \tdop_{\tilde{K}} \rightarrow T(\sigma)$ of tori over $\tilde{K}$ maps $\In_\omega(X)$ into $\Xcal \cap {Z_\tau}$. Since $\tau$ is an open face of a tropical subfan of $\Sigma$, the proof of Proposition \ref{proper intersection} and Remark \ref{multiplication fibre and initial degeneration} show that $\In_\omega(X)= q^{-1}(\Xcal \cap {Z_\tau})$ holds set theoretically. This yields
\begin{equation} \label{projection and trop}
\pi_\sigma(\Trop_0(\In_\omega(X)))= \Trop_0(\Xcal \cap {Z_\tau}).
\end{equation}
By Proposition \ref{trop trop}, we have $\Trop_0(\In_\omega(X))= {\rm LC}_\omega(\Trop_v(X))$. Since $\Trop_v(X)$ is a finite union of polyhedra which either contain $\tau$ or are disjoint from $\tau$, we get ${\rm LC}_\omega(\Trop_v(X))={\rm LC}_\tau(\Trop_v(X))$. Inserting these facts into \eqref{projection and trop}, we get the proposition.

Now we prove the proposition in the case of an arbitrary $\Gamma$-admissible fan $\Sigma$. The affine toric scheme $\Ucal_\sigma$ associated to the closure $\sigma$ of $\tau$ is an open subset of $\Ycal_\Sigma$ containing $Z_\tau$. The claim in the proposition depends only on $\Ucal_\sigma$ and hence we may change $\Sigma$ outside of $\sigma$. So we may assume that $\Sigma_1$ is a complete $\Gamma$-rational polyhedral complex containing $\sigma$.  
By Theorem \ref{existence of tropical fans} and Proposition \ref{subdivision of tropical fan}, there is a $\Gamma$-admissible fan  $\Sigma'$ which is a subdivision of $\Sigma$ and  which  has a tropical subfan. Then we have a canonical $\tdop$-equivariant morphism $\varphi:\Ycal_{\Sigma'} \rightarrow \Ycal_\Sigma$ of $\tdop$-toric schemes over $K^\circ$. It follows from Proposition \ref{proper morphism} that $\varphi$ is closed and surjective. For the closure $\Xcal'$ of $X$ in $\Ycal_{\Sigma'}$, we conclude that $\varphi(\Xcal')=\Xcal$. Relevant for our purposes is that $\tau$ has a subdivision into open faces $\tau_1, \dots, \tau_r$ of $\Sigma'$. The orbit correspondence in Proposition \ref{uniform orbit correspondence} leads to the partition of $\varphi^{-1}({Z_\tau})$ into the orbits $Z_{\tau_1}, \dots,Z_{\tau_r}$.  We conclude that 
$\Xcal \cap Z_\tau$ is the union of the sets $\varphi(\Xcal' \cap Z_{\tau_i})$. Let $\sigma_i$ be the closure of $\tau_i$ and let $\pi_i:N(\sigma_i) \rightarrow N(\sigma)$ be the canonical homomorphism. Then we get
\begin{equation} \label{Covering of the tropical variety}
\Trop_0(\Xcal \cap {Z_\tau})=\bigcup_{i=1}^r \Trop_0(\varphi(\Xcal' \cap Z_{\tau_i}))=\bigcup_{i=1}^r \pi_i(\Trop_0(\Xcal' \cap Z_{\tau_i})).
\end{equation}
Using the special case above, we have 
\begin{equation*} \label{tropical special case}
\Trop_0(\Xcal' \cap Z_{\tau_i}) = \pi_{\sigma_i}({\rm LC}_{\tau_i}(\Trop_v(X))).
\end{equation*}
Inserting this in \eqref{Covering of the tropical variety} and using $\pi_i \circ \pi_{\sigma_i}= \pi_\sigma$, we get the claim. \qed

\vspace{2mm}

For simplicity, we assume now that the closed subscheme $X$ of $T$ is of pure dimension $d$. By the Bieri--Groves Theorem \ref{Bieri-Groves theorem}, there is a finite set $S$ of $d$-dimensional $\Gamma$-rational polyhedra in $N_\rdop$ such that $\Trop_v(X)= \bigcup_{\Delta \in S} \Delta$. 


\begin{cor} \label{dimension formula for orbit intersection}
Under the hypothesis above, we have $$\dim(\Xcal \cap Z_\tau)= d- \inf\{\dim(\Delta \cap \tau)\mid \Delta \in S, \, \Delta \cap \tau \neq \emptyset\}.$$
\end{cor}

\proof By Tevelev's Lemma \ref{Tevelev's Lemma}, $\Xcal \cap Z_\tau$ is empty if and only if no $\Delta \in S$ intersects $\tau$. We see that the claim holds in this special case as the dimension of the empty set is defined as $-\infty$ and the infimum over an empty set  is $\infty$. Proposition \ref{tropical variety for orbit intersection} shows that we have
\begin{equation} \label{union formula}
\Trop_0(\Xcal \cap Z_\tau)= \bigcup_{\Delta \in S} \pi_\sigma({\rm LC}_\tau(\Delta)).
\end{equation}
For $\Delta \in S$ with $\Delta \cap \tau \neq \emptyset$, we have $\dim(\pi_\sigma({\rm LC}_\tau(\Delta)))=d-\dim(\Delta \cap \tau)$. Using this in \eqref{union formula}, we get the claim. \qed


\begin{rem} \label{dimension formula for generic orbit intersection} \rm
We recall from \ref{consequences of the orbit correspondence} that the open faces $\tau$ of $\Sigma_0$ correspond to the orbits $Z_\tau$ contained in the generic fibre $Y_{\Sigma_0}$ of $\Ycal_\Sigma$. If we use a decomposition $\Trop_0(X)= \bigcup_{\Delta \in S} \Delta$ into $d$-dimensional rational cones $\Delta$ in $N_\rdop$ (see Remark \ref{tropical variety for trivial absolute value}), then Corollary \ref{dimension formula for orbit intersection} holds also for these orbits. This follows immediately from Corollary \ref{dimension formula for orbit intersection} replacing $v$ by the trivial valuation. Then the generic fibre is equal to the special fibre.  
\end{rem}

\begin{art} \rm \label{definition of proper intersection}
Let $X$ be a closed subscheme of $T$ of pure dimension $d$ with closure $\Xcal$ in $\Ycal_\Sigma$ and let $Z_\tau$ be the orbit of $\Ycal_\Sigma$ induced by  the open face $\tau$ of $\Sigma_1$ (resp. $\Sigma_0$). We say that $\Xcal$ {\it intersects $Z_\tau$ properly} if  $\dim(\Xcal \cap Z_\tau)= d - \dim(\tau)$. We emphasize that in this case, $\Xcal \cap Z_\tau$ is not empty. Note that $\Ycal_\Sigma$ is a noetherian topological space and one can easily show that $\Xcal$ intersects $Z_\tau$ properly if and only if every irreducible component of $\Xcal \cap Z_\tau$ has  codimension  in $\Xcal$ equal to $\codim(Z_\tau, \Ycal_\Sigma)$. 
\end{art}

The following result was shown to me by Sam Payne.

\begin{prop} \label{pure dimensionality of orbit intersection}
If $\Xcal$ intersects $Z_\tau$ properly, then $\Xcal \cap Z_\tau$ is pure dimensional.
\end{prop}

\proof We may assume that $\tau \in \Sigma_1$ and hence $Z_\tau$ is contained in the special fibre of $\Ycal_\Sigma$. Indeed, the case $\tau \in \Sigma_0$ follows as usual from this replacing $v$ by the trivial valuation. 
We choose a vertex $\omega$ of $\overline{\tau}$. By Proposition \ref{special orbit closure}, the irreducible component $Y_\omega$ of $\Xcal_s$ corresponding to $\omega$ is a toric variety over $\tilde{K}$ associated to the fan ${\rm LC}_\omega(\Sigma_1):=\{ {\rm LC}_\omega(\Delta) \mid \Delta \in \Sigma_1\}$. 
We note that $Z_\tau$ is also an orbit of $Y_\omega$ and we have $\codim(Z_\tau,Y_\omega)=\dim(\tau)$. For every $z \in Z_\tau$, there is a neighbourhood $U$ of $z$ in $Y_\omega$ such that $Z_\tau \cap U$ is set theoretically the intersection of $\codim(Z_\tau,Y_\omega)$ effective Cartier divisors. This follows from Hochter's theorem which says that a toric variety is Cohen-Macaulay (see \cite{CLS}, Theorem 9.2.9). 
We conclude that every irreducible component of $\Xcal \cap Z_\tau$ has dimension at least $d- \dim(\tau)$.
\qed

\begin{prop} \label{proper orbit intersection and tropical variety}
Let $\tau$ be an open face of $\Sigma_1$. If $\Xcal$ intersects $Z_\tau$ properly, then $\tau \subset \Trop_v(X)$. 
\end{prop}

\proof 
Assuming that $\Xcal$ intersects $Z_\tau$ properly, we deduce from Corollary \ref{dimension formula for orbit intersection} the following fact which is crucial for the proof: If $\Delta$ is any $d$-dimensional polyhedron contained in $\Trop_v(X)$ and if $\Delta \cap \tau \neq \emptyset$, then we have 
\begin{equation} \label{following fact}
\dim(\Delta \cap \tau)= \dim(\tau). 
\end{equation}
By assumption, $\Xcal \cap Z_\tau$ is non-empty and hence $\tau \cap \Trop_v(X) \neq \emptyset$ by Tevelev's Lemma \ref{Tevelev's Lemma}. By the Bieri--Groves theorem and Theorem \ref{support of Groebner}, there is a complete $\Gamma$-rational polyhedral complex $\Dcal$ in  $N_\rdop$ of pure dimension $d$ with a subcomplex $\Ccal$ such that $\Trop_v(X)=|\Ccal|$.

Using these facts, the proposition will follow from elementary arguments in convex geometry. Especially important is the collection $\Ecal$ of all $\sigma \in \Dcal$ with $\dim(\relint(\sigma) \cap \tau)=\dim(\tau)$. We note that  $(\sigma \cap \tau)_{\sigma \in {\Ecal}}$ is a covering of the open face $\tau$ which is like a tiling of $\tau$. We have seen above that $\tau \cap \Trop_v(X) \neq \emptyset$ and hence there is a $d$-dimensional polyhedron $\Delta \in \Ccal$ with $\Delta \cap \tau \neq \emptyset$. Using an appropriate closed face of $\Delta$, we get the existence of a polyhedron $\sigma \in \Ccal \cap \Ecal$. 

Let $\sigma' \in \Ecal$ such that $\sigma' \cap \tau$ is a direct neighbor of $\sigma \cap \tau$ which means  that $\nu:= \sigma \cap \sigma' \cap \tau$ has dimension equal to $\dim(\tau)-1$. Using the above tiling of $\tau$, the open face $\tau$ may be covered by using successively such  neighboring $\sigma' \cap \tau$. We conclude that it is enough to show that $\sigma' \cap \tau \subset \Trop_v(X)$. 

Remembering that $\Ccal$ is a polyhedral complex of pure dimension $d$, there is a  $d$-dimensional polyhedron $\Delta \in \Ccal$ with closed face $\sigma$. Note that $\nu$ is obtained by intersecting $\tau$ with the proper closed face $\sigma \cap \sigma'$ of $\sigma$. We conclude that 
there is a closed face $\rho$ of $\Delta$ with $\dim(\rho)=d-1$ which contains $\nu$ but not $\sigma$. We have $\dim(\relint(\sigma) \cap \tau)=\dim(\tau)$ and hence $\Delta \cap \tau = \sigma \cap \tau$ contains $\rho \cap \tau$ as a proper subset. Since $\nu$ is of codimension $1$ in $\sigma \cap \tau$, we get $\nu=\rho \cap \tau$. 

We choose a hyperplane $H$ in $N_\rdop$ which contains $\rho$ but not $\Delta$. By the balancing condition in  Theorem \ref{balancing theorem}, there is a $d$-dimensional polyhedron $\Delta' \in \Ccal$ with closed face $\rho$ on the other side of $H$ than $\Delta$. Since $\nu$ is the border of $\sigma \cap \tau$ and $\sigma' \cap \tau$ in the above tiling  of $\tau$, we conclude that $\Delta' \cap \tau \subset \sigma' \cap \tau$. Using \eqref{following fact} for $\Delta'$, we get $\dim(\Delta' \cap \tau)= \dim(\tau)$ and hence $\Delta' \cap \tau = \sigma' \cap \tau$. This proves $\sigma' \cap \tau \subset \Trop_v(X)$. \qed

\vspace{2mm}

To deal with orbits in the special fibre and in the generic fibre simultaneously, one has to use the $\Gamma$-admissible fan $\Sigma$ in $N_\rdop \times \rdop_+$ and the tropical cone $\Trop_W(X)$ of $X$ in $N_\rdop \times \rdop_+$ (see Definition \ref{tropical cone of X}).

\begin{thm} \label{proper and proper intersection with orbits}
Let $\Sigma$ be a $\Gamma$-admissible fan in $N_\rdop \times \rdop_+$ and let $X$ be a closed subscheme of $T$ of pure dimension $d$. Then the following properties are equivalent for the closure $\Xcal$ of $X$ in the toric scheme $\Ycal_\Sigma$:
\begin{itemize}
\item[(a)] The special fibre $\Xcal_s$ is non-empty, proper over $\tilde{K}$, and $\Xcal$ intersects all the orbits of $\Ycal_\Sigma$ properly.
\item[(b)] The support of $\Sigma$ is equal to the tropical cone $\Trop_W(X)$.
\end{itemize} 
If the value group $\Gamma$ is divisible or discrete in $\rdop$, then (a) and (b) are also equivalent to the condition that $\Xcal$ is a  proper scheme over $K^\circ$ which intersects all the orbits properly. \end{thm}

\proof We assume that (a) holds. By Proposition \ref{proper orbit intersection and tropical variety}, the assumption that $\Xcal$ intersects all orbits properly yields that  $|\Sigma_1|$  is contained in $\Trop_v(X)$. If we replace $v$ by the trivial valuation, then the same argument shows that $|\Sigma_0| \subset \Trop_0(X)$. Since $\Trop_W(X)$  is the closed cone in $N_\rdop \times \rdop_+$ generated by $\Trop_v(X) \times \{1\}$ (see Proposition \ref{cone property}),  we conclude that $|\Sigma| \subset \Trop_W(X)$.
On the other hand, $\Xcal_s$ is a non-empty proper scheme over $\tilde{K}$ and hence $\Trop_v(X)$ is contained in $|\Sigma_1|$ by Proposition \ref{properness of closure}. We conclude that (a) yields (b).

Now we suppose that (b) holds. Then $\Sigma_1$ is a $\Gamma$-rational complex with support equal to $\Trop_v(X)$. We choose an open face $\tau$ of $\Sigma_1$ with corresponding orbit $Z_\tau$ in the special fibre of $\Ycal_\Sigma$. For any $d$-dimensional polyhedron $\Delta \in \Sigma_1$, either $\Delta \cap \tau$ is empty or $\tau$. From Corollary \ref{dimension formula for orbit intersection}, we deduce that $\Xcal$ intersects $Z_\tau$ properly. Using $\Trop_0(X)= |\Sigma_0|$, Remark \ref{dimension formula for generic orbit intersection} shows that $\Xcal$ intersects the orbits in the generic fibre of $\Ycal_\Sigma$ properly. It follows from Proposition \ref{properness of closure} that $\Xcal_s$ is a non-empty proper scheme over $\tilde{K}$.

If $\Gamma$ is divisible or discrete in $\rdop$, then the last claim follows immediately from Proposition \ref{proper closure and divisible case}. \qed

\begin{rem} \label{compare with tropical fans} \rm
We have seen in Proposition \ref{support of tropical fan} that every tropical fan satisfies the equivalent properties (a) and (b) of Theorem \ref{proper and proper intersection with orbits}. However, the converse does not hold as it was shown by Sturmfels--Tevelev in Example 3.10 of \cite{ST} and by Cartwright in Section 1 of \cite{Ca}.
\end{rem}

\appendix
\section{Convex geometry}

In this appendix, we collect the notation used from convex geometry. We denote by $\Gamma$ a subgroup of $\rdop$. We consider a free abelian group $M$ of rank $n$ with dual $N:=\Hom(M,\zdop)$ and the corresponding real vector spaces $V:=M \otimes_\zdop \rdop$ and $W:=\Hom(V, \rdop)=N \otimes_\zdop \rdop$. The natural duality between $u \in V$ and $\omega \in W$ is denoted by $\langle u ,\omega \rangle$. References: \cite{Roc}, \cite{McM}.

\begin{art} \rm \label{polyhedron}
A {\it polyhedron} $\Delta$  in $W$ is an intersection of finitely many closed half-spaces $\{\omega \in W \mid \langle u_i , \omega \rangle \geq c_i\}$. We say that $\Delta$ is {\it $\Gamma$-rational} if we may choose all $u_i \in M$ and all $c_i \in \Gamma$. If $\Gamma=\qdop$, then we say that $\Delta$ is {\it rational}. A {\it closed face} of $\Delta$ is either $\Delta$ itself or has the form $H \cap \Delta$ where $H$ is the boundary of a closed half-space containing $\Delta$. An {\it open face} of $\Delta$ is a closed face  without all its properly contained closed faces. We denote 
by $\relint(\Delta)$ the unique open face of $\Delta$ which is dense in $\Delta$. 
\end{art}

\begin{art} \rm \label{polytope}
A bounded  polyhedron is called a {\it  polytope}. This is equivalent to being the convex hull of finitely many points. Let $G:=\{\lambda \in \rdop \mid \exists m \in \ndop \setminus \{0\}, \, m \lambda \in \Gamma\}$ be the divisible hull of $\Gamma$ in $\rdop$. Simple linear algebra shows that a  polytope is $\Gamma$-rational if and only if all vertices are $G$-rational and the edges have rational slopes. Similarly, a polyhedron is $\Gamma$-rational if and only if every closed face spans an affine subspace which is a translate of a rational linear subspace by a $G$-rational vector.
\end{art}


\begin{art} \rm \label{polyhedral complex}
A {\it polyhedral complex} $\Ccal$ in $W$ is a finite set of polyhedra such that 
\begin{itemize}
\item[(a)] $\Delta \in {\Ccal} \text{ $\Rightarrow$ all closed faces of $\Delta$ are in $\Ccal$}$;
\item[(b)] $\Delta, \sigma \in {\Ccal} \text{ $\Rightarrow$ $\Delta \cap \sigma$ is either empty or a closed face of $\Delta$ and $\sigma$}$.
\end{itemize}
The polyhedral complex is called {\it $\Gamma$-rational} if every $\Delta \in \Ccal$ is $\Gamma$-rational. The {\it support}  of $\Ccal$ is defined as $$|\Ccal|:= \bigcup_{\Delta \in
\Ccal} \Delta.$$
The polyhedral complex $\Ccal$ is called {\it complete} if $|\Ccal|=W$. A {\it subcomplex} of a polyhedral complex $\Ccal$ is a polyhedral complex $\Dcal \subset \Ccal$.
\end{art}

\begin{art} \rm \label{polyhedral decomposition}
A polyhedral complex $\Dcal$ {\it subdivides} the polyhedral complex $\Ccal$ if they have the same support and if every polyhedron $\Delta$ of $\Dcal$ is contained in a polyhedron of $\Ccal$. In this case, we say that $\Dcal$ is a {\it subdivision} of $\Ccal$.
\end{art}

\begin{art} \rm \label{cone}
A {\it cone} $\sigma$ in $W$ is centered at $0$, i.e. it is characterized by $\rdop_+ \sigma = \sigma$. Its {\it dual} is defined by 
$$\check{\sigma}:= \{u \in V \mid \langle u , \omega \rangle \geq 0 \, \; \forall \omega \in \sigma\}.$$ A {\it fan} is a polyhedral complex consisting of  polyhedral  cones. A {\it subfan} is a subcomplex of a fan.
\end{art}

\begin{art} \rm \label{local cone}
The {\it local cone} ${\rm LC}_\omega(S)$ of $S\subset W$ at $\omega$ is defined by 
$${\rm LC}_\omega(S):= \{ \omega' \in W \mid 
\text{$\omega + [0,\ve){\omega'} \subset S$ for some $\ve >0$}\}.
$$
\end{art}

\begin{art} \rm \label{recession cone}
The {\it recession cone} of a polyhedron $\Delta$ is defined by 
$$\rec(\Delta) := \{ \omega \in W \mid \omega + \Delta \subset \Delta\}.$$
By the Minkowski--Weil theorem, the recession cone is the unique  convex polyhedral cone $\sigma$ such that $\Delta = \sigma + \rho$ for a polytope $\rho$ of $W$. If $\Delta$ is $\Gamma$-rational, then $\rec(\Delta)$ is a rational convex polyhedral cone.
\end{art}

\begin{art} \rm \label{pointed polyhedron}
A polyhedron $\Delta$ is called {\it pointed} if it does not contain an affine line. Note that $\Delta$ is a pointed polyhedron if and only if $\rec(\Delta)$ has the origin $0$ as a vertex. This explains the terminology. A {\it pointed polyhedral complex} is a polyhedral complex consisting of pointed polyhedra.
\end{art}


\begin{art} \rm \label{proper convex function}
We say that $f: W \rightarrow \rdop \cup \{\infty\}$ is a {\it proper polyhedral function} if the {\it epigraph} ${\rm epi}(f):=\{(\omega,s) \in W \times \rdop \mid f(\omega) \leq s\}$ is a non-empty polyhedron. Then the faces of the polyhedron ${\rm epi}(f)$ contained in the graph of $f$ form a polyhedral complex in $W \times \rdop$ called the {\it graph complex}. The projection of the graph complex onto $W$ gives a polyhedral complex in $W$. Such a complex is called a {\it coherent polyhedral complex} in $W$.

Note that $f$ is a proper polyhedral function if and only if there is a non-empty polyhedron $\Sigma$ in $W$ and a  function $f_\Sigma:\Sigma \rightarrow \rdop$ with the following properties:

\begin{itemize}
 \item[(a)] $f_\Sigma$ is continuous and piecewise affine;
 \item[(b)] $f_\Sigma$ is a convex function in the usual sense, i.e. 
\begin{equation} \label{convex}
f_\Sigma(r\omega+s\omega') \leq r f_\Sigma(\omega) +s f_\Sigma(\omega')
\end{equation}
for $\omega, \omega' \in \Sigma$ and $r,s \in [0,1]$ with $r+s=1$.
 \item[(c)] $f$ agrees with $f_\Sigma$ on $\Sigma$ and $f=\infty$ outside of $\Sigma$.
\end{itemize}

We call $\Sigma$ the {\it domain of $f$}. The {\it domains of linearity} are the maximal subsets of $W$ where $f_\Sigma$ is affine. They are just the maximal dimensional polyhedra $\Delta$ from the coherent polyhedral complex corresponding to $f$. On such a $\Delta$, we have 
$$f(\omega)= c_\Delta + \langle u_\Delta, \omega \rangle$$
for some $u_\Delta \in V$ and $c_\Delta \in \rdop$. We call $u_\Delta$ the {\it peg} of $\Delta$.  
\end{art}

\begin{art} \label{conjugate} \rm 
Let $f$ be a proper polyhedral function with associated coherent polyhedral complex $\Ccal$. Then the {\it conjugate} of $f$ is the proper polyhedral function $f^*:V \rightarrow \rdop \cup \{\infty\}$ given by
$$f^*(u):= \sup\{\langle u, \omega \rangle - f(\omega)  \mid \omega \in W\}.$$
We have $f^{**}=f$. 
\end{art}

\begin{art} \rm \label{dual complex}
Let $f$ be a proper polyhedral function on $W$ with associated coherent polyhedral complex $\Ccal$. The coherent polyhedral complex in $V$ associated to $f^*$ is called the {\it dual complex} $\Ccal^f$ of $\Ccal$. The duality is a bijective order reversing correspondence $\sigma \mapsto \sigma^f$ between polyhedra of $\Ccal$ and polyhedra of $\Ccal^f$ given by
$$\sigma^f = \{ u \in V \mid f^
*(u)= \langle u, \omega \rangle - f(\omega) \; \forall \omega \in \sigma\}$$ 
and we have 
$$\dim(\sigma) + \dim(\sigma^f) = n.$$
This follows  from \cite{McM}, Theorem 7.1 and its proof. Note also that $\Ccal^f$ is complete if and only if the support of $\Ccal$ is bounded (which is then a polytope). 
\end{art}


\begin{art} \rm \label{halfspaces for totally ordered groups}
Let $G$ be a totally ordered abelian group which is not necessarily a subgroup of $\rdop$. Then $G$-rational polyhedra can be defined as subsets of $N_G:=N \otimes_\zdop G$ in the same way as in \ref{polyhedron}. For us, compatibilities are of interest in the situation when $G$ contains $\Gamma$ as an ordered subgroup where $\Gamma$ is our given subgroup of $\rdop$. Let $H^+:=\{\omega \in W \mid \langle u , \omega \rangle \geq c\}$ be a $\Gamma$-rational half-space in $W$ with $u \in M$ and $c \in \Gamma$. By base change, the homomorphism $N \rightarrow \zdop, \omega \mapsto \langle u, \omega \rangle$ extends canonically to a homomorphism $N_G \rightarrow G$ of abelian groups and hence we get a pairing $M \times N_G \rightarrow G$ which we also denote by $\langle \phantom{a}, \phantom{b} \rangle$. Then we set $H^+(G):=\{\omega  \in N_G \mid \langle u, \omega \rangle \geq c\}$. 
\end{art}

\begin{art} \label{polyhedra for totally ordered groups}
Let $H_i^+$ be a finite family of $\Gamma$-rational half-spaces in $W$ defining the $\Gamma$-rational polyhedron $\Delta := \bigcap_i H^+_i$. 
 Using the assumptions and the notation from \ref{halfspaces for totally ordered groups}, the intersection $\bigcap_i H^+_i(G)$ depends only on $\Delta$ and not on the choice of the half-spaces.
\end{art}

\proof The given half-spaces have the form $H_i^+:=\{\omega \in W \mid \langle u_i , \omega \rangle \geq c_i\}$ with $u_i \in M$ and $c_i \in \Gamma$. It is enough to show that $H^+(G) \supset \bigcap_i H^+_i(G)$ for any $\Gamma$-rational half-space $H^+ := \{\omega \in W \mid \langle u , \omega \rangle \geq c\} \supset \Delta$. We may assume that $H^+$  is a supporting half-space of $\Delta$, i.e. there is $\omega_0 \in W$ such that $\langle u, \omega_0 \rangle =c$. Moreover, we may assume that all $H_i^+$ are supporting hyperplanes in $\omega_0$, i.e. $\langle u_i, \omega_0 \rangle =c_i$. We note that the recession cone $ \rho=\{\omega \in W \mid \langle u, \omega\rangle \geq 0 \}$ of $H^+$ contains the recession cone $\sigma = \{ \omega \in W \mid \langle u_i, \omega \rangle \geq 0 \; \forall i\}$ of $\Delta$. By the bijective correspondence between  rational polyhedral cones in $W=N_\rdop$ and finitely generated saturated semigroups of $M$ (see \ref{affine toric variety}), we see that $\check{\rho} \cap M \subset \check{\sigma} \cap M = \{ \sum_i m_i u_i \in M \mid 0 \leq m_i \in \qdop \, \forall i\}$. We conclude that $u=\sum_i m_i u_i$ with $ 0 \leq m_i \in \qdop$ for all $i$. Replacing $u$ by a positive multiple, we may assume that all $m_i \in \ndop$. For $\omega \in \bigcap_i H^+_i(G)$, we get 
$$\langle u, \omega \rangle = \sum_i m_i \langle u_i, \omega \rangle \geq \sum_i m_i c_i=\sum_i m_i\langle u_i, \omega_0 \rangle =\langle u, \omega_0 \rangle =c$$
and hence $\omega \in H^+(G)$. This proves the claim. \qed

\begin{art} \label{compatibility properties} \rm
Let $\Delta$ be a $\Gamma$-rational polyhedron in $W$ as above. Then \ref{polyhedra for totally ordered groups} shows that $\Delta(G):= \bigcap_i H_i^+(G)$ is well-defined. For another $\Gamma$-rational polyhedron $\Delta'$ in $W$, we get $\Delta(G) \cap \Delta'(G)=(\Delta \cap \Delta')(G)$. In particular, if $\Delta' \subset \Delta$, then $\Delta'(G) \subset \Delta(G)$. 

More generally, if finitely many $\Gamma$-rational polyhedra $\Delta_i$ cover $\Delta$, then $\Delta(G) \subset \bigcup_i \Delta_i(G)$. To see this, let $(H_j^+)_{j=1, \dots, r}$ be the half-spaces occurring in the definitions of the polyhedra. We denote by $H_j^-$  the half-space on the other side of the boundary of $H_j^+$. Since $G$ is a totally ordered group, the sets $\Delta(G)\cap H_1^\pm(G) \cap \dots \cap H_r^\pm(G)$ cover $\Delta(G)$. Every polyhedron $\Delta \cap H_1^\pm \cap \dots \cap H_r^\pm$ is contained in a $\Delta_i$ and we conclude from \ref{polyhedra for totally ordered groups} and the above that the sets  $\Delta_i(G)$ cover $\Delta(G)$.
\end{art}

\begin{art} \label{compatibility of affine maps} \rm
Let us still consider a totally ordered abelian group $G$ containing $\Gamma$ as an ordered subgroup and a $\Gamma$-rational polyhedron $\Delta$ in $W=N_\rdop$. 
If $\varphi:N' \rightarrow N$ is a homomorphism of free abelian groups of finite rank and if $\psi:N'_\Gamma \rightarrow  N_\Gamma$ is given by $\psi = \varphi_\Gamma + \omega_0$ for a fixed $\omega_0 \in N_\Gamma$, then we have $$(\psi^{-1}(\Delta))(G)=\psi_G^{-1}(\Delta(G)),$$
 where $\varphi_\Gamma$ and $\psi_G$ denote base changes of $\varphi$ and $\psi$. This is easily seen  for a $\Gamma$-rational half-space and then the claim follows from  \ref{polyhedra for totally ordered groups}. 
\end{art}

\bibliographystyle{amsalpha}

{\small Walter Gubler, Fakult\"at f\"ur Mathematik,  Universit\"at Regensburg,
Universit\"atsstrasse 31, D-93040 Regensburg, walter.gubler@mathematik.uni-regensburg.de}
\end{document}